\documentclass[a4paper, 11pt]{scrartcl}
\pdfoutput=1

\title{Generalised Summation-by-Parts Operators and Variable Coefficients}
\author{Hendrik Ranocha}
\date{14th November 2017}
\titlehead{\includegraphics[width=0.5\textwidth]{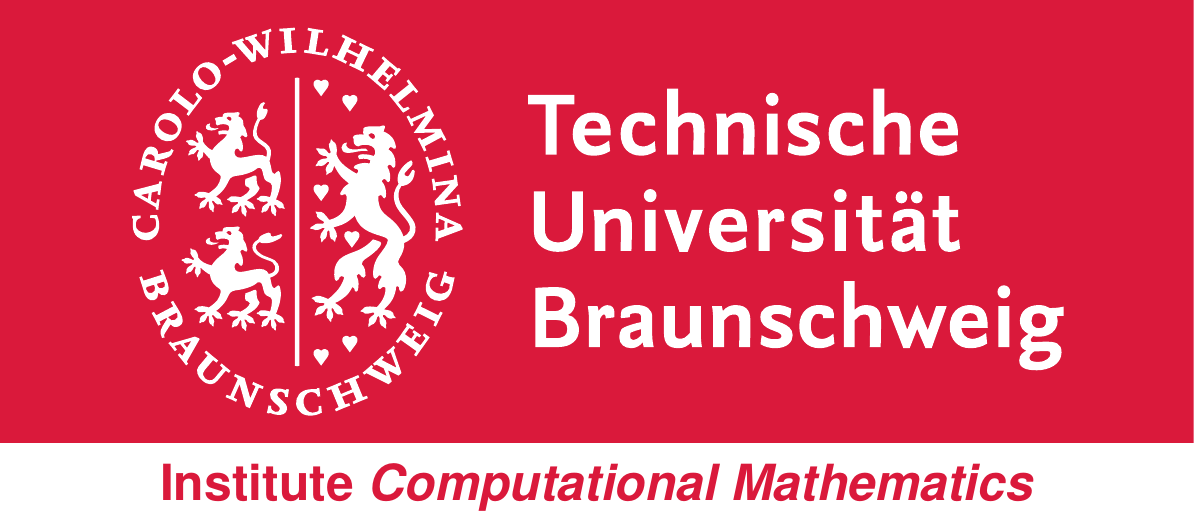}}

\usepackage[utf8]{luainputenc}
\usepackage[english]{babel}
\usepackage{csquotes}

\usepackage[a4paper, top=0.51in, bottom=0.79in, left=0.999in, right=0.99in]{geometry}

\usepackage[plainpages=false,pdfpagelabels,hidelinks,unicode]{hyperref}
\makeatletter
\hypersetup{pdfauthor={\@author}}
\hypersetup{pdftitle={\@title}}
\makeatother

\usepackage[%
  backend=biber,
  style=numeric-comp,
  firstinits=true, uniquename=init, 
  natbib=true,
  url=true,
  doi=true,
  isbn=false,
  backref=false,
  maxnames=99,
  ]{biblatex}
\usepackage{amsmath}
\allowdisplaybreaks
\usepackage{amssymb}
\usepackage{commath}
\usepackage{mathtools}
\usepackage{bbm}

\usepackage{siunitx}

\usepackage{amsthm}
\theoremstyle{plain}
  \newtheorem{theorem}{Theorem}
  
\theoremstyle{definition}
  \newtheorem{definition}[theorem]{Definition}
  \newtheorem{remark}[theorem]{Remark}
  
\usepackage{color}
\usepackage{graphicx}
\usepackage[small]{caption}
\usepackage{subcaption}

\usepackage{pgfplots}
\pgfplotsset{compat=1.11}

\begingroup\expandafter\expandafter\expandafter\endgroup
\expandafter\ifx\csname pdfsuppresswarningpagegroup\endcsname\relax
\else
\fi

\usepackage{booktabs}
\usepackage{rotating}
\usepackage{multirow}

\usepackage{multicol}
\usepackage{enumitem}

\usepackage{calc}
\usepackage{xparse}

\renewcommand{\vec}[1]{\underline{#1}}
\NewDocumentCommand{\mat}{mo}{%
  \IfValueTF{#2}{%
    \underline{\underline{#1}}{#2}
  }{%
    \underline{\underline{#1}}\,
  }%
}
\newcommand{\diag}[1]{\operatorname{diag}\left(#1\right)}

\newcommand{\scp}[2]{\left\langle{#1,\, #2}\right\rangle}

\newcommand{\inv}[1]{{#1}^{-1}}
\newcommand{\vect}[1]{\begin{pmatrix} #1 \end{pmatrix}}

\newcommand{\xmin}{{x_L}}
\newcommand{\xmax}{{x_R}}
\newcommand{\e}{\mathrm{e}}

\newcommand{\fnum}{f^{\mathrm{num}}}

\newcommand{\vecfnum}{\vec{f}^{\mathrm{num}}}

\renewcommand{\epsilon}{\varepsilon}
\renewcommand{\phi}{\varphi}
\renewcommand{\rho}{\varrho}

\newcommand{\R}{\mathbb{R}}

\newsavebox{\DelimiterBox}
\newlength{\DelimiterHeight}
\newlength{\DelimiterDepth}
\newsavebox{\ArgumentBox}
\newlength{\ArgumentHeight}
\newlength{\ArgumentDepth}
\newlength{\ResizedDelimiterHeight}
\newlength{\ResizedDelimiterDepth}

\begin{document}

\maketitle

\begin{abstract}
  
High-order methods for conservation laws can be highly efficient if their
stability is ensured. A suitable means mimicking estimates of the continuous
level is provided by summation-by-parts (SBP) operators and the weak enforcement
of boundary conditions. Recently, there has been an increasing interest in
generalised SBP operators both in the finite difference and the discontinuous
Galerkin spectral element framework.
However, if generalised SBP operators are used, the treatment of the boundaries
becomes more difficult since some properties of the continuous level are no longer
mimicked discretely --- interpolating the product of two functions will in general
result in a value different from the product of the interpolations. Thus, desired
properties such as conservation and stability are more difficult to obtain.
Here, new formulations are proposed, allowing the creation of discretisations
using general SBP operators that are both conservative and stable. Thus, several
shortcomings that might be attributed to generalised SBP operators are overcome
(cf. J.~Nordström and A.~A.~Ruggiu, ``On Conservation and Stability Properties
for Summation-By-Parts Schemes'', \emph{Journal of Computational Physics} 344
(2017), pp. 451--464, and J.~Manzanero, G.~Rubio, E.~Ferrer, E.~Valero and
D.~A.~Kopriva, ``Insights on aliasing driven instabilities for advection
equations with application to Gauss-Lobatto discontinuous Galerkin methods'',
Journal of Scientific Computing (2017), \url{https://doi.org/10.1007/s10915-017-0585-6}).

\end{abstract}

\section{Introduction}
\label{sec:introduction}

Considering the solution of hyperbolic partial differential equations (PDEs), high order methods
can be very efficient, providing accurate numerical solutions with relatively
low computational effort \cite{kreiss1972comparison}. In order to make use of
this accuracy, stability has to be established. Mimicking estimates on the
continuous level for the PDE via integration-by-parts, summation-by-parts (SBP)
operators \cite{kreiss1993stability, strand1994summation} can be used. 
In short, SBP operators are discrete derivative operators equipped with a
compatible quadrature providing a discrete analogue of the $L^2$ norm. The
compatibility of discrete integration and differentiation mimics integration-by-parts
on a discrete level.
Combined with the weak enforcement of boundary conditions via simultaneous approximation
terms (SATs) \cite{carpenter1993stability}, highly efficient and stable
semidiscretisations can be obtained, as described also in \cite{fernandez2014review,
svard2014review, nordstrom2015new, gustafsson2013time} and references cited
therein.

During the last years, there has been an enduring and increasing interest in the
basic ideas of SBP operators and their application in various frameworks including 
finite volume (FV) \cite{nordstrom2001finite, nordstrom2003finite}, 
discontinuous Galerkin (DG) \cite{gassner2013skew, gassner2014kinetic, kopriva2014energy,
gassner2016well, ortleb2016kinetic},
and the recent flux reconstruction / correction procedure via reconstruction
framework \cite{huynh2007flux, wang2009unifying, vincent2011newclass, huynh2014high,
witherden2014pyfr, vincent2015extended} as described in \cite{ranocha2016summation}.

While a mesh of equidistant nodes including the boundaries of the domain is 
classically used in finite difference methods, multiple elements and nonuniform
point distributions possibly not including the boundary nodes are used in
discontinuous Galerkin and flux reconstruction schemes. However, the discrete
differentiation and quadrature rules can still be compatible, mimicking 
integration-by-parts.
In this context, the notion of generalised summation-by-parts operators 
\cite{fernandez2014generalized} covers many stable semidiscretisations and
also operators on unstructured meshes in multiple space dimensions
\cite{hicken2016multidimensional}. Moreover, using the framework of entropy
stable numerical fluxes of \citet{tadmor1987numerical, tadmor2003entropy},
entropy stable semidiscretisations can be constructed for systems of nonlinear
conservation laws \cite{fisher2013highJCP, carpenter2014entropy, parsani2015entropy,
yamaleev2017family, gassner2016split, wintermeyer2017entropy, ranocha2017shallow,
ranocha2017comparison}.
These methods can be interpreted as generalisations of split forms of the
underlying PDE. For classical SBP operators in finite difference methods,
these split forms are well known to have desirable conservation and stability
properties \cite{nordstrom2006conservative, fisher2013discretely}.

However, if generalised SBP operators are used, several problems have to be
handled for non-diagonal norms \cite{svard2004coordinate} and nodal bases not
including boundary nodes, as can be seen in recent investigations by
\citet{nordstrom2017conservation} and \citet{manzanero2017insights}.
In order to remedy these problems, corrections for generalised SBP operators
will be presented in this article, extending the work of \cite{ranocha2016sbp,
ranocha2017extended}.

Therefore, the article is structured as follows. In \autoref{sec:SBP},
the concept of summation-by-parts operators and generalised SBP operators is
introduced. 
Afterwards, a conservation law with variable coefficients, introducing aliasing
issues, is studied in \autoref{sec:standard-L2} using the $L^2$ framework of
\citet{nordstrom2017conservation}. Introducing new corrections, generalised SBP
operators are proven to be both conservative and stable. 
Similarly, the weighted $L^2$ framework of \citet{manzanero2017insights} is applied 
in \autoref{sec:weighted-L2}, yielding new conservative and stable discretisations
using general SBP operators. Thereafter, numerical results for the new methods
are presented in \autoref{sec:numerical-results}, including convergence studies
and eigenvalue analyses. Finally, some conclusions are presented in \autoref{sec:summary}.
Moreover, since the notations used in the finite difference and finite element
framework are sufficiently different, some translation rules and comparisons with
previous works are provided in the appendix.

\section{Summation-by-Parts Operators}
\label{sec:SBP}

In this section, a general notion of summation-by-parts (SBP) operators for 
semidiscretisations of conservation laws will be presented using the notation of
\cite{ranocha2016summation, ranocha2017extended}. Tables containing translations 
to the notations used in the finite difference community \citep{nordstrom2017conservation}
and the finite element framework \citep{manzanero2017insights} can be found in
\ref{sec:translation-rules}.

\subsection{General Setting}
\label{subsec:general-setting}

In order to discretise a scalar conservation law in one dimension
\begin{equation}
\label{eq:CL}
\begin{aligned}
  \partial_t u(t,x) + \partial_x f\big(t,x,u(t,x)\big) &= 0,
  \qquad t > 0, \; x \in (\xmin, \xmax),
\end{aligned}
\end{equation}
equipped with appropriate initial and boundary conditions, the domain 
$(\xmin, \xmax)$ is partitioned into several intervals
$(x_0, x_1), (x_1, x_2), \dots, (x_{N-1}, x_N)$, where
$\xmin = x_0 < x_1 < \dots < x_N = \xmax$.
On each element $(x_{i-1}, x_i)$, the numerical solution is a vector,
represented in some finite-dimensional basis by $\vec{u} = (u_0, \dots, u_p)^T$.
In finite difference (FD) and nodal discontinuous Galerkin (DG) methods,
nodal bases are used, i.e. the representation $\vec{u}$ is given by the values
of $u$ at pairwise different points $\xi_0, \dots, \xi_p$. Other choices are
also possible, e.g. modal bases \citep{ranocha2016sbp, ranocha2017extended}.

With respect to the chosen basis, (an approximation of) the derivative is
represented by the matrix $\mat{D}$. Moreover, a discrete scalar product is
represented by the symmetric and positive definite mass / norm matrix $\mat{M}$,
approximating the $L^2$ scalar product, i.e.
\begin{equation}
  \mat{D} \vec{u} \approx \vec{\partial_x u},
  \qquad
  \vec{u}^T \mat{M} \vec{v} = \scp{\vec{u}}{\vec{v}}_{\mat{M}}
  \approx
  \scp{u}{v}_{L^2} = \int_{x_{i-1}}^{x_i} u \, v.
\end{equation}
Additionally, a restriction operator is represented by the matrix $\mat{R}$,
approximating the interpolation of a function to the boundary points
$\set{x_{i-1}, x_i}$. Furthermore, the diagonal boundary matrix
$\mat{B} = \diag{-1,1}$ gives the difference of boundary values, i.e.
\begin{equation}
\label{eq:R-and-B}
  \mat{R} \vec{u} \approx \vect{u(x_{i-1}) \\ u(x_i)},
  \qquad
  (u_L, u_R) \cdot \mat{B} \cdot \vect{v_L \\ v_R} = u_R v_R - u_L v_L.
\end{equation}
Finally, these operators fulfil the summation-by-parts property
\begin{equation}
\label{eq:SBP}
  \mat{M} \mat{D} + \mat{D}[^T] \mat{M} = \mat{R}[^T] \mat{B} \mat{R},
\end{equation}
mimicking integration-by-parts via
\begin{equation}
\label{eq:SBP-IBP}
\begin{array}{*3{>{\displaystyle}c}}
  \underbrace{
    \vec{u}^T \mat{M} \mat{D} \vec{v} 
    + \vec{u}^T \mat{D}[^T] \mat{M} \vec{v}
  }
  & = &
  \underbrace{
    \vec{u}^T \mat{R}[^T] \mat{B} \mat{R} \vec{v},
  }
  \\
  \rotatebox{90}{$\!\approx\;$}
  &&
  \rotatebox{90}{$\!\!\approx\;$}
  \\
  \overbrace{
    \int_{x_{i-1}}^{x_i} u \, (\partial_x v) 
    + \int_{x_{i-1}}^{x_i} (\partial_x u) \, v
  }
  & = &
  \overbrace{
    u \, v \big|_{x_{i-1}}^{x_i}
  }.
\end{array}
\end{equation}
If linear equations with constant coefficients are considered, these operators
suffice to describe appropriate semidiscretisations. For varying coefficients or
nonlinear equations, additional multiplication operators have to be used.
If the function $u$ is represented by $\vec{u}$, the discrete multiplication
operator approximating the linear operator $v \mapsto u v$ is represented by
the matrix $\mat{u}$, mapping $\vec{v}$ to $\mat{u} \vec{v}$.

In a nodal basis, the standard multiplication operators consider pointwise
multiplication, i.e. they are diagonal $\mat{u} = \diag{\vec{u}} = \diag{u_0, \dots, u_p}$.
In modal bases, exact multiplication followed by an $L^2$ projection can be
considered \citep{ranocha2016sbp, ranocha2017extended}.

\subsection{Analytical Setting}
\label{subsec:analytical-setting}

In nodal DG methods, the numerical solution is represented as a polynomial of
degree $\leq p$ on each element. Thus, the matrices $\mat{M}, \mat{D}, \mat{R}$
can be computed such that they represent the corresponding operations exactly
and the SBP property \eqref{eq:SBP} will be fulfilled since it is exactly given
by integration-by-parts.

However, the exactness of the mass matrix $\mat{M}$ can be relaxed, since
only products of one polynomial of degree $\leq p$ and another polynomial of
degree $\leq p-1$ are integrated in \eqref{eq:SBP-IBP}.
Considering nodal bases with diagonal mass matrices
$\mat{M} = \diag{\omega_0, \dots, \omega_p}$, a corresponding quadrature rule
is given by the positive weights $\omega_0, \dots, \omega_p$. Thus, the SBP
property \eqref{eq:SBP} is fulfilled, if the quadrature given by $\mat{M}$ is
exact for polynomials of degree $\leq 2p-1$, i.e. Gauss-Lobatto, Gauss-Radau,
or Gauss quadrature, since these are exact for polynomials of degree
$\leq 2p-1$, $\leq 2p$, and $\leq 2p+1$, respectively \citep{kopriva2010quadrature,
hicken2013summation}.

\subsection{Numerical Setting}
\label{subsec:numerical-setting}

There are classical finite difference operators that can be interpreted in an
analytical setting similar to the one described above \citep{gassner2013skew}.
However, in general --- to the authors knowledge --- it is not known whether
finite difference SBP operators correspond to an analytical basis. Nevertheless,
the basic requirement of the SBP property \eqref{eq:SBP} can be enhanced by
accuracy conditions in order to get a useful definition of (numerical) SBP
operators, cf. \cite{fernandez2014generalized, nordstrom2017conservation}.

\begin{definition}
  Using a numerical representation $\vec{u} = (u_0, \dots, u_p)^T$ of a function
  $u$ in a nodal basis, the operators $\mat{D}, \mat{R}, \mat{M}$, and $\mat{B}$
  described above form a $q$th order SBP discretisation, if
  \begin{enumerate}
    \item 
    the derivative matrix $\mat{D}$ is exact for polynomials of degree $q$,
    
    \item 
    the mass / norm matrix $\mat{M}$ is symmetric and positive definite,
    
    \item 
    the SBP property \eqref{eq:SBP} is fulfilled, and
    
    \item 
    $\vec{u}^T \mat{R}[^T] \mat{B} \mat{R} \vec{v} = u \, v \big|_{x_{i-1}}^{x_i}$
    is exact for polynomials $u,v$ with degrees
    $\operatorname{deg}(u), \operatorname{deg}(v) \leq r$, where $r \geq q$.
  \end{enumerate}
\end{definition}

These kind of SBP semidiscretisations together with diagonal multiplication operators
$\mat{u} = \diag{\vec{u}}$ will be considered in the following. However, modal
bases can be used as well \citep{ranocha2016sbp, ranocha2017extended}.

\subsection{Special Classes of SBP Discretisations}
\label{subsec:special-SBP}

There are several properties of some SBP discretisations that should be considered
separately since they yield special simplifications.

\subsubsection{Diagonal Norms Versus Non-Diagonal Norms}
\label{sec:diag-vs-nondiag-norms}

If the nodal basis is equipped with a diagonal mass matrix
$\mat{M} = \diag{\omega_0, \dots, \omega_p}$, the positive weights $\omega_i$
correspond to a quadrature rule, cf. \cite{hicken2013summation}.
Moreover, the multiplication operators commute with the mass matrix, since both
are diagonal. Thus, multiplication operators $\mat{u}$ are self-adjoint with
respect to the scalar product induced by $\mat{M}$, similarly to multiplication
operators in the continuous setting (if the correct domain is chosen).
This exact mimicking of properties at the discrete level is desirable, allowing
proofs of semidiscrete results analogously to the continuous counterparts.

However, corrections are sometimes possible, if multiplication operators are
not self-adjoint. In this case, the $\mat{M}$-adjoint operator
\begin{equation}
\label{eq:M-adjoint}
  \mat{u}[^*] = \mat{M}[^{-1}] \mat{u}[^T] \mat{M}
\end{equation}
can be used to construct appropriate corrections \cite{ranocha2016sbp,
ranocha2017extended, ranocha2017shallow}. This has to be used since the discrete
multiplication operators are not exact.

\subsubsection{Inclusion of Boundary Nodes Versus Bases Exclusion of Boundary Nodes}
\label{sec:including-boundary-nodes-vs-not-including-boundary-nodes}

Another useful property of SBP discretisations is the inclusion of the boundary
points $\set{x_{i-1}, x_i}$ into the nodes of the basis. In this case, the restriction
operators are simply
\begin{equation}
  \mat{R}
  =
  \begin{pmatrix}
    1 & 0 & \dots & 0 & 0 \\
    0 & 0 & \dots & 0 & 1
  \end{pmatrix},
  \qquad
  \mat{R}[^T] \mat{B} \mat{R}
  =
  \diag{-1, 0, \dots, 0, 1}.
\end{equation}
Thus, restriction to the boundary and multiplication commute, i.e.
\begin{equation}
  \left( \mat{R} \vec{u} \right) \cdot \left( \mat{R} \vec{v} \right)
  =
  \vect{ u_0 v_0 \\ u_p v_p }
  =
  \mat{R} \mat{u} \vec{v} = \mat{R} \mat{v} \vec{u}.
\end{equation}
This property is fulfilled at the continuous level and therefore also
desirable for the semidiscretisation.
If the boundary nodes are not included, restriction and multiplication will in
general not commute. In this case, some corrections can be applied, e.g. the
application of a linear combination of
$\left( \mat{R} \vec{u} \right) \cdot \left( \mat{R} \vec{v} \right)$
and $\mat{R} \mat{u} \vec{v}$ in order to get the desired results
\cite{ranocha2016sbp, ranocha2017extended, ranocha2017shallow}.
Again, the reason for this is the inexactness of discrete multiplication operators.
However, the construction of correction terms can become cumbersome. In other
situations, it might be unknown whether suitable correction terms can be found
\cite{ranocha2017comparison}.

\subsection{Simultaneous Approximation Terms and Numerical Fluxes}
\label{subsec:SQTs-and-fnum}

In order to enforce boundary conditions either between neighbouring elements or
at the physical boundary of the domain, several methods can be applied.
For finite difference methods, these include the injection method 
(combining the boundary and differential operators directly), the projection
method and the usage of simultaneous approximation terms (SATs) 
\citep{olsson1995summationI, olsson1995summationII, mattsson2003boundary}.
In discontinuous Galerkin methods, boundary conditions are usually enforced via
numerical fluxes, resulting in stable schemes. For linear problems, the weak
enforcement of boundary conditions via numerical fluxes is equivalent to the use
of SATs, as will be shown for the examples in the following sections.

For a conservation law $\partial_t u(t,x) + \partial_x f\big( u(t,x) \big) = 0$
with flux $f$ not depending explicitly on space and time, a (two-point) numerical
flux $\fnum = \fnum(u_-, u_+)$ will be considered to be a Lipschitz continuous
mapping consistent with the continuous flux $f$ of the conservation law, i.e. 
$\fnum(u, u) = f(u)$, as common in conservative finite difference and finite
volume methods, see e.g. \cite{leveque1992numerical, leveque2002finite}.

\begin{figure}
  \centering
  \begin{tikzpicture}
  \begin{axis}[
      width=10cm, height=5cm,
      xlabel=$x$, ylabel=$u$, ticks=none, xlabel near ticks, ylabel near ticks,
      axis lines=left, axis line style={->},
      xmin=-1.1, xmax=1.1, ymin=-2.5, ymax=4.0,
    ]
    \addplot[solid, samples=500, domain=-1:0] {x^3 - 2*x^2 - x};
    \addplot[solid, samples=500, domain=0:1] {x^3 + 4*x^2 - x - 0.5};
    
    \node[inner sep=0,label={right:$u^{(i-1)}_L$}] at (axis cs:-1,-2) {$\bullet$};
    \node[inner sep=0,label={above left:$u^{(i-1)}_R = u_-$}] at (axis cs:0,0) {$\bullet$};
    \node[inner sep=0,label={below right:$u^{(i)}_L = u_+$}] at (axis cs:0,-0.5) {$\bullet$};
    \node[inner sep=0,label={left:$u^{(i)}_R$}] at (axis cs:1,3.5) {$\bullet$};
    
    \draw[dotted] (axis cs:-1,-2.5) -- (axis cs:-1,4);
    \draw[dotted] (axis cs:0,-2.5) -- (axis cs:0,4);
    \draw[dotted] (axis cs:1,-2.5) -- (axis cs:1,4);
    
    \node at (axis cs:-0.5,3.5) {Cell $i-1$};
    \node at (axis cs:0.5,3.5) {Cell $i$};
  \end{axis}
  \end{tikzpicture}
  \caption{Visualisation of the notation for numerical fluxes $\fnum(u_-,u_+)$
           between two cells.}
  \label{fig:notation-numerical-fluxes}
\end{figure}
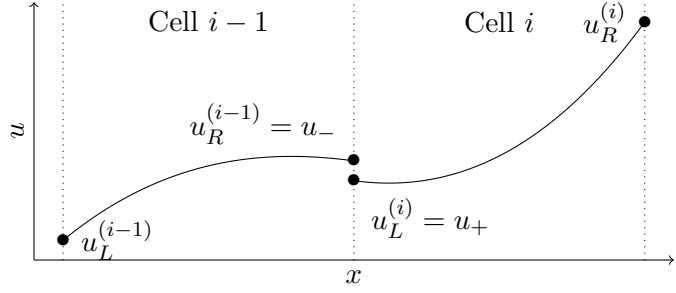

In methods using numerical fluxes $\fnum$, these fluxes are evaluated
at the boundaries between two cells. In order to introduce the notation, consider
two one-dimensional cells $i-1$ and $i$. On each cell~$j$, the numerical solution 
has a value at the left and right hand side of this cell, written as $u^{(j)}_L$ and
$u^{(j)}_R$, respectively. At the boundary between the cells $i$ and $i+1$, there
are two values of the numerical solutions from the neighbouring cells, i.e.
$u^{(i-1)}_R =: u_-$ and $u^{(i)}_L =: u_+$. Here, the indices $-$ and $+$ refer to
values computed from the left and right of a given boundary, respectively.
Thus, the numerical flux between the cells $i-1$ and $i$ is evaluated as
$\fnum = \fnum(u_-, u_+) = \fnum\Big( u^{(i-1)}_R, u^{(i)}_L \Big)$, as visualised in
\autoref{fig:notation-numerical-fluxes}.
Finally, the two numerical fluxes $\fnum_L$ and $\fnum_R$ at the left and right 
hand side of a cell are written as components of the vector 
$\vecfnum = \left( \fnum_L, \fnum_R \right)^T$.

In the following sections, the linear advection equation with variable coefficients
\begin{equation}
\label{eq:var-lin-adv}
\begin{aligned}
  \partial_t u(t,x) + \partial_x \left( a(x) u(t,x) \right) &= 0,
  \qquad&& t > 0, \; x \in (\xmin, \xmax),
  \\
  u(t,\xmin) &= g_L(t),
  \qquad&& t \geq 0,
  \\
  u(0,x) &= u_0(x),
  \qquad&& x \in [\xmin,\xmax],
\end{aligned}
\end{equation}
with variable speed $a(x) > 0$ and compatible initial and boundary conditions
$u_0$, $g_L$ will be considered. Since the flux $f$ depends explicitly on the 
space coordinate $x$ via the variable coefficients $a(x)$, the notion of numerical 
fluxes has to be adapted correspondingly. There are at least two possibilities 
concerning the evaluation of the variable coefficients:
\begin{enumerate}[label=\roman*)]
  \item \label{itm:use-edge-fnum}
  Use the speed $a(x)$ evaluated at the position of the boundary where the
  numerical flux is evaluated.
  
  \item \label{itm:use-interpolated-fnum}
  Use some interpolation of the discretised function $\vec{a}$ in the numerical
  flux.
\end{enumerate}
If boundary nodes are included and the coefficients $\vec{a}$ of the discrete
version of the function $a$ are obtained by evaluating $a$ at these nodes, both
possibilities are identical. However, if boundary nodes are not included, the
two approaches are different. Moreover, there are several possible versions of
the interpolations that can be used. In the remaining sections, the following
numerical fluxes for \eqref{eq:var-lin-adv} will be considered:
\begin{align}
  \label{eq:fnum-edge-central}
  \text{Edge based central flux}\qquad
  \fnum(u_-,u_+) &= a(x) \frac{u_- + u_+}{2},
  \\
  \label{eq:fnum-split-central}
  \text{Split central flux}\qquad
  \fnum(u_-,u_+) &= \frac{a_- u_- + a_+ u_+}{2},
  \\
  \label{eq:fnum-unsplit-central}
  \text{Unsplit central flux}\qquad
  \fnum(u_-, u_+) &= \frac{ (a u)_- + (a u)_+ }{2},
  \\
  \label{eq:fnum-edge-upwind}
  \text{Edge based upwind flux}\qquad
  \fnum(u_-,u_+) &= a(x) \, u_-,
  \\
  \label{eq:fnum-split-upwind}
  \text{Split upwind flux}\qquad
  \fnum(u_-,u_+) &= a_- u_-,
  \\
  \label{eq:fnum-unsplit-upwind}
  \text{Unsplit upwind flux}\qquad
  \fnum(u_-, u_+) &= (a u)_-.
\end{align}
Here, $a_- u_-$ is the product of the interpolated values of $a$ and $u$ from
the cell on the left hand side. Similarly, $(a u)_-$ is the interpolated value 
of the product $a u$. In the edge based fluxes, the velocity $a(x)$ is evaluated 
at the position of the boundary.
In the edge based numerical fluxes \eqref{eq:fnum-edge-central} and 
\eqref{eq:fnum-edge-upwind}, the possibility \ref{itm:use-edge-fnum} from above
is used. The split and unsplit fluxes use possibility \ref{itm:use-interpolated-fnum}
with two different choices of the interpolation.

In order to help the reader avoid misunderstandings on the notation, some
translation rules to the finite difference setting of \cite{nordstrom2017conservation}
and the discontinuous Galerkin spectral element framework of \cite{kopriva2014energy, 
manzanero2017insights} are presented in \ref{sec:translation-rules}.

\section{Standard \texorpdfstring{$L^2$}{L²} Estimates}
\label{sec:standard-L2}

In this section, the standard (i.e. unweighted) $L^2$ estimates of 
\citet{nordstrom2017conservation}
will be reviewed and suitable formulations for general SBP discretisations will
be presented, resulting in both conservative and stable methods.

\subsection{Continuous Setting}
\label{subsec:standard-continuous}

Consider the linear advection equation \eqref{eq:var-lin-adv}
with variable speed $a(x) > 0$ and compatible initial and boundary conditions
$u_0$, $g_L$. The desired properties of this model problem are conservation and 
stability \citep[Section 2]{nordstrom2017conservation}.
In order to investigate conservation, the rate of change of the total mass of $u$
can be expressed as
\begin{equation}
\label{eq:standard-conservation}
  \od{}{t} \int_\xmin^\xmax u
  =
  \int_\xmin^\xmax \partial_t u
  =
  - \int_\xmin^\xmax \partial_x (a \, u)
  =
  - a \, u \big|_\xmin^\xmax
  =
  a(\xmin) g_L - a(\xmax) u(\xmax).
\end{equation}
Thus, the change of the total mass is given by the flux through the boundaries.

Similarly, to study stability, the standard $L^2$ estimate can be obtained by
application of the product rule and integration by parts
\begin{equation}
\label{eq:standard-stability}
\begin{aligned}
  \od{}{t} \int_\xmin^\xmax u^2
  &=
  2 \int_\xmin^\xmax u \, \partial_t u
  =
  -2 \int_\xmin^\xmax u \, \partial_x (a \, u)
  =
  - \int_\xmin^\xmax \left(
    u \, \partial_x (a \, u) 
    + a \, u \, \partial_x u
    + u^2 \, \partial_x a
  \right)  
  \\
  &=
  - a \, u^2 \big|_\xmin^\xmax
  - \int_\xmin^\xmax u^2 \, \partial_x a
  =
  a(\xmin) g_L^2 - a(\xmax) u(\xmax)^2 - \int_\xmin^\xmax u^2 \, \partial_x a.
\end{aligned}
\end{equation}
Thus, if $\partial_x a$ is bounded, the energy estimate
\begin{equation}
\label{eq:standard-continuous-energy-estimate}
  \norm{ u(t) }_{L^2}^2
  \leq
  \e^{t \, \norm{\partial_x a}_{L^\infty}}
  \left(
    \norm{ u_0 }_{L^2}^2
    + \int_0^t \e^{-\tau \, \norm{\partial_x a}_{L^\infty}}
      \left(
        a(\xmin) g_L(\tau)^2 - a(\xmax) u(\tau,\xmax)^2
      \right) \dif \tau
  \right)
\end{equation}
follows, cf. \cite[Section 2]{nordstrom2017conservation}. Hence, in the rest of
this section, $\partial_x a \in L^\infty$ will be assumed.

\begin{remark}
\label{rem:positive-speed}
  In order to enable the above computations, it suffices to consider a speed $a(x)$
  that is positive at the boundaries $x_L, x_R$. If the speed is positive everywhere,
  similar properties hold on subintervals of $(x_L, x_R)$ --- corresponding to
  multi-block or multi-element discretisations.
  
  In order to simplify the following investigations, $a(x)$ is assumed to be
  positive everywhere. The case of negative speed or speed with varying sign can
  be handled analogously. If $a$ is negative at the boundaries, a boundary condition
  has to be prescribed only at $x_R$ instead of $x_L$. Similarly, the upwind
  numerical fluxes \eqref{eq:fnum-edge-upwind}, \eqref{eq:fnum-split-upwind}, and
  \eqref{eq:fnum-unsplit-upwind} use the values from the right-hand side (with
  index $+$ instead of $-$) for negative speed. The central numerical fluxes
  \eqref{eq:fnum-edge-central}, \eqref{eq:fnum-split-central}, and
  \eqref{eq:fnum-unsplit-central} remain unchanged.
\end{remark}

\subsection{Diagonal Norm Discretisations Including Boundary Nodes}
\label{subsec:standard-simplified}

In this section, diagonal norm SBP discretisations including the boundary nodes
are considered. In this case, multiplication operators are self-adjoint and restriction
to the boundary and multiplication commute.
Moreover, as mentioned in section~\ref{subsec:SQTs-and-fnum}, the numerical fluxes
(edge based, split, and unsplit) are equivalent for bases including the boundary 
nodes.

Since the product rule has been used for the continuous energy estimate, a split
form of the equation will be used in order to get a similar estimate at the
semidiscrete level.
A standard split form SBP semidiscretisation of \eqref{eq:var-lin-adv} can be
written using numerical fluxes on one element as
\begin{equation}
\label{eq:standard-simplified}
  \partial_t \vec{u}
  =
  - \frac{1}{2} \mat{D} \mat{a} \vec{u}
  - \frac{1}{2} \mat{a} \mat{D} \vec{u}
  - \frac{1}{2} \mat{u} \mat{D} \vec{a}
  - \mat{M}[^{-1}] \mat{R}[^T] \mat{B} \left(
    \vecfnum
    - \mat{R} \mat{a} \vec{u}
  \right).
\end{equation}
Here, the first three terms on the right hand side approximate the split form
$-\frac{1}{2} \left( \partial_x(a u) + a (\partial_x u) + (\partial_x a) u \right)$
of the flux term $-\partial_x (a u)$ of \eqref{eq:var-lin-adv}. The last term
$- \mat{M}[^{-1}] \mat{R}[^T] \mat{B} \left(\vecfnum - \mat{R} \mat{a} \vec{u} \right)$
is consistent with zero, since both $\vecfnum$ and $\mat{R} \mat{a} \vec{u}$
approximate the product of $a$ and $u$ at the boundaries of an element. It
corresponds to a simultaneous approximation term in finite difference methods, 
see also \ref{sec:translation-rules} and \ref{subsec:comparison-nordstrom-ruggiu}.
Thus, it is used to enforce the boundary conditions weakly.
A translation of \eqref{eq:standard-simplified} to the notation of
\cite{nordstrom2017conservation} is given by
\eqref{eq:standard-diagonal-incl-boundaries-nordstrom}.

Here, the numerical flux can be either an interior flux, evaluated at the boundary
between two elements in a multi-block discretisation, or an exterior flux, evaluated
at the boundary of the domain $(x_L, x_R)$. In the first case, $\fnum$ is evaluated
using the values of $a$ and $u$ from the two neighbouring cells. In the second
case, the boundary condition is inserted as one argument instead. If periodic
boundary conditions are used, all boundaries are treated as interior boundaries.

\subsubsection{Conservation}

Investigating conservation, the semidiscrete equation \eqref{eq:standard-simplified}
is multiplied with $\vec{1}^T \mat{M}$, where $\vec{1} = (1, \dots, 1)^T$ is the
representation of the constant function $x \mapsto 1$. Thus, using the SBP
property \eqref{eq:SBP},
\begin{equation}
\begin{aligned}
  &
  \vec{1}^T \mat{M}\partial_t \vec{u}
  =
  - \frac{1}{2} \vec{1}^T \mat{M} \mat{D} \mat{a} \vec{u}
  - \frac{1}{2} \vec{1}^T \mat{M} \mat{a} \mat{D} \vec{u}
  - \frac{1}{2} \vec{1}^T \mat{M} \mat{u} \mat{D} \vec{a}
  - \vec{1}^T \mat{R}[^T] \mat{B} \left(
    \vecfnum
    - \mat{R} \mat{a} \vec{u}
  \right)
  \\
  =&
  - \frac{1}{2} \vec{1}^T \mat{R}[^T] \mat{B} \mat{R} \mat{a} \vec{u}
  + \frac{1}{2} \underbrace{\vec{1}^T \mat{D}[^T]}_{=\vec{0}^T} \mat{M} \vec{a}
  - \frac{1}{2} \vec{1}^T \mat{M} \mat{a} \mat{D} \vec{u}
  - \frac{1}{2} \vec{1}^T \mat{M} \mat{u} \mat{D} \vec{a}
  - \vec{1}^T \mat{R}[^T] \mat{B} \left(
    \vecfnum
    - \mat{R} \mat{a} \vec{u}
  \right).
\end{aligned}
\end{equation}
Since the diagonal multiplication operators are self-adjoint with respect to $\mat{M}$,
$\mat{M} \mat{a} = \mat{a} \mat{M}$ and $\mat{M} \mat{u} = \mat{u} \mat{M}$.
Hence, the SBP property \eqref{eq:SBP} can be applied again, yielding
\begin{equation}
\begin{aligned}
  &
  \vec{1}^T \mat{M}\partial_t \vec{u}
  =
  - \frac{1}{2} \vec{1}^T \mat{a} \mat{M} \mat{D} \vec{u}
  - \frac{1}{2} \vec{1}^T \mat{u} \mat{M} \mat{D} \vec{a}
  - \vec{1}^T \mat{R}[^T] \mat{B} \vecfnum
  + \frac{1}{2} \vec{1}^T \mat{R}[^T] \mat{B} \mat{R} \mat{a} \vec{u}
  \\
  =&
  - \frac{1}{2} \vec{a}^T \mat{R}[^T] \mat{B} \mat{R} \vec{u}
  + \frac{1}{2} \vec{a}^T \mat{D}[^T] \mat{M} \vec{u}
  - \frac{1}{2} \vec{u}^T \mat{M} \mat{D} \vec{a}
  - \vec{1}^T \mat{R}[^T] \mat{B} \vecfnum
  + \frac{1}{2} \vec{1}^T \mat{R}[^T] \mat{B} \mat{R} \mat{a} \vec{u}.
\end{aligned}
\end{equation}
Since boundary nodes are included,
\begin{equation}
  \vec{1}^T \mat{R}[^T] \mat{B} \mat{R} \mat{a} \vec{u}
  =
  1 \cdot a_p u_p - 1 \cdot a_0 u_0
  =
  a_p \cdot u_p - a_0 \cdot u_0
  =
  \vec{a}^T \mat{R}[^T] \mat{B} \mat{R} \vec{u}.
\end{equation}
Therefore, the desired conservation property follows
\begin{equation}
\label{eq:standard-simplified-conservation}
  \vec{1}^T \mat{M}\partial_t \vec{u}
  =
  - \vec{1}^T \mat{R}[^T] \mat{B} \vecfnum.
\end{equation}
Since the numerical flux $\fnum$ is defined per boundary, the contributions of
two neighbouring
elements cancel out and after summing up all elemental contributions, the two
boundary terms
\begin{equation}
  \fnum\big(g_L, u(\xmin)\big) - \fnum\big(u(\xmax),*\big)
\end{equation}
remain. Using the upwind numerical flux $\fnum(u_-,u_+) = a(x) \, u_-$, these boundary
terms become
\begin{equation}
  \fnum\big(g_L, u(\xmin)\big) - \fnum\big(u(\xmax),*\big)
  =
  a(\xmin) g_L - a(\xmax) u(\xmax),
\end{equation}
just as in the continuous case \eqref{eq:standard-conservation}. Here, the upwind
numerical fluxes \eqref{eq:fnum-edge-upwind}, \eqref{eq:fnum-split-upwind},
and \eqref{eq:fnum-unsplit-upwind} are equivalent, since boundary nodes are included.
Therefore, $a$ is treated as a constant in the following computations.

\subsubsection{Stability}

Mimicking the continuous estimate, the semidiscretisation \eqref{eq:standard-simplified}
is multiplied with $\vec{u}^T \mat{M}$, resulting due to the SBP property 
\eqref{eq:SBP} in
\begin{equation}
\begin{aligned}
  \vec{u}^T \mat{M} \partial_t \vec{u}
  &=
  - \frac{1}{2} \vec{u}^T \mat{M} \mat{D} \mat{a} \vec{u}
  - \frac{1}{2} \vec{u}^T \mat{a} \mat{M} \mat{D} \vec{u}
  - \frac{1}{2} \vec{u}^T \mat{u} \mat{M} \mat{D} \vec{a}
  - \vec{u}^T \mat{R}[^T] \mat{B} \left(
    \vecfnum
    - \mat{R} \mat{a} \vec{u}
  \right)
  \\
  &=
  - \frac{1}{2} \vec{u}^T \mat{u} \mat{M} \mat{D} \vec{a}
  - \vec{u}^T \mat{R}[^T] \mat{B} \vecfnum
  + \frac{1}{2} \vec{u}^T \mat{R}[^T] \mat{B}\mat{R} \mat{a} \vec{u}.
\end{aligned}
\end{equation}
Inserting the upwind numerical flux $\fnum(u_-,u_+) = a \, u_-$ 
(where $a$ is treated as a constant, since boundary nodes are included)
for a semidiscretisation using only one element results in
\begin{equation}
\label{eq:standard-simplified-stability}
\begin{aligned}
  \od{}{t} \norm{ \vec{u} }_{\mat{M}}^2
  =&\,
  2 \vec{u}^T \mat{M} \partial_t \vec{u}
  =
  - \vec{u}^T \mat{u} \mat{M} \mat{D} \vec{a}
  - 2 \vec{u}^T \mat{R}[^T] \mat{B} \vecfnum
  + \vec{u}^T \mat{R}[^T] \mat{B}\mat{R} \mat{a} \vec{u}
  \\
  =&
  - \vec{u}^T \mat{u} \mat{M} \mat{D} \vec{a}
  - 2 \left( a(\xmax) u(\xmax)^2 - a(\xmin) u(\xmin) g_L \right)
  + a(\xmax) u(\xmax)^2 - a(\xmin) u(\xmin)^2
  \\
  =&
  + a(\xmin) g_L^2
  - a(\xmax) u(\xmax)^2
  - \vec{u}^T \mat{u} \mat{M} \mat{D} \vec{a}
  - a(\xmin) \left( u(\xmin) - g_L \right)^2.
\end{aligned}
\end{equation}
The first three terms mimic the continuous counterpart \eqref{eq:standard-stability} 
exactly and the fourth one is an additional stabilisation term, since $a > 0$.

If multiple elements are used, the boundary terms
$- 2 \vec{u}^T \mat{R}[^T] \mat{B} \vecfnum
+ \vec{u}^T \mat{R}[^T] \mat{B}\mat{R} \mat{a} \vec{u}$
in \eqref{eq:standard-simplified-stability} of two neighbouring elements can
be added. 
More precisely, the terms of right hand side of the element to the left
of a given interior boundary as well as the terms of the left hand side of the
element to the right can be added, see also \autoref{fig:notation-numerical-fluxes}.
Thus, the numerical flux remains the same and the indices of the interpolated values
($\vec{u}^T \mat{R}[^T]$ and $\vec{u}^T \mat{R}[^T] \mat{B}\mat{R} \mat{a} \vec{u}$) 
are transformed via $R \mapsto -, L \mapsto +$. This results in the following 
contribution to the rate of change of the energy by one inter-element boundary
(treating $a$ again as a constant, given at the boundary)
\begin{equation}
\label{eq:standard-stability-diagonal-including-boundary-one-interior-boundary}
  2 (u_+ - u_-) \fnum(u_-,u_+) - \left( a u_+^2 - a u_-^2 \right).
\end{equation}
Here, the indices $-$ and $+$ indicate contributions from the elements to the
left and right of a given inter-element boundary, respectively. If the central
flux $\fnum(u_-,u_+) = a (u_- + u_+) / 2$ is used, this contribution becomes
\begin{equation}
  2 (u_+ - u_-) \fnum(u_-,u_+) - \left( a u_+^2 - a u_-^2 \right)
  =
  (u_+ - u_-) a (u_- + u_+) - \left( a u_+^2 - a u_-^2 \right)
  =
  0,
\end{equation}
where $a$ is treated as a constant as before.
Similarly, if the upwind flux $\fnum(u_-,u_+) = a \, u_-$ is used, the contribution is
\begin{equation}
  2 (u_+ - u_-) \fnum(u_-,u_+) - \left( a u_+^2 - a u_-^2 \right)
  =
  2 (u_+ - u_-) a u_- - \left( a u_+^2 - a u_-^2 \right)
  =
  - a (u_- - u_+)^2 \leq 0,
\end{equation}
resulting in an additional stabilisation.
This proves
\begin{theorem}[cf. Propositions 4.4, 5.2, 6.1, and 7.1 of \citet{nordstrom2017conservation}]
\label{thm:standard-simplified}
  Using diagonal norm SBP discretisations including the boundary nodes, the
  semidiscretisation \eqref{eq:standard-simplified} of \eqref{eq:var-lin-adv}
  is both conservative and stable across elements if the upwind numerical flux
  $\fnum(u_-, u_+) = a \, u_-$ is used at the exterior boundaries.
  
  If multiple elements are used, the numerical flux at inter-element boundaries
  can be chosen to be the upwind one (adding additional dissipation) or the central
  flux $\fnum(u_-,u_+) = a \frac{u_- + u_+}{2}$ (without additional dissipation).
\end{theorem}
Since boundary nodes are included, the upwind numerical fluxes \eqref{eq:fnum-edge-upwind},
\eqref{eq:fnum-split-upwind}, \eqref{eq:fnum-unsplit-upwind} as well as the 
central fluxes \eqref{eq:fnum-edge-central}, \eqref{eq:fnum-split-central},
\eqref{eq:fnum-unsplit-central} are equivalent, respectively. Therefore, $a$ can
be treated as a constant.

\subsection{General SBP Discretisations}
\label{subsec:standard-general}

Here, a general SBP discretisation will be considered. Thus, multiplication
operators are not necessarily self-adjoint. Furthermore, restriction to the boundary and
multiplication do not commute in general, since boundary nodes may not be included.
Therefore, corrections of the semidiscretisation \eqref{eq:standard-simplified}
have to be introduced. A suitable choice is
\begin{equation}
\label{eq:standard-general}
  \partial_t \vec{u}
  =
  - \frac{1}{2} \mat{D} \mat{a} \vec{u}
  - \frac{1}{2} \mat{a}[^*] \mat{D} \vec{u}
  - \frac{1}{2} \mat{u}[^*] \mat{D} \vec{a}
  - \mat{M}[^{-1}] \mat{R}[^T] \mat{B} \left(
    \vecfnum
    - \frac{1}{2} \mat{R} \mat{a} \vec{u}
    - \frac{1}{2} \left(\mat{R} \vec{a}\right) \cdot \left(\mat{R} \vec{u}\right)
  \right).
\end{equation}
Compared to \eqref{eq:standard-simplified}, two types of corrections have been
performed. At first, the multiplication operators $\mat{a}$ and $\mat{u}$ applied
to derivatives are substituted by the corresponding $\mat{M}$-adjoint operators
in order to allow non-diagonal norms. Secondly, the interpolation term
$\mat{R} \mat{a} \vec{u}$ has been replaced by a split form similar to the one used
for the flux derivative. Here and in the following, 
$\left(\mat{R} \vec{a}\right) \cdot \left(\mat{R} \vec{u}\right)$ denotes the 
componentwise product of the two vectors $\mat{R} \vec{a}$ and $\mat{R} \vec{u}$.
Due to the accuracy of the interpolation operator $\mat{R}$, the order of accuracy 
of the boundary terms should be the same as for the simplified semidiscretisation 
\eqref{eq:standard-simplified}.

As mentioned in section~\ref{subsec:SQTs-and-fnum}, the numerical fluxes (edge based,
split, and unsplit) that are equivalent for bases including the boundary nodes
will behave differently in the general context. Thus, the various choices have
to be considered separately.

\subsubsection{Conservation}

Multiplication of \eqref{eq:standard-general} with $\vec{1}^T \mat{M}$ results in
\begin{equation}
\begin{aligned}
  \vec{1}^T \mat{M} \partial_t \vec{u}
  =&
  - \frac{1}{2} \vec{1}^T \mat{M} \mat{D} \mat{a} \vec{u}
  - \frac{1}{2} \vec{1}^T \mat{M} \mat{a}[^*] \mat{D} \vec{u}
  - \frac{1}{2} \vec{1}^T \mat{M} \mat{u}[^*] \mat{D} \vec{a}
  \\&
  - \vec{1}^T \mat{R}[^T] \mat{B} \left(
    \vecfnum
    - \frac{1}{2} \mat{R} \mat{a} \vec{u}
    - \frac{1}{2} \left(\mat{R} \vec{a}\right) \cdot \left(\mat{R} \vec{u}\right)
  \right).
\end{aligned}
\end{equation}
Inserting the definition of the adjoint
$\mat{u}[^*] = \mat{M}[^{-1}] \mat{u}[^T] \mat{M}$ \eqref{eq:M-adjoint} and the
SBP property \eqref{eq:SBP} yield
\begin{equation}
\label{eq:standard-general-conservation}
\begin{aligned}
  \vec{1}^T \mat{M} \partial_t \vec{u}
  =&
  - \frac{1}{2} \vec{1}^T \mat{M} \mat{D} \mat{a} \vec{u}
  - \frac{1}{2} \vec{a}^T \mat{M} \mat{D} \vec{u}
  - \frac{1}{2} \vec{u}^T \mat{M} \mat{D} \vec{a}
  \\&
  - \vec{1}^T \mat{R}[^T] \mat{B} \left(
    \vecfnum
    - \frac{1}{2} \mat{R} \mat{a} \vec{u}
    - \frac{1}{2} \left(\mat{R} \vec{a}\right) \cdot \left(\mat{R} \vec{u}\right)
  \right)
  \\
  =&
  - \frac{1}{2} \vec{a}^T \mat{R}[^T] \mat{B} \mat{R} \vec{u}
  - \vec{1}^T \mat{R}[^T] \mat{B} \left(
    \vecfnum
    - \frac{1}{2} \left(\mat{R} \vec{a}\right) \cdot \left(\mat{R} \vec{u}\right)
  \right)
  =
  - \vec{1}^T \mat{R}[^T] \mat{B} \vecfnum,
\end{aligned}
\end{equation}
since $\vec{a}^T \mat{R}[^T] \mat{B} \mat{R} \vec{u} = a_R u_R - a_L u_L = 
\vec{1}^T \mat{R}[^T] \mat{B} \left(\mat{R} \vec{a}\right) \cdot 
\left(\mat{R} \vec{u}\right)$.
This is the same result as in the simplified case
\eqref{eq:standard-simplified-conservation}. Using again an upwind numerical flux,
the continuous property \eqref{eq:standard-conservation} is mimicked on a
semidiscrete level.

\subsubsection{Stability}

As in the case of the simplified semidiscretisation \eqref{eq:standard-simplified},
equation \eqref{eq:standard-general} is multiplied with $\vec{u}^T \mat{M}$,
resulting due to the SBP property \eqref{eq:SBP} in
\begin{equation}
\label{eq:standard-general-stability-first-step}
\begin{aligned}
  \vec{u}^T \mat{M} \partial_t \vec{u}
  =&
  - \frac{1}{2} \vec{u}^T \mat{M} \mat{D} \mat{a} \vec{u}
  - \frac{1}{2} \vec{u}^T \mat{a}[^T] \mat{M} \mat{D} \vec{u}
  - \frac{1}{2} \vec{u}^T \mat{u}[^T] \mat{M} \mat{D} \vec{a}
  \\&
  - \vec{u}^T \mat{R}[^T] \mat{B} \left(
    \vecfnum
    - \frac{1}{2} \mat{R} \mat{a} \vec{u}
    - \frac{1}{2} \left(\mat{R} \vec{a}\right) \cdot \left(\mat{R} \vec{u}\right)
  \right)
  \\
  =&
  - \frac{1}{2} \vec{u}^T \mat{u}[^T] \mat{M} \mat{D} \vec{a}
  - \vec{u}^T \mat{R}[^T] \mat{B} \left(
    \vecfnum
    - \frac{1}{2} \left(\mat{R} \vec{a}\right) \cdot \left(\mat{R} \vec{u}\right)
  \right).
\end{aligned}
\end{equation}
In order to enforce the boundary conditions as before, the upwind numerical
fluxes \eqref{eq:fnum-edge-upwind}, \eqref{eq:fnum-split-upwind}, and 
\eqref{eq:fnum-unsplit-upwind} will be considered in the following.
\begin{itemize}
  \item 
  Using the edge based upwind flux $\fnum = a(x) \, u_-$ \eqref{eq:fnum-edge-upwind},
  \eqref{eq:standard-general-stability-first-step} becomes
  \begin{equation}
  \begin{aligned}
    \od{}{t} \norm{ \vec{u} }_{\mat{M}}^2
    &=
    2 \vec{u}^T \mat{M} \partial_t \vec{u}
    =
    - \vec{u}^T \mat{u}[^T] \mat{M} \mat{D} \vec{a}
    - 2 \vec{u}^T \mat{R}[^T] \mat{B} \vecfnum
    + \vec{u}^T \mat{R}[^T] \mat{B} \left(\mat{R} \vec{a}\right) 
      \cdot \left(\mat{R} \vec{u}\right)
    \\
    &=
    - \vec{u}^T \mat{u}[^T] \mat{M} \mat{D} \vec{a}
    - 2 \left( a(\xmax) u_R^2 - a(\xmin) u_L g_L \right)
    + a_R u_R^2 - a_L u_L^2,
  \end{aligned}
  \end{equation}
  where the indices $L,R$ indicate interpolations to the left and right boundary,
  respectively. If the interpolation of $a$ is exact, i.e.
  \begin{equation}
  \label{eq:interpolated-speed-exact}
    \mat{R} \vec{a} = (a_L, a_R)^T \stackrel{!}{=} (a(\xmin), a(\xmax))^T,
  \end{equation}
  where the first equality holds by definition and the second one describes the
  desired exactness, this can be rewritten as
  \begin{equation}
    \od{}{t} \norm{ \vec{u} }_{\mat{M}}^2
    =
    + a(\xmin) g_L^2
  - a(\xmax) u_R^2
  - \vec{u}^T \mat{u}[^T] \mat{M} \mat{D} \vec{a}
  - a(\xmin) \left( u_L - g_L \right)^2,
  \end{equation}
  mimicking again the continuous counterpart \eqref{eq:standard-stability} with an
  additional stabilising term.

  \item
  Using instead the split upwind flux $\fnum = a_- u_-$ \eqref{eq:fnum-split-upwind},
  \eqref{eq:standard-general-stability-first-step} becomes
  \begin{equation}
  \begin{aligned}
    \od{}{t} \norm{ \vec{u} }_{\mat{M}}^2
    &=
    - \vec{u}^T \mat{u}[^T] \mat{M} \mat{D} \vec{a}
    - 2 \left( a_R u_R^2 - a_L u_L g_L \right)
    + a_R u_R^2 - a_L u_L^2
    \\
    &=
    + a_L g_L^2
    - a_R u_R^2
    - \vec{u}^T \mat{u}[^T] \mat{M} \mat{D} \vec{a}
    - a_L \left( u_L - g_L \right)^2,
  \end{aligned}
  \end{equation}
  mimicking the continuous counterpart \eqref{eq:standard-stability} with an
  additional stabilising term if the interpolated speeds $a_L, a_R$ are positive.
  This property does not hold in general, since (Lagrange) interpolants of
  positive functions can be negative at some points. Moreover, in order to
  mimic \eqref{eq:standard-stability} reliably, the interpolated speeds should
  be exact \eqref{eq:interpolated-speed-exact}.

  \item
  Finally, if the unsplit upwind flux $\fnum = (a u)_-$ \eqref{eq:fnum-unsplit-upwind}
  is used, \eqref{eq:standard-general-stability-first-step} becomes
  \begin{equation}
    \od{}{t} \norm{ \vec{u} }_{\mat{M}}^2
    =
    - \vec{u}^T \mat{u}[^T] \mat{M} \mat{D} \vec{a}
    - 2 \left( (a u)_R u_R - (a u)_L g_L \right)
    + a_R u_R^2 - a_L u_L^2.
  \end{equation}
  However, since the restriction $(a u)_{L/R}$ of the product can not be compared
  to the product $a_{L/R} u_{L/R}$ of the restrictions, the continuous counterpart
  \eqref{eq:standard-stability} will not be mimicked in general.
\end{itemize}
Consequently, the unsplit flux should not be used. The edge based and the split
numerical fluxes can be used if the interpolated speeds are exact. In this case,
both fluxes are equivalent.

\begin{remark}
\label{rem:use-Lobatto-nodes-and-interpolate-to-Gauss}
  The exactness of the interpolated speed \eqref{eq:interpolated-speed-exact} can
  be achieved in nodal DG methods (corresponding to the analytical setting
  described in section~\ref{subsec:analytical-setting}) as follows. 
  If a nodal basis using $p+1$ points is used to represent polynomials of 
  degree~$\leq p$, the standard procedure to compute the representation $\vec{a}$
  is to evaluate the function $a$ at the collocation nodes, i.e. 
  $\vec{a}_i = a(\xi_i)$. Instead, the speed $a(x)$ can be evaluated at $p+1$ 
  points including the boundary nodes (e.g. Gauss-Lobatto nodes). Afterwards, 
  the unique interpolating polynomial can be evaluated at the nodes used in the 
  basis not including the boundaries (e.g. Gauss nodes).
\end{remark}

Again, if multiple elements are used, the boundary terms
$- 2 \vec{u}^T \mat{R}[^T] \mat{B} \vecfnum
+ \vec{u}^T \mat{R}[^T] \mat{B} \left(\mat{R} \vec{a}\right) 
\cdot \left(\mat{R} \vec{u}\right)$
in~\eqref{eq:standard-general-stability-first-step} of two neighbouring elements
can be added, resulting similarly to 
\eqref{eq:standard-stability-diagonal-including-boundary-one-interior-boundary}
in the following contribution of one boundary to the rate of change of the energy, 
\begin{equation}
\label{eq:standard-general-stability-a-general-boundary-contribution}
  2 (u_+ - u_-) \fnum(u_-,u_+) - \left( a_+ u_+^2 - a_- u_-^2 \right),
\end{equation}
where the indices $-$ and $+$ indicate again contributions from the elements to
the left and right of a given inter-element boundary, respectively.

If the interpolation of $a$ is continuous across boundaries, i.e. $a_+ = a_- =: a$, 
and the (edge based or split) central flux $\fnum(u_-,u_+) = a (u_- + u_+) / 2$
is used, this contribution vanishes again.
If the split upwind flux $\fnum = a_- \, u_-$ \eqref{eq:fnum-split-upwind} is used, 
the contribution is
\begin{equation}
\begin{aligned}
  &
  2 (u_+ - u_-) \fnum(u_-,u_+) - \left( a_+ u_+^2 - a_- u_-^2 \right)
  =
  2 (u_+ - u_-) a_- u_- - \left( a_+ u_+^2 - a_- u_-^2 \right)
  \\
  =&
  - a_- u_-^2 + 2 a_- u_- u_+ - a_+ u_+^2
  =
  - a_- (u_- - u_+)^2 + (a_- - a_+) u_+^2.
\end{aligned}
\end{equation}
If $a_- \neq a_+$, this can be positive, e.g. for $u_- = u_+ > 0$ and $a_- - a_+ > 0$.
Thus, $a$ should be discretised as being continuous across boundaries, i.e.
$a_- = a_+$, in order to get the desired stabilisation. 
This proves
\begin{theorem}
\label{thm:standard-general}
  Consider the semidiscretisation \eqref{eq:standard-general} of \eqref{eq:var-lin-adv}
  using general SBP operators where the speed $a > 0$ is discretised such that the
  restrictions $a_{L/R}$ to the boundary are positive, the discretised speed is
  continuous across elements, and the interpolations at the boundaries of the 
  domain $(\xmin, \xmax)$ are exact.
  
  Then, the edge based fluxes \eqref{eq:fnum-edge-central}, \eqref{eq:fnum-edge-upwind}
  are equivalent to the split ones  \eqref{eq:fnum-split-central}, \eqref{eq:fnum-split-upwind},
  respectively.
  \begin{itemize}
    \item 
    If the upwind flux $\fnum = a_- u_-$ is used at the exterior boundaries,
    the semidiscretisation \eqref{eq:standard-general} is both conservative and
    stable.
    
    \item
    If multiple elements are used, the numerical flux at the interior boundaries
    can be chosen to be the upwind one (adding additional dissipation) or the 
    central flux $\fnum = (a_- u_- + a_+ u_+) / 2$ (without additional dissipation).
  \end{itemize}
\end{theorem}
In general, the assumptions on the interpolated speed are very strong. They can
be achieved easily if the speed $a$ is a polynomial and the interpolations are
exact for polynomials of the same degree as $a$. Otherwise, some care has to
be taken. In nodal DG methods, the assumptions on the interpolated speed can be
enforced as described in Remark~\ref{rem:use-Lobatto-nodes-and-interpolate-to-Gauss}.

A comparison of these results with the ones of \cite{nordstrom2017conservation}
is given in \ref{subsec:comparison-nordstrom-ruggiu}.

\begin{remark}
  It might be tempting to choose the numerical flux $\fnum$ such that the
  contribution \eqref{eq:standard-general-stability-a-general-boundary-contribution} 
  of one inter-element boundary vanishes, cf. the discussion before and after
  \eqref{eq:standard-stability-diagonal-including-boundary-one-interior-boundary}.
  Naively solving the resulting equation for $\fnum$ results in
  \begin{equation}
  \label{eq:naively-solving-for-entropy-conservative-flux}
  \begin{aligned}
    \fnum(u_-,u_+)
    &=
    \frac{a_+ u_+^2 - a_- u_-^2}{2 (u_+ - u_-)}
    =
    \frac{a_+ u_+^2 - a_+ u_+ u_- + a_- u_+ u_- - a_- u_-^2}{2 (u_+ - u_-)}
    + \frac{a_+ u_+ u_- - a_- u_+ u_-}{2 (u_+ - u_-)}
    \\
    &=
    \frac{a_+ u_+ + a_- u_-}{2} + \frac{1}{2} (a_+ - a_-) \frac{u_+ u_-}{u_+ - u_-}.
  \end{aligned}
  \end{equation}
  However, this is only possible for $a_+ = a_-$ in general, since otherwise the
  case $u_+ = u_-$ cannot be handled.
  Similarly, a local Lax-Friedrichs flux
  $\fnum(u_-,u_+) = \frac{a_- u_- + a_+ u_+}{2} - \frac{\lambda}{2} (u_+ - u_-)$,
  $\lambda = \max\set{ \abs{a_-}, \abs{a_+} }$, does not yield the desired estimate, since
  \begin{equation}
  \begin{aligned}
    &
    2 (u_+ - u_-) \fnum(u_-,u_+) - \left( a_+ u_+^2 - a_- u_-^2 \right)
    \\
    =&\,
    (u_+ - u_-) (a_+ u_+ + a_- u_-) - \lambda (u_+ - u_-)^2
    - \left( a_+ u_+^2 - a_- u_-^2 \right)
    \\
    =&\,
    a_+ u_+^2 - a_+ u_+ u_- + a_- u_+ u_- - a_- u_-^2 - \lambda (u_+ - u_-)^2
    - a_+ u_+^2 + a_- u_-^2
    \\
    =&\,
    - a_+ u_+ u_- + a_- u_+ u_- - \lambda (u_+ - u_-)^2.
  \end{aligned}
  \end{equation}
  This can become positive, e.g. if $u_+ = u_- > 0$ and $a_- - a_+ > 0$.
  Both problems apply only if $a_- \neq a_+$. Thus, if the discretised speed is
  continuous as supposed in \autoref{thm:standard-general}, both the local
  Lax-Friedrichs flux and \eqref{eq:naively-solving-for-entropy-conservative-flux} 
  can be used,   stressing the importance of the condition $a_- = a_+$.
  Note the analogy to the continuous energy estimate 
  \eqref{eq:standard-continuous-energy-estimate},
  requiring $\partial_x a \in L^\infty$.
\end{remark}

\section{Weighted \texorpdfstring{$L^2$}{L²} Estimates}
\label{sec:weighted-L2}

In this section, the weighted $L^2$ estimates of \citet{manzanero2017insights}
will be reviewed and suitable formulations for general SBP discretisations will
be presented, resulting in both conservative and stable methods.
Here, only the estimates for the conservative equation~\eqref{eq:var-lin-adv}
will be considered. Results for a nonconservative equation can be found in
\ref{sec:weighted-nonconservative}.

Consider again the linear advection equation
\begin{equation}
\label{eq:conservative}
\begin{aligned}
  \partial_t u(t,x) + \partial_x \left( a(x) u(t,x) \right) &= 0,
  \qquad&& t > 0, \; x \in (\xmin, \xmax),
  \\
  u(t,\xmin) &= g_L(t),
  \qquad&& t \geq 0,
  \\
  u(0,x) &= u_0(x),
  \qquad&& x \in [\xmin,\xmax],
\end{aligned}
\end{equation}
with variable speed $a(x) > 0$ and compatible initial and boundary conditions
$u_0$, $g_L$. Here, stability will be investigated in a weighted $L^2$ norm
given by the scalar product
\begin{equation}
\label{eq:L2-a}
  \scp{u}{v}_{L^2_a}
  =
  \int_\xmin^\xmax a \, u \, v.
\end{equation}

\subsection{Continuous Estimates}
\label{subsec:weighted-conservative-continuous}

The conservation property \eqref{eq:standard-conservation} is the same as in the
previous \autoref{sec:standard-L2}. Multiplying equation \eqref{eq:conservative}
with~$a$ and integrating results due to integration-by-parts in
\begin{equation}
\label{eq:weighted-conservative-stability}
\begin{gathered}
  \od{}{t} \int_\xmin^\xmax a \, u^2
  =
  2 \int_\xmin^\xmax a \, u \, \partial_t u
  =
  - 2 \int_\xmin^\xmax a \, u \, \partial_x(a \, u)
  =
  - a^2 \, u^2 \big|_\xmin^\xmax,
  \\
  \implies
  \od{}{t} \norm{u}_{L^2_a}^2 
  = a(\xmin)^2 g_L^2 - a(\xmax)^2 u(\xmax)^2 
  \leq a(\xmin)^2 g_L^2.
\end{gathered}
\end{equation}
Thus, the weighted $L^2$ norm fulfils
\begin{equation}
\label{eq:weighted-conservative-energy-estimate}
  \norm{u(t)}_{L^2_a}^2
  \leq
  \norm{u_0}_{L^2_a}^2 + \int_0^t a(\xmin)^2 g_L(\tau)^2 \dif \tau.
\end{equation}
As described in \cite[Section 4]{manzanero2017insights}, this can be translated
by the equivalence of norms
\begin{equation}
  \min_{x \in [\xmin,\xmax]} \set{a(x)} \norm{u}_{L^2}^2
  \leq
  \norm{u(t)}_{L^2_a}^2
  \leq
  \max_{x \in [\xmin,\xmax]} \set{a(x)} \norm{u}_{L^2}^2
\end{equation}
to the following $L^2$ bound on the solution $u$.
\begin{equation}
  \min_{x \in [\xmin,\xmax]} \set{a(x)} \norm{u(t)}_{L^2}^2
  \leq
  \max_{x \in [\xmin,\xmax]} \set{a(x)} \norm{u_0}_{L^2}^2
  + \int_0^t a(\xmin)^2 g_L(\tau)^2 \dif \tau.
\end{equation}
Note that no continuity of the speed $a$ has to be assumed for these calculations,
contrary to the unweighted $L^2$ estimates in section~\ref{subsec:standard-continuous}.

\subsection{Semidiscrete Estimates}
\label{subsec:weighted-conservative-semidiscrete}

In this section, general SBP semidiscretisations of \eqref{eq:conservative} are 
considered. Since no product rule has been used for the estimate
\eqref{eq:weighted-conservative-stability} in the continuous setting, the following
(unsplit) form of the semidiscretisation will be considered.
\begin{equation}
\label{eq:weighted-conservative-semidiscretisation}
  \partial_t \vec{u}
  =
  - \mat{D} \mat{a} \vec{u}
  - \mat{M}[^{-1}] \mat{R}[^T] \mat{B} \left( \vecfnum - \mat{R} \mat{a} \vec{u} \right).
\end{equation}
Investigating conservation, multiplying the semidiscretisation
\eqref{eq:weighted-conservative-semidiscretisation} with $\vec{1}^T \mat{M}$
results due to the SBP property \eqref{eq:SBP} in
\begin{equation}
\label{eq:weighted-conservative-semidiscretisation-conservation}
  \vec{1}^T \mat{M} \partial_t \vec{u}
  =
  - \vec{1}^T \mat{M} \mat{D} \mat{a} \vec{u}
  - \vec{1}^T \mat{R}[^T] \mat{B} \left( \vecfnum - \mat{R} \mat{a} \vec{u} \right)
  =
  - \vec{1}^T \mat{R}[^T] \mat{B} \vecfnum,
\end{equation}
just as in the cases \eqref{eq:standard-simplified-conservation} and
\eqref{eq:standard-general-conservation} considered before.

Turning to stability, the discrete ``norm'' corresponding to $\norm{\cdot}_{L^2_a}$ is 
given by $\norm{\vec{u}}_{\mat{a}[^T] \mat{M}}^2 = \vec{u}^T \mat{a}[^T] \mat{M} \vec{u}$.
If $\mat{a}[^T] \mat{M}$ is symmetric and positive definite, this is a norm,
see also Remark~\ref{rem:norms} below. Assuming that $\norm{\cdot}_{\mat{a}[^T] \mat{M}}$
is indeed a norm, the semidiscrete energy estimate is obtained by multiplying
\eqref{eq:weighted-conservative-semidiscretisation} with $\vec{u}^T \mat{a}[^T] \mat{M}$
and inserting the SBP property \eqref{eq:SBP}, resulting in
\begin{equation}
\label{eq:weighted-conservative-semidiscretisation-stability}
\begin{aligned}
  \od{}{t} \norm{ \vec{u} }_{\mat{a}[^T]\mat{M}}^2
  &=
  2 \vec{u}^T \mat{a}[^T] \mat{M} \partial_t \vec{u}
  =
  - 2 \vec{u}^T \mat{a}[^T] \mat{M} \mat{D} \mat{a} \vec{u}
  - 2 \vec{u}^T \mat{a}[^T]  \mat{R}[^T] \mat{B} \left( \vecfnum - \mat{R} \mat{a} \vec{u} \right)
  \\
  &=
  \vec{u}^T \mat{a}[^T]  \mat{R}[^T] \mat{B} \mat{R} \mat{a} \vec{u}
  - 2 \vec{u}^T \mat{a}[^T]  \mat{R}[^T] \mat{B} \vecfnum.
\end{aligned}
\end{equation}
Thus, the contribution of one boundary between the cells on the left (index $-$)
and on the right (index $+$) to the rate of change of the total energy is given by
\begin{equation}
\label{eq:energy-contribution-conservative}
  2 \left( (a u)_+ - (a u)_- \right) \fnum 
  - \left( (a u)_+^2 - (a u)_-^2 \right).
\end{equation}
Similarly to section~\ref{sec:standard-L2}, several forms (edge-based, split, unsplit)
of the numerical fluxes could be used. However, since no split form of the equation has
been used in the semidiscretisation \eqref{eq:weighted-conservative-semidiscretisation},
only the unsplit flux can be expected to yield a stable scheme and is investigated
below. The other fluxes are considered in
Remark~\ref{rem:weighted-conservative-simplified} below.

Using the unsplit upwind flux $\fnum = (a u)_-$ \eqref{eq:fnum-unsplit-upwind},
the contribution \eqref{eq:energy-contribution-conservative} of one boundary
to the rate of change of the energy becomes
\begin{equation}
\begin{aligned}
  &
  2 \left( (a u)_+ - (a u)_- \right) \fnum 
  - \left( (a u)_+^2 - (a u)_-^2 \right)
  \\
  =&\,
  2 \left( (a u)_+ - (a u)_- \right) (a u)_-
  - \left( (a u)_+^2 - (a u)_-^2 \right)
  =
  - \left( (a u)_- - (a u)_+ \right)^2
  \leq 0.
\end{aligned}
\end{equation}
Similarly, using the unsplit central flux $\fnum = (a_- u_- + a_+ u_+) / 2$ 
\eqref{eq:fnum-unsplit-central} results in
\begin{equation}
\begin{aligned}
  &
  2 \left( (a u)_+ - (a u)_- \right) \fnum 
  - \left( (a u)_+^2 - (a u)_-^2 \right)
  \\
  =&\,
  \left( (a u)_+ - (a u)_- \right) \left( (a u)_- + (a u)_+ \right)
  - \left( (a u)_+^2 - (a u)_-^2 \right)
  =
  0.
\end{aligned}
\end{equation}
Considering the total rate of change of the energy in a bounded domain, using
the unsplit upwind numerical flux \eqref{eq:fnum-split-upwind} at the boundaries 
results due to \eqref{eq:weighted-conservative-semidiscretisation-stability} in
\begin{equation}
\begin{aligned}
  \od{}{t} \norm{ \vec{u} }_{\mat{a}[^T]\mat{M}}^2
  &=
  \vec{u}^T \mat{a}[^T]  \mat{R}[^T] \mat{B} \mat{R} \mat{a} \vec{u}
  - 2 \vec{u}^T \mat{a}[^T]  \mat{R}[^T] \mat{B} \vecfnum
  \\
  &=
  (a u)_R^2 - (a u)_L^2 - 2 (a u)_R^2 + 2 (a u)_L (a g)_L
  \\
  &=
  (a g)_L^2 - (a u)_R^2 - \left( (a g)_L - (a u)_L \right)^2
  \leq (a g)_L^2 - (a u)_R^2 \leq (a g)_L^2.
\end{aligned}
\end{equation}
This mimics the continuous counterpart \eqref{eq:weighted-conservative-stability}
with an additional stabilising term $- \left( (a g)_L - (a u)_L \right)^2 \leq 0$.
In the implementation of the numerical flux at the left (exterior) boundary,
the term $(a g)_L$ is computed as $a(\xmin) g_L$ in general, exactly as in the
continuous estimate \eqref{eq:weighted-conservative-stability}. This proves
\begin{theorem}
\label{thm:weighted-conservative}
  Using general SBP discretisations such that $\mat{a}[^T]\mat{M}$ is symmetric
  and positive definite, the semidiscretisation
  \eqref{eq:weighted-conservative-semidiscretisation} of \eqref{eq:conservative}
  is both conservative and stable (in the discrete norm 
  $\norm{ \cdot }_{\mat{a}[^T]\mat{M}}$) across elements if the unsplit upwind 
  numerical flux \eqref{eq:fnum-unsplit-upwind} is used at the exterior boundaries.
  
  If multiple elements are used, the numerical flux at inter-element boundaries
  can be chosen to be the unsplit upwind one (adding additional dissipation) 
  or the unsplit central flux \eqref{eq:fnum-unsplit-central} (without additional
  dissipation).
\end{theorem}

\begin{remark}
\label{rem:conservation-in-thm:weighted-conservative}
  If boundary nodes are not included, the interpolated values $(a u)_{L/R}$ will
  be different from $a(x_{L/R}) u_{L/R}$ in general. Thus, since the unsplit
  upwind flux \eqref{eq:fnum-unsplit-upwind} has to be used for the energy
  estimate, the conservation property \eqref{eq:standard-conservation} is not 
  mimicked perfectly, as can be seen in the numerical results of 
  section~\ref{subsec:conservation-properties}. Nevertheless, convergence has 
  been observed in all examples.
\end{remark}

\begin{remark}
\label{rem:weighted-conservative-simplified}
  If boundary nodes are included, $(a u)_\pm = a_\pm u_\pm$, since multiplication
  and restriction commute. Furthermore the three forms (edge, split, and unsplit)
  of the central (\eqref{eq:fnum-edge-central},\eqref{eq:fnum-split-central}, and 
  \eqref{eq:fnum-unsplit-central}) and upwind (\eqref{eq:fnum-edge-upwind},
  \eqref{eq:fnum-split-upwind}, and \eqref{eq:fnum-unsplit-upwind}) numerical 
  fluxes are equivalent, respectively. Thus, the other forms of the numerical
  fluxes can be used as well as described in \autoref{thm:weighted-conservative}.
  
  If boundary nodes are not included, the other fluxes can result in an energy growth.
  As an example, consider polynomials of degree $p = 1$ in a nodal basis using
  the two Gauss nodes $x_0 = -1 / \sqrt{3}, x_1 = 1 / \sqrt{3}$.
  The corresponding Lagrange polynomials are 
  $\phi_{0/1}(x) = \mp \frac{\sqrt{3}}{2} x + \frac{1}{2}$, satisfying
  $\phi_i(x_j) = \delta_{ij}$.
  If the coefficients of $\vec{a}$ are $a_0 = 1, a_1 = 2$, the difference 
  $a_R u_R - (a u)_R$ is given as 
  \begin{equation}
    \left( a_0 \phi_0(1) + a_1 \phi_1(1) \right) 
    \left( u_0 \phi_-(1) + u_1 \phi_1(1) \right) 
    - 
    \left( a_0 u_0 \phi_0(1) + a_1 u_1 \phi_1(1) \right) 
    = \frac{u_1 - u_0}{2}.
  \end{equation}
  Thus, there can be an arbitrary difference between the product of the interpolants
  and the interpolant of the product.
  Inserting the values $(a u)_+ = 0$ and $a_+ u_+ = 0$ into the contribution 
  \eqref{eq:energy-contribution-conservative} of one boundary to the energy rate 
  results in
  \begin{equation}
    2 \left( (a u)_+ - (a u)_- \right) \fnum 
    - \left( (a u)_+^2 - (a u)_-^2 \right)
    =
    \begin{cases}
      (a u)_-^2 - (a u)_- a_- u_-, 
      & \text{split central flux \eqref{eq:fnum-split-central}}, 
      \\
      (a u)_-^2 - 2 (a u)_- a_- u_-, 
      & \text{split upwind flux \eqref{eq:fnum-split-upwind}}.
    \end{cases}
  \end{equation}
  Thus, for the coefficients $u_0 = 3$ and $u_1 = 1$ of the polynomial in the left
  cell (index $-$), the energy rate contributions become
  \begin{equation}
    2 \left( (a u)_+ - (a u)_- \right) \fnum 
    - \left( (a u)_+^2 - (a u)_-^2 \right)
    =
    \begin{cases}
      \frac{5}{2} - \frac{\sqrt{3}}{2} > 0, 
      & \text{split central flux \eqref{eq:fnum-split-central}}, 
      \\
      \frac{3 \sqrt{3}}{2} - 2 > 0, 
      & \text{split upwind flux \eqref{eq:fnum-split-upwind}}.
    \end{cases}
  \end{equation}
\end{remark}

\begin{remark}
\label{rem:norms}
  Here, stability is given with respect to the discrete ``norm''
  $\norm{ \cdot }_{\mat{a}[^T]\mat{M}}$. In order to give a reliable stability
  estimate, this has to be a real norm, i.e. the matrix $\mat{a}[^T]\mat{M}$ has
  to be symmetric and positive definite. This is the case for positive $a$ and
  nodal SBP discretisations with diagonal mass matrix $\mat{M}$, e.g. nodal
  discontinuous Galerkin schemes using Gauss nodes. 
  Another possibility is given by modal bases using exact multiplication followed
  by $L^2$ projection \cite[Section~4.2, around equation (42)]{ranocha2017extended}.
\end{remark}

\begin{remark}
\label{rem:weighted-conservative-compared-to-standard}
  Compared to \autoref{thm:standard-general} in \autoref{subsec:standard-general},
  the assumption that $a$ is discretised as being continuous with positive values
  at the interfaces can be dropped. Contrary, the assumptions on the SBP operators/bases
  have been strengthened, since $\mat{a}[^T] \mat{M}$ has to be symmetric
  and positive definite. While \autoref{thm:standard-general} allows nodal
  bases with dense norms $\mat{M}$, they are in general not allowed in 
  \autoref{thm:weighted-conservative}. As an example, consider the roots of the
  Chebyshev polynomials of first kind, i.e. 
  $\xi_{k} = \cos\Big( \frac{2k+1}{2p+2} \pi \Big), k \in \set{0,\dots,p}$.
  The corresponding mass matrix is
  $\mat{M} = \frac{1}{6} \begin{psmallmatrix} 5  & 1 \\ 1 & 5 \end{psmallmatrix}$.
  The product $\mat{a} \mat{M}$ of $\mat{a} = \diag{a_0, a_1}$ and $\mat{M}$ is 
  in general not symmetric. Furthermore, even its symmetric part
  $\frac{1}{2} \left( \mat{a} \mat{M} + (\mat{a} \mat{M})^T \right) = \frac{1}{6}
  \begin{psmallmatrix} 10 a_0 & a_0 + a_1 \\ a_0 + a_1 & 10 a_1 \end{psmallmatrix}$
  is not positive definite in general, even if $a_0,a_1 > 0$. Indeed, choosing
  $a_0 = 1$ and $a_1 = 100$, the eigenvalues of $\mat{a} \mat{M} + (\mat{a} \mat{M})^T$
  are $\lambda_\pm = \big(505 \pm \sqrt{255226} \big)/6$ and 
  $\lambda_- \approx -0.03$ is negative.
\end{remark}

A comparison of these results with the ones of \cite{manzanero2017insights}
is given in \ref{subsec:comparison-manzanero-conservative}.

\section{Numerical Results}
\label{sec:numerical-results}

In this section, some numerical experiments using the schemes constructed in the
previous sections~\ref{sec:standard-L2} and~\ref{sec:weighted-L2} will be conducted, 
including convergence studies and eigenvalue analyses.

\subsection{Convergence Studies}
\label{subsec:convergence-studies}

Consider the linear advection equation
\begin{equation}
\label{eq:var-lin-adv-without-boundary}
\begin{aligned}
  \partial_t u(t,x) + \partial_x \left( a(x) u(t,x) \right) &= 0,
  \qquad&& t > 0, \; x \in \R,
  \\
  u(0,x) &= u_0(x),
  \qquad&& x \in \R,
\end{aligned}
\end{equation}
with smooth speed $a(x) > 0$ and initial condition $u_0$. The solution
$u$ of this equation can be obtained by the method of characteristics, see e.g.
\cite[Chapter 3]{bressan2000hyperbolic}. 
Choosing the speed $a(x) = 1 + \cosh(x)$, the solution of 
\eqref{eq:var-lin-adv-without-boundary} can be written as
\begin{equation}
\label{eq:var-lin-adv-without-boundary-solution}
  u(t,x)
  =
  \frac{
    u_0\bigg( 2 \operatorname{artanh}\Big(- t + \tanh\big( \frac{x}{2} \big) \Big) \bigg)
  }{
    1 
    + 2 t \sinh\big( \frac{x}{2} \big) \cosh\big( \frac{x}{2} \big)
    - t^2 \big( \cosh \frac{x}{2} \big)^2
  }.
\end{equation}
The value of this analytical solution of the initial value problem 
\eqref{eq:var-lin-adv-without-boundary} will be imposed as boundary condition of
the initial boundary value problem \eqref{eq:var-lin-adv}, i.e. $g_L(t) = u(t,\xmin)$
with $u(t,x)$ as in \eqref{eq:var-lin-adv-without-boundary-solution}. Here, the 
initial condition $u_0(x) = \sin(\pi x)$ is used.

The domain $(\xmin,\xmax) = (-1,1)$ is divided into $N$ elements. On each element,
a polynomial of degree~$p$ is used to approximate the solution $u$. 
In the following, both Lobatto nodes and Gauss nodes (i.e. Legendre Gauss-Lobatto points
and Legendre Gauss points in the nomenclature of \cite[Chapter 3]{kopriva2009implementing}) 
will be used to represent the polynomials. While Lobatto nodes include the
boundary, Gauss nodes yield generalised SBP operators not including boundary nodes.
For Gauss nodes, the discretisation of the speed~$a$ has been performed both
directly on Gauss nodes and via interpolation from Lobatto nodes as described in
Remark~\ref{rem:use-Lobatto-nodes-and-interpolate-to-Gauss}.

Moreover, both the split form \eqref{eq:standard-general} and the unsplit form
\eqref{eq:weighted-conservative-semidiscretisation} will be considered, accompanied
with the corresponding central (\eqref{eq:fnum-split-central}, 
\eqref{eq:fnum-unsplit-central}) and upwind (\eqref{eq:fnum-split-upwind}, 
\eqref{eq:fnum-unsplit-upwind}) fluxes. While the numerical flux at the interior
boundaries will be varied, the flux at the exterior boundaries is always chosen
to be the correct upwind one.

The numerical solution is computed in the time interval $[0, 0.5]$ using the
ten stage, fourth order strong stability preserving method of \cite{ketcheson2008highly}.
In order to decrease the influence of the time discretisation, small time steps 
$\Delta t = 1 / \big( 100 (2p+1) N \big)$ have been used for the convergence tests.

\begin{table}
\centering
  \caption{Convergence results for \eqref{eq:var-lin-adv} using the speed 
           $a(x) = 1 + \cosh(x)$ in the domain $(-1,1)$ with initial condition 
           $u_0(x) = \sin(\pi x)$ and boundary data given by the analytical
           solution \eqref{eq:var-lin-adv-without-boundary-solution}.
          }
  \label{tab:convergence}
  \begin{tabular}{rrr | cccccccc}
    \toprule
    &&& \multicolumn{4}{c}{Split form \eqref{eq:standard-general}}
      & \multicolumn{4}{c}{Unsplit form \eqref{eq:weighted-conservative-semidiscretisation}}
    \\
    &&& \multicolumn{2}{c}{Central flux} & \multicolumn{2}{c}{Upwind flux}
      & \multicolumn{2}{c}{Central flux} & \multicolumn{2}{c}{Upwind flux}
    \\
    & $p$  & $N$
    & $\norm{u - u_\mathrm{ana}}$ & EOC & $\norm{u - u_\mathrm{ana}}$ & EOC
    & $\norm{u - u_\mathrm{ana}}$ & EOC & $\norm{u - u_\mathrm{ana}}$ & EOC
    \\
    \midrule 
    \multirow{12}{*}{\rotatebox{90}{Lobatto nodes (Lobatto for $a$)}} 
    & $ 5$ & $   8$  & 4.05e-02 &        & 4.06e-02 &        & 4.24e-02 &        & 4.25e-02 &       \\ 
    &      & $  16$  & 1.16e-03 &  5.13  & 1.18e-03 &  5.10  & 1.08e-03 &  5.29  & 1.11e-03 &  5.26 \\ 
    &      & $  32$  & 2.15e-04 &  2.43  & 2.25e-04 &  2.39  & 2.15e-04 &  2.34  & 2.25e-04 &  2.29 \\ 
    &      & $  64$  & 8.76e-06 &  4.62  & 8.80e-06 &  4.68  & 8.83e-06 &  4.60  & 8.90e-06 &  4.66 \\ 
    &      & $ 128$  & 2.42e-07 &  5.17  & 1.94e-07 &  5.51  & 2.45e-07 &  5.17  & 1.96e-07 &  5.50 \\ 
    &      & $ 256$  & 6.82e-09 &  5.15  & 3.41e-09 &  5.83  & 6.88e-09 &  5.15  & 3.46e-09 &  5.83 \\ 
    \cmidrule{3-11} 
    & $ 6$ & $   8$  & 4.32e-03 &        & 4.32e-03 &        & 4.73e-03 &        & 4.73e-03 &       \\ 
    &      & $  16$  & 8.73e-04 &  2.31  & 8.66e-04 &  2.32  & 8.71e-04 &  2.44  & 8.64e-04 &  2.45 \\ 
    &      & $  32$  & 4.16e-05 &  4.39  & 4.13e-05 &  4.39  & 4.22e-05 &  4.37  & 4.19e-05 &  4.37 \\ 
    &      & $  64$  & 6.88e-07 &  5.92  & 7.27e-07 &  5.83  & 7.02e-07 &  5.91  & 7.40e-07 &  5.82 \\ 
    &      & $ 128$  & 6.52e-09 &  6.72  & 7.70e-09 &  6.56  & 6.66e-09 &  6.72  & 7.87e-09 &  6.56 \\ 
    &      & $ 256$  & 6.35e-11 &  6.68  & 7.59e-11 &  6.66  & 6.45e-11 &  6.69  & 7.72e-11 &  6.67 \\ 
    \midrule 
    \multirow{12}{*}{\rotatebox{90}{Gauss nodes (Lobatto for $a$)}} 
    & $ 5$ & $   8$  & 1.36e-02 &        & 1.35e-02 &        & 1.53e-02 &        & 1.53e-02 &       \\ 
    &      & $  16$  & 5.30e-04 &  4.68  & 4.89e-04 &  4.79  & 5.77e-04 &  4.73  & 5.27e-04 &  4.86 \\ 
    &      & $  32$  & 4.48e-05 &  3.56  & 5.34e-05 &  3.19  & 4.51e-05 &  3.68  & 5.46e-05 &  3.27 \\ 
    &      & $  64$  & 2.63e-06 &  4.09  & 2.44e-06 &  4.45  & 2.68e-06 &  4.07  & 2.52e-06 &  4.44 \\ 
    &      & $ 128$  & 9.71e-08 &  4.76  & 5.62e-08 &  5.44  & 9.88e-08 &  4.76  & 5.82e-08 &  5.44 \\ 
    &      & $ 256$  & 3.16e-09 &  4.94  & 1.01e-09 &  5.80  & 3.22e-09 &  4.94  & 1.04e-09 &  5.80 \\ 
    \cmidrule{3-11} 
    & $ 6$ & $   8$  & 3.56e-03 &        & 3.55e-03 &        & 3.93e-03 &        & 3.93e-03 &       \\ 
    &      & $  16$  & 1.67e-04 &  4.42  & 1.60e-04 &  4.47  & 1.72e-04 &  4.51  & 1.65e-04 &  4.57 \\ 
    &      & $  32$  & 1.06e-05 &  3.98  & 1.06e-05 &  3.92  & 1.11e-05 &  3.95  & 1.11e-05 &  3.90 \\ 
    &      & $  64$  & 1.63e-07 &  6.02  & 2.08e-07 &  5.67  & 1.72e-07 &  6.02  & 2.18e-07 &  5.67 \\ 
    &      & $ 128$  & 1.07e-09 &  7.25  & 2.29e-09 &  6.50  & 1.13e-09 &  7.25  & 2.40e-09 &  6.50 \\ 
    &      & $ 256$  & 3.65e-11 &  4.88  & 4.14e-11 &  5.79  & 3.65e-11 &  4.95  & 4.19e-11 &  5.84 \\ 
    \midrule 
    \multirow{12}{*}{\rotatebox{90}{Gauss nodes (Gauss for $a$)}} 
    & $ 5$ & $   8$  & 1.36e-02 &        & 1.35e-02 &        & 1.53e-02 &        & 1.53e-02 &       \\ 
    &      & $  16$  & 5.30e-04 &  4.68  & 4.89e-04 &  4.79  & 5.77e-04 &  4.73  & 5.27e-04 &  4.86 \\ 
    &      & $  32$  & 4.48e-05 &  3.56  & 5.34e-05 &  3.19  & 4.51e-05 &  3.68  & 5.46e-05 &  3.27 \\ 
    &      & $  64$  & 2.63e-06 &  4.09  & 2.44e-06 &  4.45  & 2.68e-06 &  4.07  & 2.52e-06 &  4.44 \\ 
    &      & $ 128$  & 9.71e-08 &  4.76  & 5.62e-08 &  5.44  & 9.88e-08 &  4.76  & 5.82e-08 &  5.44 \\ 
    &      & $ 256$  & 3.16e-09 &  4.94  & 1.01e-09 &  5.80  & 3.22e-09 &  4.94  & 1.04e-09 &  5.80 \\ 
    \cmidrule{3-11} 
    & $ 6$ & $   8$  & 3.56e-03 &        & 3.55e-03 &        & 3.93e-03 &        & 3.93e-03 &       \\ 
    &      & $  16$  & 1.67e-04 &  4.42  & 1.60e-04 &  4.47  & 1.72e-04 &  4.51  & 1.65e-04 &  4.57 \\ 
    &      & $  32$  & 1.06e-05 &  3.98  & 1.06e-05 &  3.92  & 1.11e-05 &  3.95  & 1.11e-05 &  3.90 \\ 
    &      & $  64$  & 1.63e-07 &  6.02  & 2.08e-07 &  5.67  & 1.72e-07 &  6.02  & 2.18e-07 &  5.67 \\ 
    &      & $ 128$  & 1.07e-09 &  7.25  & 2.29e-09 &  6.50  & 1.13e-09 &  7.25  & 2.40e-09 &  6.50 \\ 
    &      & $ 256$  & 3.65e-11 &  4.88  & 4.14e-11 &  5.79  & 3.65e-11 &  4.95  & 4.19e-11 &  5.84 \\ 
    \bottomrule 
  \end{tabular}
\end{table}

Numerical results of these convergence tests using polynomials of degree $p=5$
and $p=6$ are listed in \autoref{tab:convergence}. 
There, the errors $\norm{u - u_\mathrm{ana}}$ have been calculated using Gauss
quadrature on $p+1$ nodes. The experimental order of convergence (EOC) has been
calculated as 
\begin{equation}
  \mathrm{EOC}(N_1,N_2)
  =
  - \frac{ \log\big( \mathrm{error}(N_1) / \mathrm{error}(N_2) \big) }
         { \log\big( N_1 / N_2 \big) }.
\end{equation}
As can be seen in \autoref{tab:convergence}, all schemes seem to converge 
numerically to the analytical solution. The experimental orders of convergence
show some variations. In most cases, they are in the range $[4, 6]$ for $p = 5$
and in $[4, 7]$ for $p = 6$.

There are only small deviations (up to approximately \si{10}{\%}) between the 
errors of the split forms and the corresponding unsplit forms. Similarly, the
errors for the central fluxes are of the same order as the errors for the upwind
fluxes. However, there are deviations by a factor up to approximately $3$,
e.g. for Gauss nodes (Lobatto for $a$), $p=5$, $N=256$, unsplit form.
Especially for the odd polynomial degree $p=5$, the upwind flux seems to result
in slightly more accurate schemes. This can be regarded as related to the results
of \cite{nordstrom2007error, kopriva2017error}, where the error of numerical
solutions of linear hyperbolic problems has been investigated and related to
dissipation by upwind numerical fluxes.

In general, Gauss nodes yield smaller errors than Lobatto nodes. In most cases,
the error using Gauss nodes is smaller by a factor of $2$ to $4$. Contrary, the 
choice of the basis used to compute the discretised speed~$a$ does not seem
to influence the errors in \autoref{tab:convergence}. However, there are some
slight differences of order $10^{-12}$ (not visible in \autoref{tab:convergence}).

\begin{figure}
\centering
  \begin{subfigure}{0.49\textwidth}
    \includegraphics[width=\textwidth]{%
    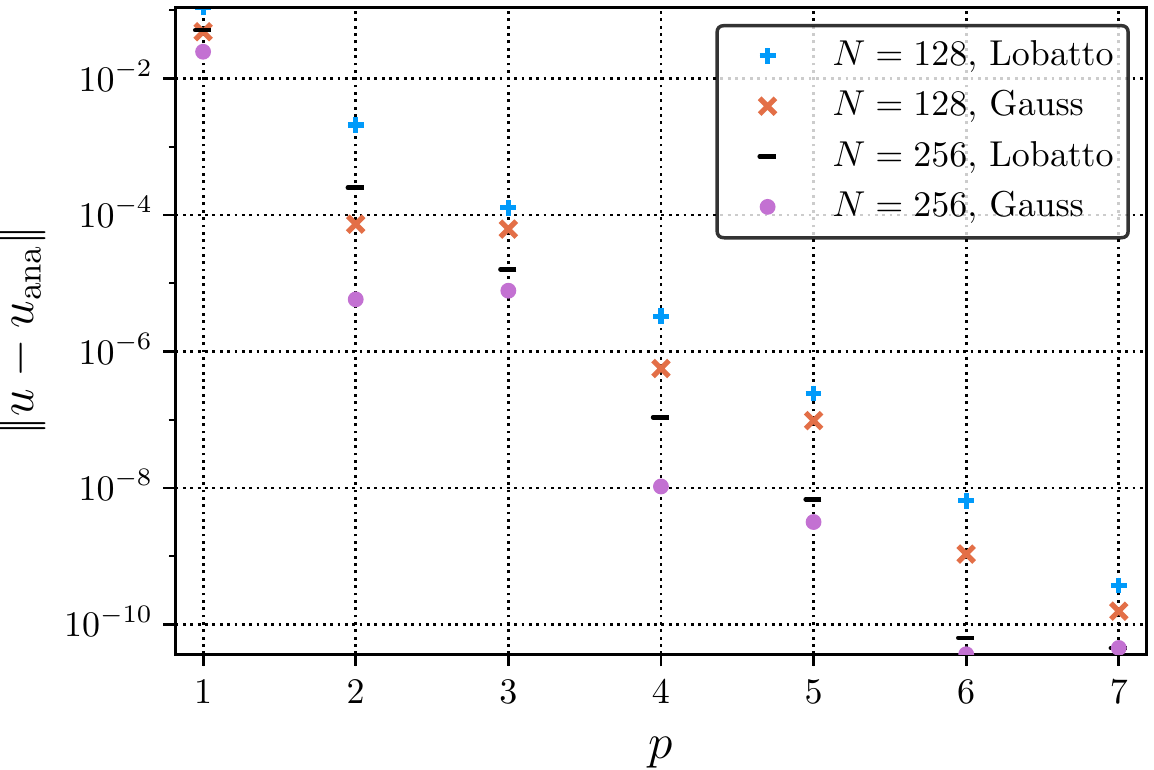}
    \caption{Central flux.}
    \label{fig:conservative_nonperiodic__1_p_cosh_x__sinpi_x__exponential_convergence_central}
  \end{subfigure}%
  \begin{subfigure}{0.49\textwidth}
    \includegraphics[width=\textwidth]{%
    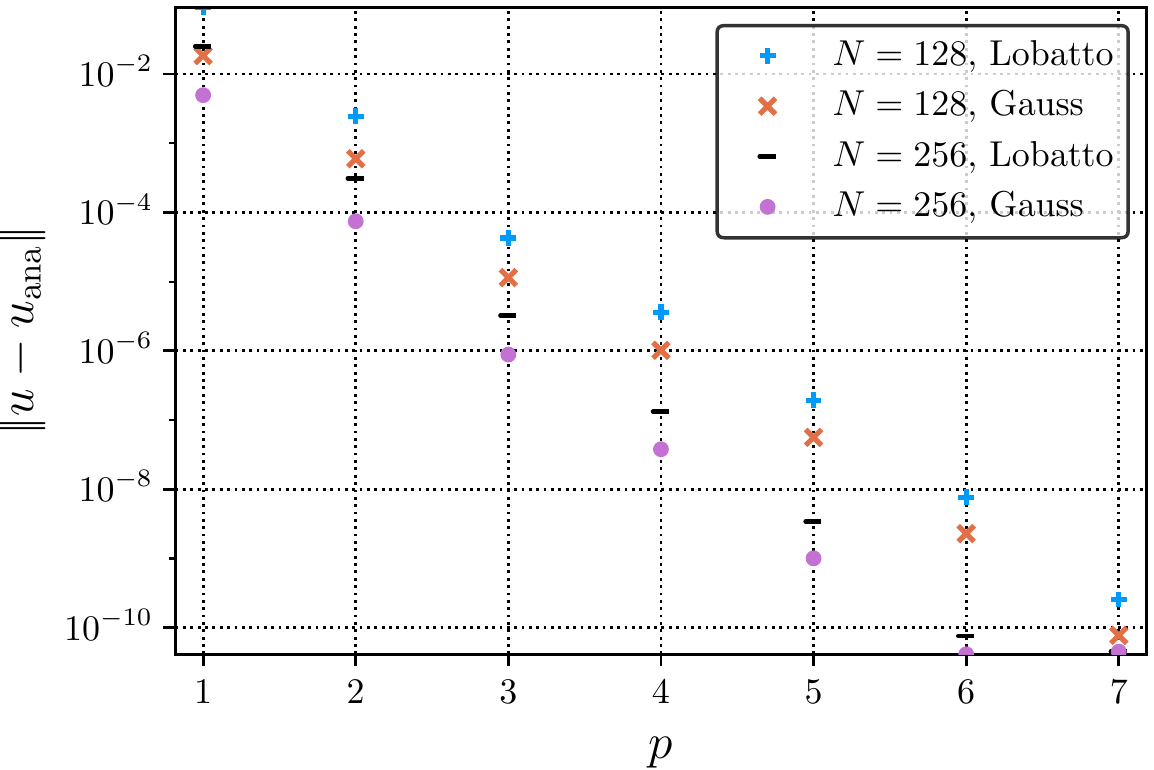}
    \caption{Upwind flux.}
    \label{fig:conservative_nonperiodic__1_p_cosh_x__sinpi_x__exponential_convergence_upwind}
  \end{subfigure}%
  \caption{Convergence results for \eqref{eq:var-lin-adv} using the speed 
           $a(x) = 1 + \cosh(x)$ in the domain $(-1,1)$ with initial condition 
           $u_0(x) = \sin(\pi x)$ and boundary data given by the analytical
           solution \eqref{eq:var-lin-adv-without-boundary-solution}.
           The semidiscretisation uses $N$ elements with polynomials of degree~$p$
           and the split form \eqref{eq:standard-general}.}
  \label{fig:conservative_nonperiodic__1_p_cosh_x__sinpi_x__exponential_convergence}
\end{figure}

Considering varying polynomial degrees $p$ for fixed number of elements $N$,
exponential convergence rates are obtained, as visualised in 
\autoref{fig:conservative_nonperiodic__1_p_cosh_x__sinpi_x__exponential_convergence}.
There, the split form \eqref{eq:standard-general} has been used with both the 
central flux \eqref{eq:fnum-split-central} and the upwind flux \eqref{eq:fnum-split-upwind}
at interior boundaries.

\subsection{Conservation Properties}
\label{subsec:conservation-properties}

Here, the linear advection equation \eqref{eq:var-lin-adv-without-boundary} with
speed $a(x) = \cos\big( \frac{\pi}{2} x \big)$ is considered in the domain $(-1,1)$
with initial condition $u_0(x) = 1 + \frac{1}{2} \cos(\pi x)$. Since the speed
vanishes at the boundaries, zero boundary data are prescribed.

Although only positive speeds $a(x) > 0$ have been considered, all computations
in the standard $L^2$ setting of \autoref{sec:standard-L2} are valid for nonnegative
speed $a(x) \geq 0$. 
However, the discrete weighted ``norm'' $\norm{ \cdot }_{\mat{a}[^T]\mat{M}}$
of \autoref{sec:weighted-L2} is not positive definite if boundary nodes are included.

\begin{table}
\centering
  \caption{Conservation results for \eqref{eq:var-lin-adv} using the speed 
           $a(x) = \cos\big( \frac{\pi}{2} x \big)$ in the domain $(-1,1)$ with 
           initial condition $u_0(x) = 1 + \frac{1}{2} \cos(\pi x)$ and zero
           boundary data.
          }
  \label{tab:conservation}
  \begin{tabular}{rrr | cccccccc}
    \toprule
    &&& \multicolumn{4}{c}{Split form \eqref{eq:standard-general}}
      & \multicolumn{4}{c}{Unsplit form \eqref{eq:weighted-conservative-semidiscretisation}}
    \\
    &&& \multicolumn{2}{c}{Central flux} & \multicolumn{2}{c}{Upwind flux}
      & \multicolumn{2}{c}{Central flux} & \multicolumn{2}{c}{Upwind flux}
    \\
    & $p$  & $N$
    & $\abs{ \int u - \int u_0 }$ & EOC & $\abs{ \int u - \int u_0 }$ & EOC
    & $\abs{ \int u - \int u_0 }$ & EOC & $\abs{ \int u - \int u_0 }$ & EOC
    \\
    \midrule 
    \multirow{12}{*}{\rotatebox{90}{Lobatto nodes (Lobatto for $a$)}} 
    & $ 3$ & $   8$  & 5.62e-15 &        & 1.50e-15 &        & 1.23e-15 &        & 2.54e-15 &       \\ 
    &      & $  16$  & 3.57e-15 &        & 6.01e-15 &        & 5.43e-15 &        & 2.27e-16 &       \\ 
    &      & $  32$  & 4.03e-15 &        & 6.94e-16 &        & 1.43e-15 &        & 5.32e-16 &       \\ 
    &      & $  64$  & 2.87e-15 &        & 1.19e-15 &        & 4.27e-15 &        & 3.02e-15 &       \\ 
    &      & $ 128$  & 1.17e-16 &        & 4.41e-15 &        & 3.60e-15 &        & 8.55e-16 &       \\ 
    &      & $ 256$  & 6.10e-16 &        & 3.56e-15 &        & 4.55e-15 &        & 4.75e-15 &       \\ 
    \cmidrule{3-11} 
    & $ 4$ & $   8$  & 4.28e-15 &        & 1.92e-15 &        & 3.24e-15 &        & 1.33e-15 &       \\ 
    &      & $  16$  & 2.47e-15 &        & 1.46e-15 &        & 1.51e-15 &        & 1.17e-15 &       \\ 
    &      & $  32$  & 1.11e-16 &        & 4.52e-15 &        & 4.35e-16 &        & 4.45e-15 &       \\ 
    &      & $  64$  & 1.05e-15 &        & 1.05e-15 &        & 4.88e-16 &        & 5.80e-16 &       \\ 
    &      & $ 128$  & 5.98e-15 &        & 1.58e-15 &        & 7.47e-15 &        & 3.30e-15 &       \\ 
    &      & $ 256$  & 6.99e-15 &        & 5.00e-15 &        & 7.22e-15 &        & 9.44e-15 &       \\ 
    \midrule 
    \multirow{12}{*}{\rotatebox{90}{Gauss nodes (Lobatto for $a$)}} 
    & $ 3$ & $   8$  & 5.41e-16 &        & 5.31e-15 &        & 5.34e-04 &        & 5.36e-04 &       \\ 
    &      & $  16$  & 3.91e-15 &        & 7.92e-15 &        & 2.26e-05 &  4.56  & 2.27e-05 &  4.56 \\ 
    &      & $  32$  & 9.55e-15 &        & 6.28e-15 &        & 7.58e-07 &  4.90  & 7.63e-07 &  4.90 \\ 
    &      & $  64$  & 1.45e-14 &        & 1.76e-14 &        & 2.41e-08 &  4.98  & 2.42e-08 &  4.98 \\ 
    &      & $ 128$  & 2.41e-14 &        & 2.12e-14 &        & 7.55e-10 &  4.99  & 7.60e-10 &  4.99 \\ 
    &      & $ 256$  & 5.04e-14 &        & 5.06e-14 &        & 2.36e-11 &  5.00  & 2.37e-11 &  5.00 \\ 
    \cmidrule{3-11} 
    & $ 4$ & $   8$  & 1.36e-15 &        & 7.15e-16 &        & 2.59e-05 &        & 2.59e-05 &       \\ 
    &      & $  16$  & 2.68e-15 &        & 4.70e-16 &        & 2.48e-06 &  3.39  & 2.48e-06 &  3.38 \\ 
    &      & $  32$  & 2.60e-15 &        & 2.29e-15 &        & 9.56e-08 &  4.70  & 9.56e-08 &  4.70 \\ 
    &      & $  64$  & 4.89e-16 &        & 4.95e-15 &        & 3.14e-09 &  4.93  & 3.14e-09 &  4.93 \\ 
    &      & $ 128$  & 5.67e-15 &        & 4.39e-15 &        & 9.94e-11 &  4.98  & 9.95e-11 &  4.98 \\ 
    &      & $ 256$  & 1.42e-14 &        & 1.40e-14 &        & 3.12e-12 &  4.99  & 3.12e-12 &  4.99 \\ 
    \midrule 
    \multirow{12}{*}{\rotatebox{90}{Gauss nodes (Gauss for $a$)}} 
    & $ 3$ & $   8$  & 8.41e-07 &        & 8.41e-07 &        & 5.34e-04 &        & 5.35e-04 &       \\ 
    &      & $  16$  & 2.64e-08 &  5.00  & 2.64e-08 &  5.00  & 2.26e-05 &  4.56  & 2.27e-05 &  4.56 \\ 
    &      & $  32$  & 8.25e-10 &  5.00  & 8.25e-10 &  5.00  & 7.57e-07 &  4.90  & 7.62e-07 &  4.90 \\ 
    &      & $  64$  & 2.58e-11 &  5.00  & 2.58e-11 &  5.00  & 2.40e-08 &  4.98  & 2.42e-08 &  4.98 \\ 
    &      & $ 128$  & 8.32e-13 &  4.95  & 8.33e-13 &  4.95  & 7.54e-10 &  4.99  & 7.59e-10 &  4.99 \\ 
    &      & $ 256$  & 8.04e-14 &  3.37  & 8.84e-14 &  3.24  & 2.35e-11 &  5.00  & 2.37e-11 &  5.00 \\ 
    \cmidrule{3-11} 
    & $ 4$ & $   8$  & 1.16e-07 &        & 1.16e-07 &        & 2.58e-05 &        & 2.58e-05 &       \\ 
    &      & $  16$  & 3.65e-09 &  4.98  & 3.65e-09 &  4.98  & 2.48e-06 &  3.38  & 2.48e-06 &  3.38 \\ 
    &      & $  32$  & 1.14e-10 &  5.00  & 1.14e-10 &  5.00  & 9.55e-08 &  4.70  & 9.55e-08 &  4.70 \\ 
    &      & $  64$  & 3.58e-12 &  5.00  & 3.58e-12 &  5.00  & 3.14e-09 &  4.93  & 3.14e-09 &  4.93 \\ 
    &      & $ 128$  & 1.16e-13 &  4.94  & 1.10e-13 &  5.02  & 9.93e-11 &  4.98  & 9.94e-11 &  4.98 \\ 
    &      & $ 256$  & 1.05e-14 &  3.46  & 9.88e-15 &  3.48  & 3.13e-12 &  4.99  & 3.12e-12 &  4.99 \\ 
    \bottomrule 
  \end{tabular}
\end{table}

Using the same time discretisation as in section~\ref{subsec:convergence-studies},
numerical results for polynomials degrees $p \in \set{3,4}$ are presented in
\autoref{tab:conservation}. There, the errors $\abs{\int u - \int u_0 }$ in the
conserved variables are computed using the natural quadrature rule associated
with the mass matrix $\mat{M}$.

As can be seen in \autoref{tab:conservation}, all schemes converge numerically
to the correct result $\int u = \int u_0$. If Lobatto nodes are used, the
conservation error vanishes numerically, since boundary nodes are included and
the numerical fluxes at the exterior boundaries are zero. If the split form
\eqref{eq:standard-general} is used with Gauss nodes and $a$ is discretised via
Lobatto nodes, the interpolated values of $a$ are exact at the boundaries and
the conservation error vanishes again.

However, using the split form \eqref{eq:standard-general} with Gauss nodes only,
the interpolated speed will in general not vanish at the boundaries. Nevertheless,
the error goes to zero with experimental order of convergence approximately $5$
for both $p = 3$ and $p = 4$.

Similarly, since unsplit fluxes are used for the unsplit form 
\eqref{eq:weighted-conservative-semidiscretisation}, the interpolated values of
$a u$ will in general not vanish at the boundaries if Gauss nodes are used,
independently of the discretisation of the speed~$a$, see also
Remark~\ref{rem:conservation-in-thm:weighted-conservative}. Thus, there are some
conservation errors, vanishing at approximately the same rate as those of the
split form \eqref{eq:standard-general}. However, the conservation errors of the
split form \eqref{eq:standard-general} are approximately three orders of magnitude
smaller than the ones of the unsplit form 
\eqref{eq:weighted-conservative-semidiscretisation}.

\subsection{Eigenvalues}
\label{subsec:eigenvalues}

Applying a semidiscretisation to the linear advection equation \eqref{eq:var-lin-adv}
results in a system of ordinary differential equations of the form
\begin{equation}
  \od{}{t} u(t) = L \cdot u(t) + b(t),
\end{equation}
where $L$ is a linear operator and $b(t)$ describes the boundary condition.
Here, the speed $a(x) = 1 + \cosh(x)$ and the domain $(\xmin,\xmax) = (-1,1)$ are
the same as in the convergence studies in section~\ref{subsec:convergence-studies}.

Again, both the split form \eqref{eq:standard-general} and the unsplit form
\eqref{eq:weighted-conservative-semidiscretisation} are considered, corresponding
central and numerical fluxes are used and both Lobatto and Gauss nodes are
investigated. For Gauss nodes, the speed~$a$ is discretised via Lobatto nodes
as described in Remark~\ref{rem:use-Lobatto-nodes-and-interpolate-to-Gauss} if
the split form is used. Otherwise, $a$ is interpolated on Gauss nodes.

\begin{figure}
\centering
  \begin{subfigure}{0.49\textwidth}
    \includegraphics[width=\textwidth]{%
    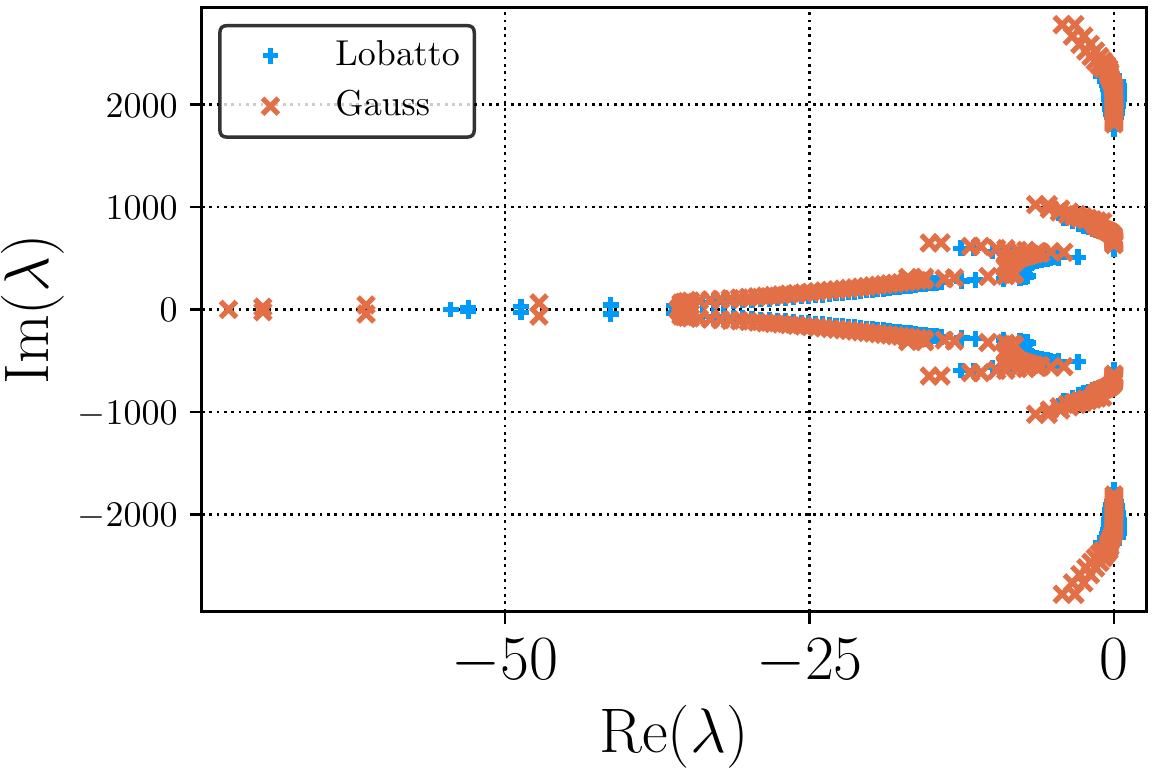}
    \caption{Split form \eqref{eq:standard-general}, 
             central flux.}
    \label{fig:conservative_nonperiodic__1_p_cosh_x__eigenvalues_split_central}
  \end{subfigure}%
  \begin{subfigure}{0.49\textwidth}
    \includegraphics[width=\textwidth]{%
    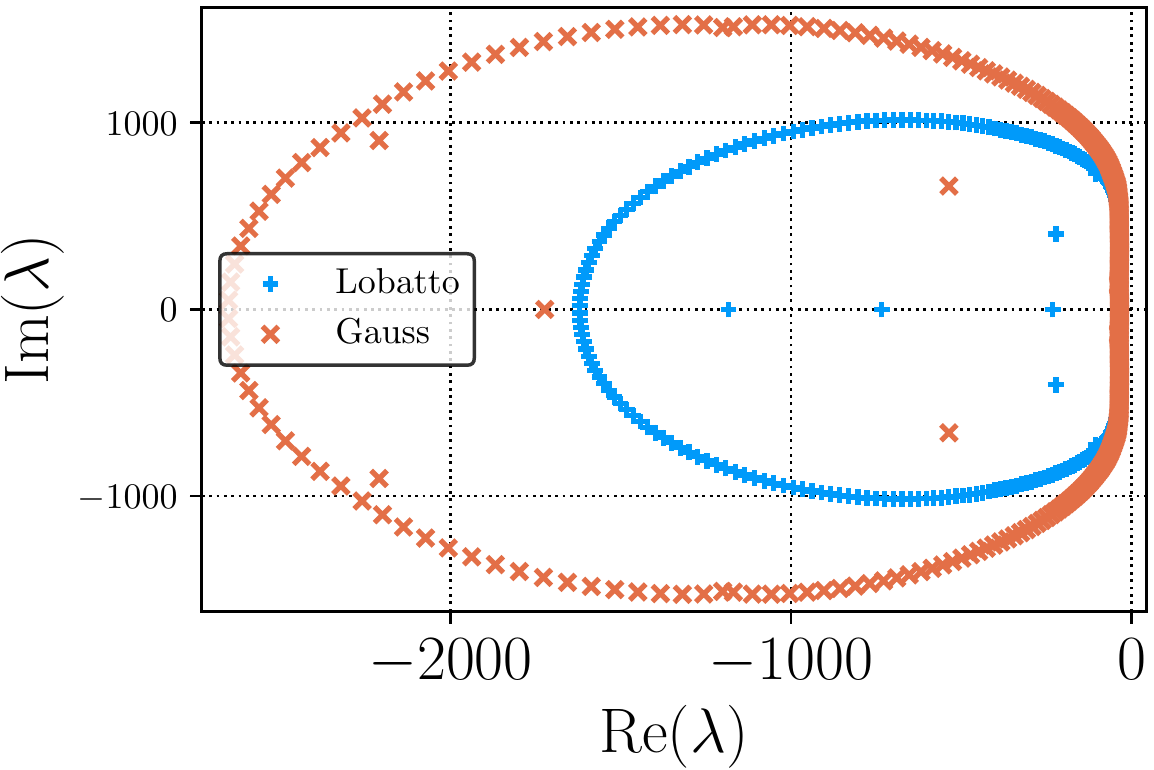}
    \caption{Split form \eqref{eq:standard-general},
             upwind flux.}
    \label{fig:conservative_nonperiodic__1_p_cosh_x__eigenvalues_split_upwind}
  \end{subfigure}%
  \\
  \begin{subfigure}{0.49\textwidth}
    \includegraphics[width=\textwidth]{%
    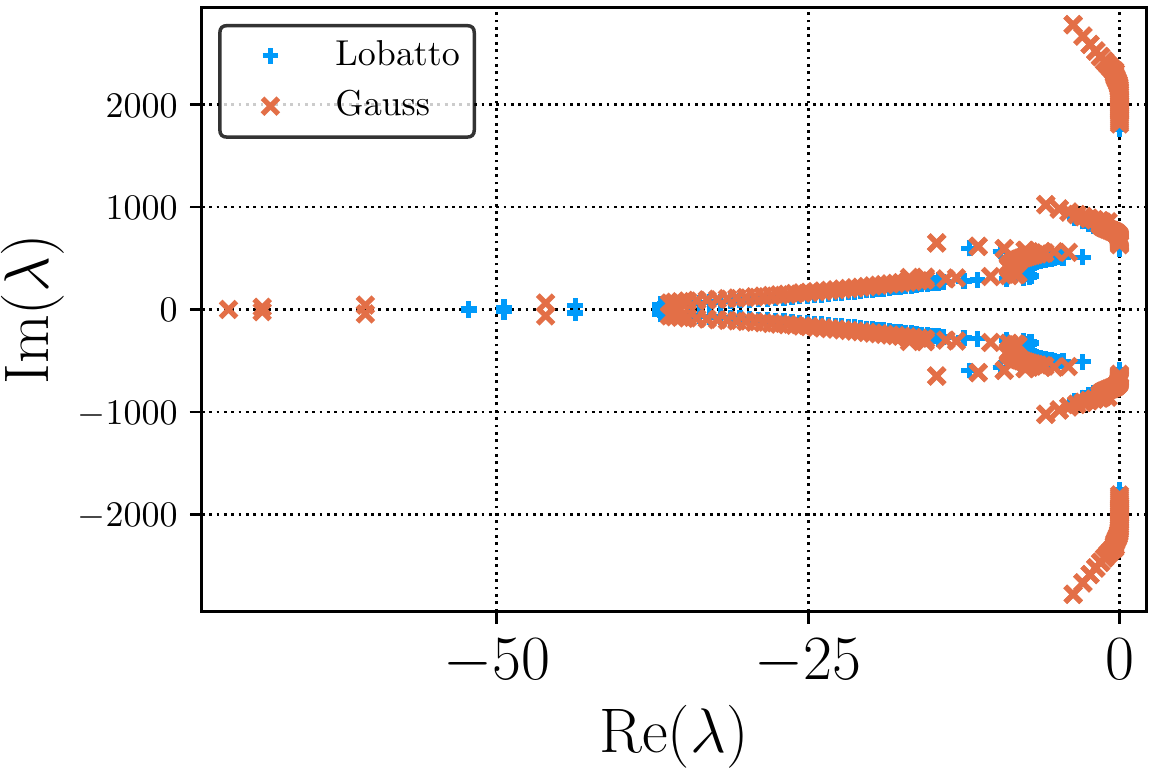}
    \caption{Unsplit form \eqref{eq:weighted-conservative-semidiscretisation},
             central flux.}
    \label{fig:conservative_nonperiodic__1_p_cosh_x__eigenvalues_unsplit_central}
  \end{subfigure}%
  \begin{subfigure}{0.49\textwidth}
    \includegraphics[width=\textwidth]{%
    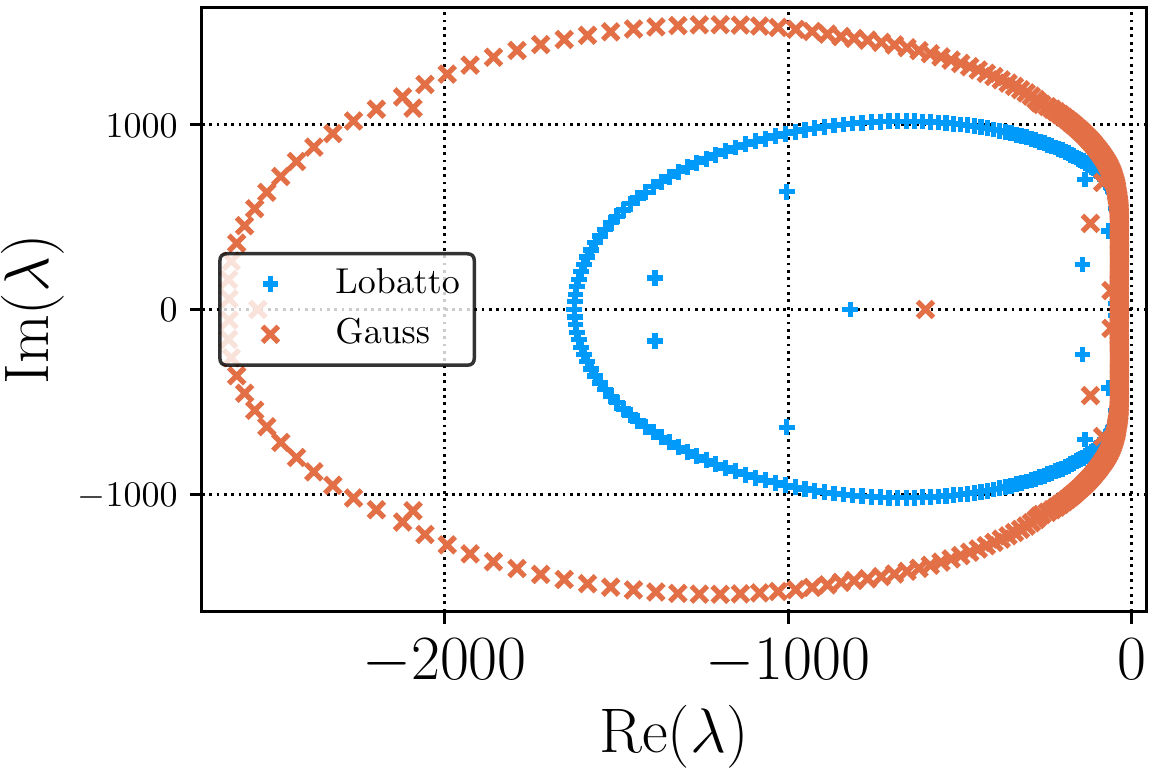}
    \caption{Unsplit form \eqref{eq:weighted-conservative-semidiscretisation},
             upwind flux.}
    \label{fig:conservative_nonperiodic__1_p_cosh_x__eigenvalues_unsplit_upwind}
  \end{subfigure}%
  \caption{Eigenvalues of linear operators of semidiscretisations of the nonperiodic
           advection equation \eqref{eq:var-lin-adv} using the speed $a(x) = 1 + \cosh(x)$
           in the domain $(-1,1)$. The semidiscretisations use $N = 50$ elements
           with polynomials of degree~$p = 7$.}
  \label{fig:conservative_nonperiodic__1_p_cosh_x__eigenvalues}
\end{figure}

The eigenvalues of the resulting linear operators $L$ are visualised in 
\autoref{fig:conservative_nonperiodic__1_p_cosh_x__eigenvalues}. In these plots,
there are only minor differences between the split form \eqref{eq:standard-general}
and the unsplit form \eqref{eq:weighted-conservative-semidiscretisation}.
The (numerical) spectra using Gauss and Lobatto nodes have similar shapes, but
the magnitude of the eigenvalues is smaller for Lobatto nodes. 
Comparing the eigenvalues with minimal real part, the absolute values for Lobatto 
nodes are approximately $0.6$ to $0.8$ times their counterparts for Gauss nodes.
Considering the greatest absolute value of the eigenvalues, the factor remains
the same. 
Thus, the generalised SBP operators using Gauss nodes are more ``stiff'' than
their counterparts using Lobatto nodes. This result can also be found similarly
in \cite{gassner2011comparison} for a linear problem with constant coefficients.
It can be attributed to the different mass matrices of the SBP operators. Indeed,
considering a change of basis to a modal Legendre basis of polynomials up to
degree $p$, the mass matrix from Gauss nodes is exact, i.e.
$\mat{\hat M}[_\mathrm{Gauss}] = \diag{ 2, \frac{2}{3}, \dots, \frac{2}{2p+1}}$,
whereas the mass matrix for Lobatto nodes is not completely exact; the last
entry differs and is $\frac{2}{p}$ instead of the correct value $\frac{2}{2p+1}$
\cite[equation (1.136)]{kopriva2009implementing}
\cite[equation (2.3.13)]{canuto2006spectral}. Since the inverse of the mass
matrix is used in the semidiscretisation, the difference of the spectra can be
explained. Some investigations in order to reduce the CFL condition for DG methods
can be found in \cite{warburton2008taming}.

\subsubsection*{Periodic Boundary Conditions}

Considering periodic boundary conditions yields a system of ordinary differential 
equations of the form
\begin{equation}
  \od{}{t} u(t) = L \cdot u(t),
\end{equation}
where $L$ is again a linear operator. The periodic boundary conditions are
enforced weakly by coupling of the left boundary of the first element with 
the right boundary of the last element as in the case of the other interior
boundaries.
Here, the speed $a(x) = 2 + \sin(\pi x)$ is the same as in 
section~\ref{subsec:conservation-properties}.

\begin{figure}
\centering
  \begin{subfigure}{0.49\textwidth}
    \includegraphics[width=\textwidth]{%
    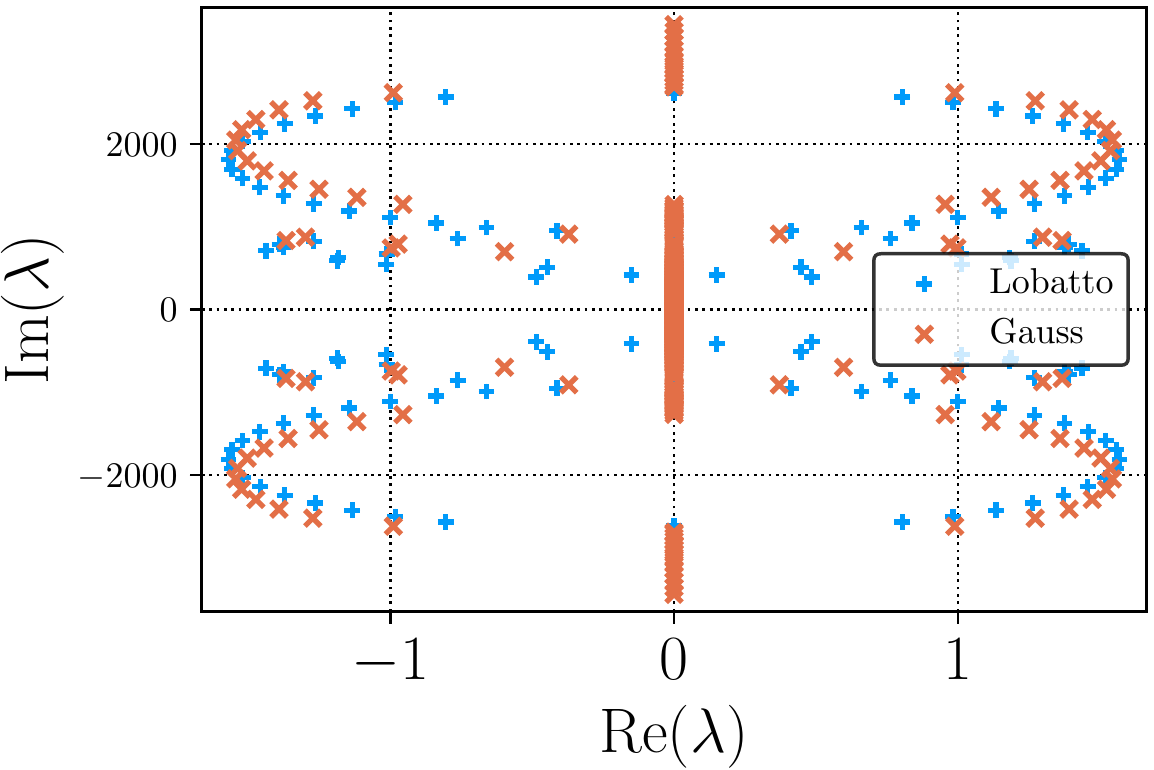}
    \caption{Split form \eqref{eq:standard-general}, 
             central flux.}
    \label{fig:conservative_periodic__2_p_sinpi_x__eigenvalues_split_central}
  \end{subfigure}%
  \begin{subfigure}{0.49\textwidth}
    \includegraphics[width=\textwidth]{%
    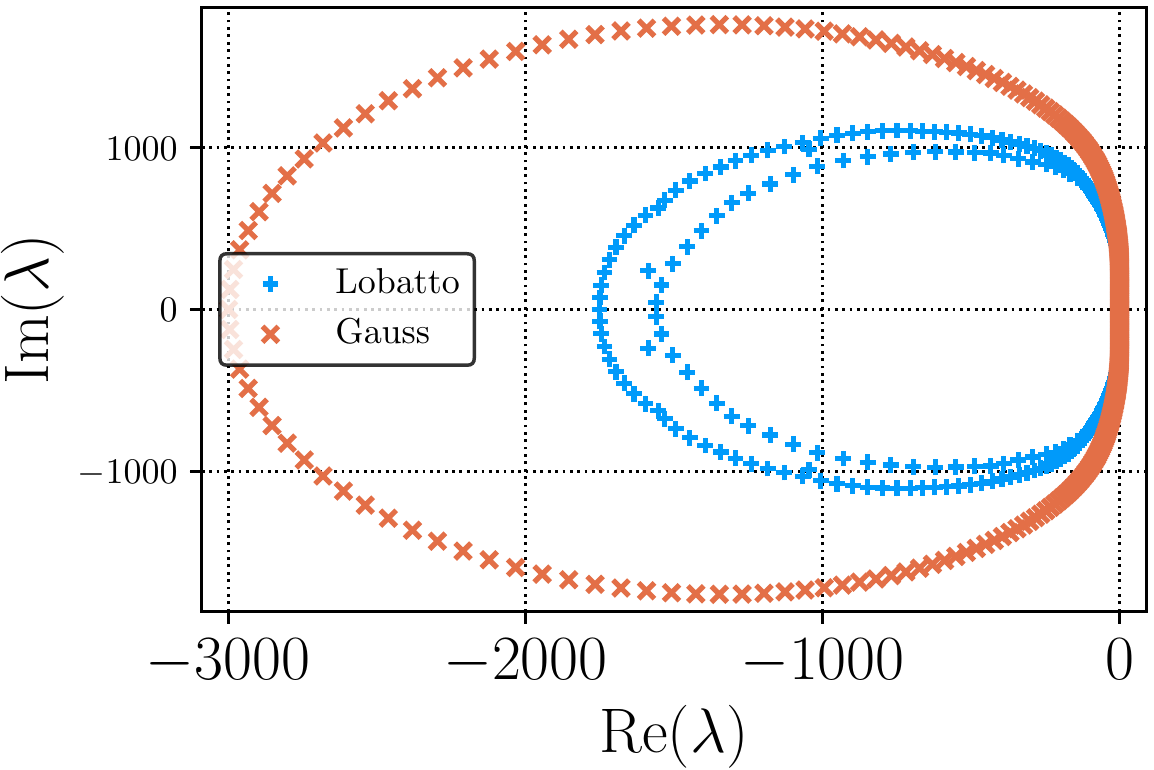}
    \caption{Split form \eqref{eq:standard-general},
             upwind flux.}
    \label{fig:conservative_periodic__2_p_sinpi_x__eigenvalues_split_upwind}
  \end{subfigure}%
  \\
  \begin{subfigure}{0.49\textwidth}
    \includegraphics[width=\textwidth]{%
    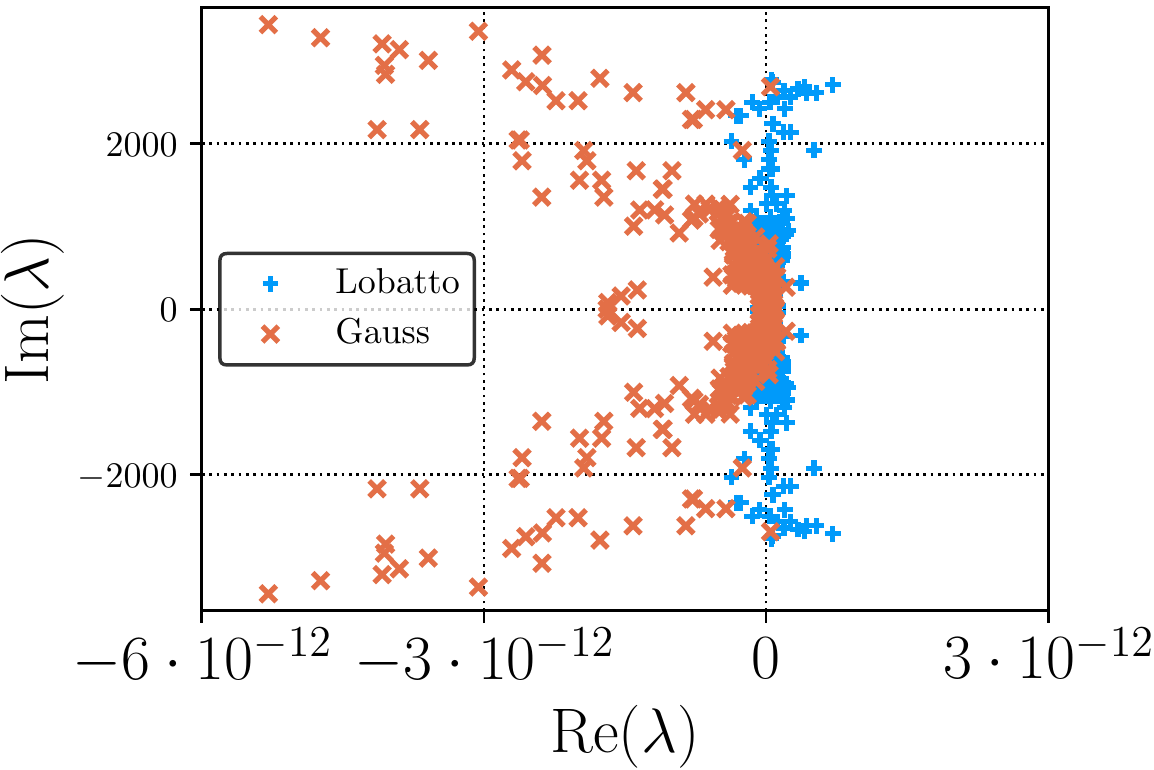}
    \caption{Unsplit form \eqref{eq:weighted-conservative-semidiscretisation},
             central flux.}
    \label{fig:conservative_periodic__2_p_sinpi_x__eigenvalues_unsplit_central}
  \end{subfigure}%
  \begin{subfigure}{0.49\textwidth}
    \includegraphics[width=\textwidth]{%
    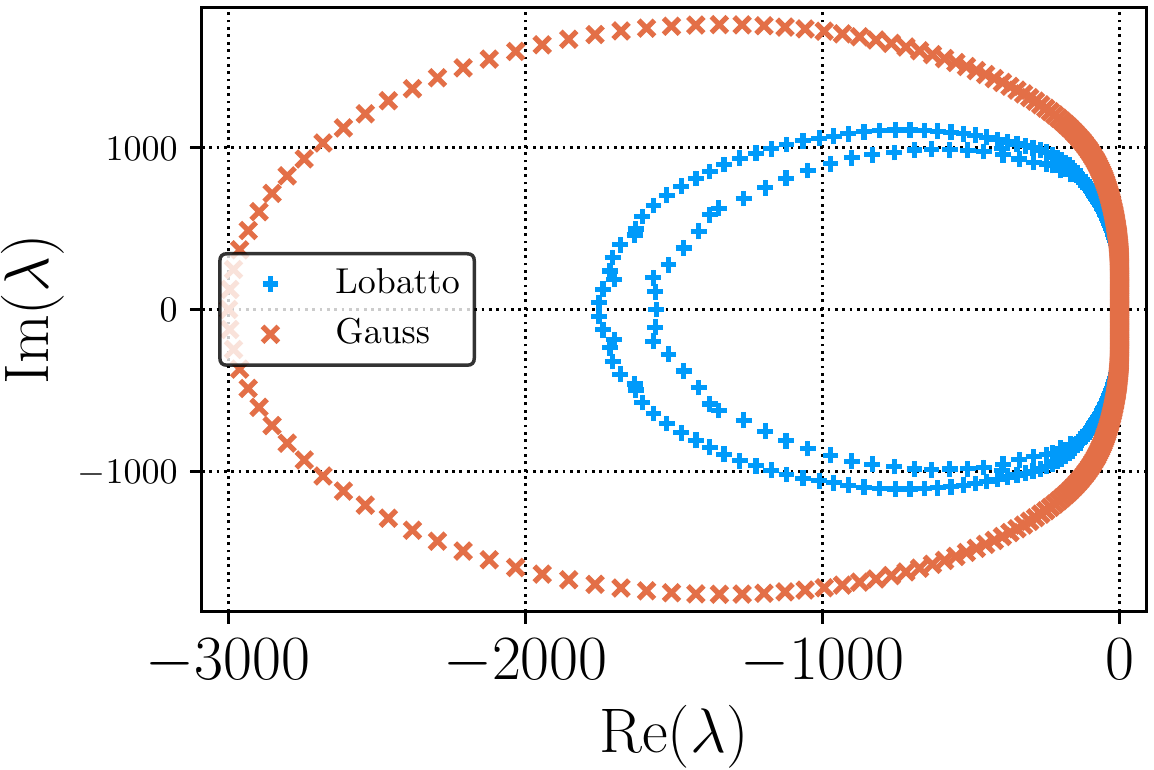}
    \caption{Unsplit form \eqref{eq:weighted-conservative-semidiscretisation},
             upwind flux.}
    \label{fig:conservative_periodic__2_p_sinpi_x__eigenvalues_unsplit_upwind}
  \end{subfigure}%
  \caption{Eigenvalues of linear operators of semidiscretisations of the advection 
           equation with periodic boundary conditions using the speed 
           $a(x) = 2 + \sin(\pi x)$ in the domain $(-1,1)$. The semidiscretisations 
           use $N = 50$ elements with polynomials of degree~$p = 7$.}
  \label{fig:conservative_periodic__2_p_sinpi_x__eigenvalues}
\end{figure}

Using the same numerical parameters as for the nonperiodic equation, the eigenvalues
of the linear operators $L$ for the periodic problem are visualised in
\autoref{fig:conservative_periodic__2_p_sinpi_x__eigenvalues}.
If the upwind flux is used, the spectra for the split form \eqref{eq:standard-general}
are very similar to those of the unsplit form 
\eqref{eq:weighted-conservative-semidiscretisation}. However, using the central
flux results in eigenvalues with positive real part for the split form while
the unsplit form results in eigenvalues with numerically (nearly) vanishing real
part as in the numerical experiments of \cite{manzanero2017insights}, see also
\ref{sec:comparison-manzanero}.

There might seem to be an error, since the split form \eqref{eq:standard-general}
has been proven to be stable in \autoref{thm:standard-general}, but there are
eigenvalues with positive real part in 
\autoref{fig:conservative_periodic__2_p_sinpi_x__eigenvalues_split_central} and
eigenvalues with positive real part are usually identified with unstable schemes.
However, stability of the split forms \eqref{eq:standard-simplified} and
\eqref{eq:standard-general} as stated in Theorems~\ref{thm:standard-simplified}
and~\ref{thm:standard-general} refers to the corresponding energy rate of change
\eqref{eq:standard-stability} and the energy estimate
\eqref{eq:standard-continuous-energy-estimate}. In a periodic setting, the boundary
terms can be dropped. However, due to the term $\int_\xmin^\xmax u^2 \, \partial_x a$
in the rate of change of the energy, the energy estimate contains the factor
$\exp\left(t \, \norm{\partial_x a}_{L^\infty} \right)$, allowing an increase
of the energy. This corresponds to eigenvalues with positive real part in
\autoref{fig:conservative_periodic__2_p_sinpi_x__eigenvalues_split_central}.

Contrary, the stability of the unsplit form
\eqref{eq:weighted-conservative-semidiscretisation} proven in
\autoref{thm:weighted-conservative} is not the same kind of stability. Indeed,
another kind of stability is considered, i.e. not the classical $L^2$ norm, but
the $a$-weighted $L^2$ norm. Considering this kind of stability, no additional
term allowing an increase of the energy appears in the corresponding energy rate
\eqref{eq:weighted-conservative-stability} and energy estimate
\eqref{eq:weighted-conservative-energy-estimate}.%

\subsection{CFL Condition}
\label{subsec:cfl}

Here, time step restrictions of the proposed semidiscretisations will be investigated
numerically. Again, the linear advection equation \eqref{eq:var-lin-adv} with
speed $a(x) = 1 + \cosh(x)$ is considered in the domain $(-1,1)$ as in the previous
sections~\ref{subsec:convergence-studies} and \ref{subsec:eigenvalues}. Using the
ten stage, fourth order strong stability preserving method of \cite{ketcheson2008highly}
as time integrator, different time steps $\Delta t$ are investigated.

Having in mind the CFL condition $\abs{c} \frac{\Delta t}{\Delta x} \leq \frac{1}{2p+1}$
given by \citet[Section~2.2]{cockburn2001runge} for DG methods using polynomials
of degree $p$ and a $p$th order explicit Runge-Kutta method, time steps proportional to
$\frac{1}{2p+1}$ have been used. Here, the width of one element is $\Delta x = \frac{2}{N}$
and the maximal speed is $\abs{c} = \max_{x\in[-1,1]} a(x) = 1 + \cosh(1) \approx 2.5$.
Thus, it can be expected that a time step $\Delta t \leq \frac{1}{(2p+1) N}$ results
in a stable scheme. Of course, the exact restrictions depend on the concrete
discretisations, i.e. on the choice of the (split or unsplit) form, the nodes,
and the numerical fluxes.

\begin{figure}
\centering
  \begin{subfigure}{0.49\textwidth}
    \includegraphics[width=\textwidth]{%
    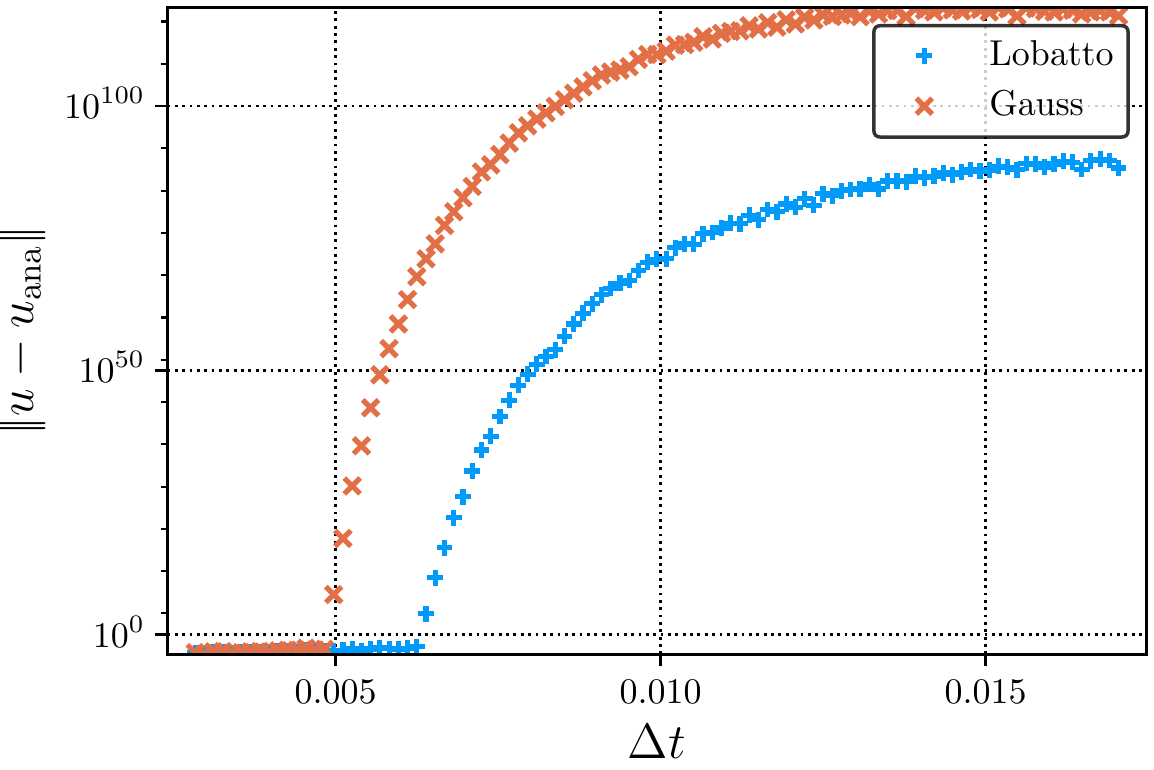}
    \caption{Central flux.}
    \label{fig:conservative_nonperiodic__1_p_cosh_x__sinpi_x__CFL_central}
  \end{subfigure}%
  \begin{subfigure}{0.49\textwidth}
    \includegraphics[width=\textwidth]{%
    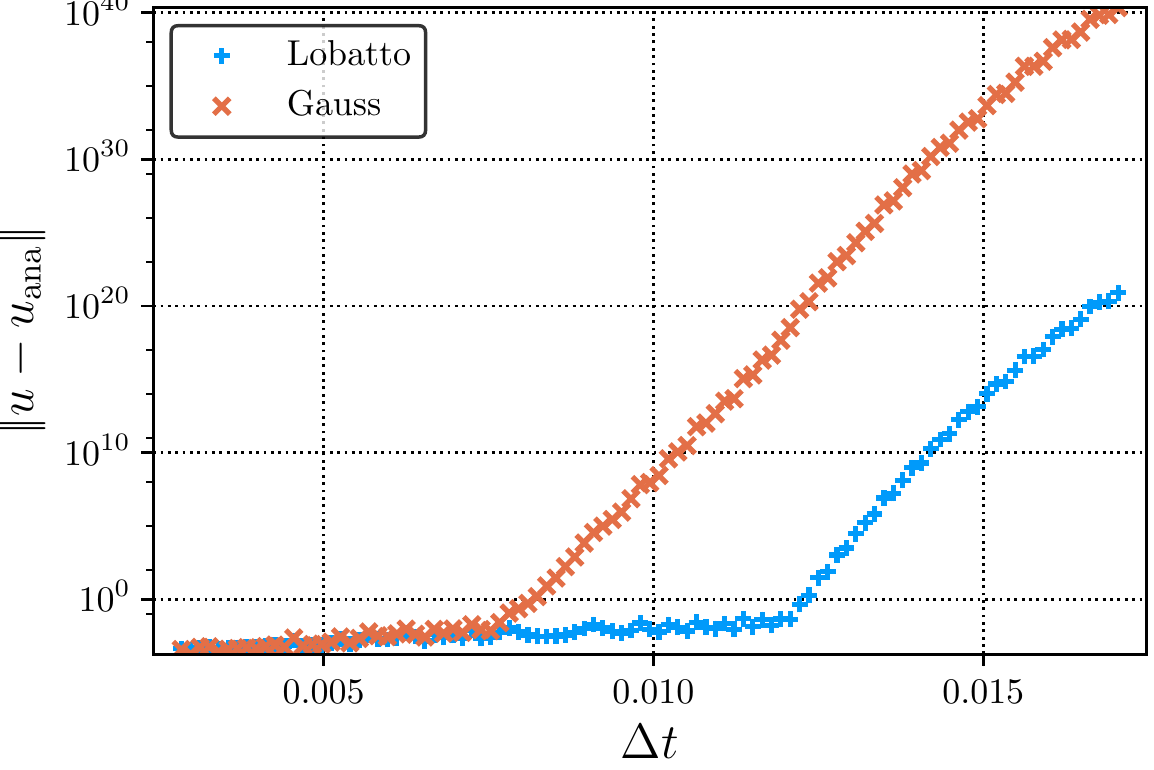}
    \caption{Upwind flux.}
    \label{fig:conservative_nonperiodic__1_p_cosh_x__sinpi_x__CFL_upwind}
  \end{subfigure}%
  \caption{Convergence results for \eqref{eq:var-lin-adv} using the speed 
           $a(x) = 1 + \cosh(x)$ in the domain $(-1,1)$ with initial condition 
           $u_0(x) = \sin(\pi x)$ and boundary data given by the analytical
           solution \eqref{eq:var-lin-adv-without-boundary-solution}.
           The semidiscretisation uses $N$ elements with polynomials of degree~$p$
           and the split form \eqref{eq:standard-general}.}
  \label{fig:conservative_nonperiodic__1_p_cosh_x__sinpi_x__CFL}
\end{figure}

In order to compute the maximal stable time step, several numerical experiments
with varying $\Delta t$ are conducted. A typical result using polynomials of
degree $p=5$ and $N=32$ elements for the split form \eqref{eq:standard-general}
is visualised in \autoref{fig:conservative_nonperiodic__1_p_cosh_x__sinpi_x__CFL}.
As can be seen there, the error is relatively small if $\Delta t$ is small enough.
However, the error starts to grow (exponentially) fast at a certain time step size.
This time step is considered to be the maximal stable $\Delta t$ in a CFL 
condition. It is approximately computed by considering the threshold $0.1$ for
the errors $\norm{u - u_\mathrm{ana}}$.

\begin{table}
\centering
  \caption{Experimental CFL conditions for \eqref{eq:var-lin-adv} using the speed 
           $a(x) = 1 + \cosh(x)$ in the domain $(-1,1)$ with initial condition 
           $u_0(x) = \sin(\pi x)$ and boundary data given by the analytical
           solution \eqref{eq:var-lin-adv-without-boundary-solution}.
          }
  \label{tab:cfl}
  \begin{tabular}{rrr | cccc}
    \toprule
    &&& \multicolumn{4}{c}{Maximal stable $\Delta t \cdot (2p+1) N$}
    \\
    &&& \multicolumn{2}{c}{Split form \eqref{eq:standard-general}}
      & \multicolumn{2}{c}{Unsplit form \eqref{eq:weighted-conservative-semidiscretisation}}
    \\
    & $p$ & $N$ & Central flux & Upwind flux & Central flux & Upwind flux
    \\
        \midrule 
    \multirow{20}{*}{\rotatebox{90}{Lobatto nodes (Lobatto for $a$)}} 
    & $ 3$ & $  32$  &      3.3  &      5.3  &      3.3  &      5.3 \\ 
    &      & $  64$  &      3.1  &      5.2  &      3.1  &      5.2 \\ 
    &      & $ 128$  &      3.1  &      5.0  &      3.1  &      5.0 \\ 
    &      & $ 256$  &      3.0  &      5.0  &      3.0  &      5.0 \\ 
    \cmidrule{3-7} 
    & $ 4$ & $  32$  &      2.6  &      4.6  &      2.6  &      4.6 \\ 
    &      & $  64$  &      2.6  &      4.5  &      2.6  &      4.5 \\ 
    &      & $ 128$  &      2.5  &      4.4  &      2.5  &      4.4 \\ 
    &      & $ 256$  &      2.5  &      4.3  &      2.5  &      4.3 \\ 
    \cmidrule{3-7} 
    & $ 5$ & $  32$  &      2.2  &      4.3  &      2.2  &      4.3 \\ 
    &      & $  64$  &      2.1  &      4.1  &      2.1  &      4.1 \\ 
    &      & $ 128$  &      2.1  &      4.0  &      2.1  &      4.0 \\ 
    &      & $ 256$  &      2.1  &      4.0  &      2.1  &      4.0 \\ 
    \cmidrule{3-7} 
    & $ 6$ & $  32$  &      1.9  &      4.0  &      1.9  &      4.0 \\ 
    &      & $  64$  &      1.8  &      3.9  &      1.8  &      3.9 \\ 
    &      & $ 128$  &      1.8  &      3.8  &      1.8  &      3.8 \\ 
    &      & $ 256$  &      1.8  &      3.7  &      1.8  &      3.7 \\ 
    \cmidrule{3-7} 
    & $ 7$ & $  32$  &      1.6  &      3.8  &      1.6  &      3.8 \\ 
    &      & $  64$  &      1.6  &      3.7  &      1.6  &      3.7 \\ 
    &      & $ 128$  &      1.6  &      3.6  &      1.6  &      3.6 \\ 
    &      & $ 256$  &      1.6  &      3.5  &      1.6  &      3.5 \\ 
    \midrule 
    \multirow{20}{*}{\rotatebox{90}{Gauss nodes (Lobatto for $a$)}} 
    & $ 3$ & $  32$  &      2.3  &      3.4  &      2.3  &      3.4 \\ 
    &      & $  64$  &      2.2  &      3.3  &      2.2  &      3.3 \\ 
    &      & $ 128$  &      2.1  &      3.1  &      2.1  &      3.1 \\ 
    &      & $ 256$  &      2.1  &      3.1  &      2.1  &      3.1 \\ 
    \cmidrule{3-7} 
    & $ 4$ & $  32$  &      2.0  &      3.0  &      2.0  &      3.0 \\ 
    &      & $  64$  &      1.9  &      3.0  &      1.9  &      3.0 \\ 
    &      & $ 128$  &      1.8  &      2.9  &      1.8  &      2.9 \\ 
    &      & $ 256$  &      1.8  &      2.8  &      1.8  &      2.7 \\ 
    \cmidrule{3-7} 
    & $ 5$ & $  32$  &      1.7  &      2.8  &      1.7  &      2.8 \\ 
    &      & $  64$  &      1.6  &      2.7  &      1.6  &      2.7 \\ 
    &      & $ 128$  &      1.6  &      2.6  &      1.6  &      2.6 \\ 
    &      & $ 256$  &      1.6  &      2.5  &      1.6  &      2.5 \\ 
    \cmidrule{3-7} 
    & $ 6$ & $  32$  &      1.5  &      2.5  &      1.5  &      2.5 \\ 
    &      & $  64$  &      1.4  &      2.5  &      1.4  &      2.5 \\ 
    &      & $ 128$  &      1.4  &      2.4  &      1.4  &      2.4 \\ 
    &      & $ 256$  &      1.4  &      2.3  &      1.4  &      2.3 \\ 
    \cmidrule{3-7} 
    & $ 7$ & $  32$  &      1.4  &      2.4  &      1.4  &      2.4 \\ 
    &      & $  64$  &      1.3  &      2.3  &      1.3  &      2.3 \\ 
    &      & $ 128$  &      1.3  &      2.3  &      1.3  &      2.3 \\ 
    &      & $ 256$  &      1.3  &      2.1  &      1.3  &      2.1 \\ 
    \bottomrule 
  \end{tabular}
\end{table}

\autoref{tab:cfl} displays the experimentally maximal stable time step sizes
$\Delta t$ multiplied with $(2p+1) N$ for both Lobatto and Gauss nodes using
varying polynomial degrees $p$, numbers of elements $N$, and numerical fluxes.
Both the split form \eqref{eq:standard-general} and the unsplit form
\eqref{eq:weighted-conservative-semidiscretisation} are used to compute the
numerical solutions.

As can be seen in \autoref{tab:cfl}, the CFL conditions are less stringent for
the upwind flux (adding additional dissipation) compared to the central one
(without additional dissipation). In most cases, there is a factor $\approx 2$
between the maximally stable time steps.
Moreover, there are no differences in the maximal stable $\Delta t$ between
the split form \eqref{eq:standard-general} and the unsplit form
\eqref{eq:weighted-conservative-semidiscretisation}.

If Gauss nodes are used, the discretisation of the speed~$a$ (either via Lobatto
nodes as described in Remark~\ref{rem:use-Lobatto-nodes-and-interpolate-to-Gauss}
or directly on Gauss nodes) does not influence the CFL conditions. Thus, the
results for Gauss nodes using Gauss nodes for $a$ are not shown in \autoref{tab:cfl}.

However, there is a clear difference between Lobatto nodes and Gauss nodes.
In general, the time step has to be be made smaller by a factor between approximately 
$0.6$ and $0.8$ if Gauss nodes are used. This factor is comparable to the
range of the spectrum of the operators considered in section~\ref{subsec:eigenvalues}.
Considering the dependency of the overall efficiency of Lobatto nodes and Gauss 
nodes on both the spatial and the temporal accuracy requirements, the implementation, 
the boundary conditions, the numerical fluxes etc., the different CFL conditions
could result in comparable efficiencies, see also \cite{gassner2011comparison}.

\section{Nonlinear Equations}
\label{sec:burgers}

In order to augment the investigations of the previous sections, a semidiscretisation
of a nonlinear conservation law will be analysed here. As in the previous sections,
the analysis will be performed at first in the continuous setting and thereafter
for the semidiscretisations using general SBP operators.

\subsection{Continuous Setting}

A nonlinear conservation law can be written as
\begin{equation}
\label{eq:hcl-nonlinear}
  \partial_t u + \partial_x f(u) = 0.
\end{equation}
Investigating the rate of change of a convex entropy $U$ in the physical element
$\Omega$, the conservation law \eqref{eq:hcl-nonlinear} is multiplied with the
entropy variables $w = U'(u)$ and integrated over $\Omega$, resulting for a
sufficiently smooth solution in
\begin{equation}
  \od{}{t} \int_\Omega U(u)
  =
  \int_{\Omega} w \cdot \partial_t u
  =
  -\int_{\Omega} w \cdot \partial_x f(u)
  =
  -\int_{\Omega} \underbrace{w \cdot f'(u)}_{=F'(u)} \cdot \partial_x u
  =
  - F(u) \big|_{\partial\Omega},
\end{equation}
where the entropy flux $F$ fulfilling $F'(u) = U'(u) \cdot f'(u)$ has been inserted,
cf. \cite{tadmor2003entropy}. With periodic boundary conditions or compactly
supported data, this results in the entropy equality $\od{}{t} \int_\Omega U(u) = 0$
and an entropy inequality $\od{}{t} \int_\Omega U(u) \leq 0$ will be used for
general solutions.

For general conservation laws, entropy conservative fluxes in the sense of
Tadmor \cite{tadmor1987numerical, tadmor2003entropy} can be used to construct
high-order semidiscretisations on periodic \cite{lefloch2002fully} or bounded
domains using diagonal norm SBP operators including the boundaries
\cite{fisher2013highJCP}. If the entropy conservative numerical fluxes can be
written as products of arithmetic mean values, the semidiscretisation corresponds
to a split form \cite{gassner2016split}. For such split forms, appropriate boundary
terms have been constructed for some specific examples of nonlinear balance laws
\cite{ranocha2016summation, ranocha2017extended, ranocha2017shallow, ortleb2016kinetic}.

Burgers' equation
\begin{equation}
  \partial_t u + \partial_x \frac{u^2}{2} = 0
\end{equation}
is a special example in this theory, since the skew-symmetric form
\begin{equation}
  \partial_t u + \frac{1}{3} u \partial_x u + \frac{1}{3} \partial_x u^2 = 0
\end{equation}
has been known for a long time to allow proofs of conservation and $L^2$ stability
using integration by parts \cite[equation (6.40)]{richtmyer1967difference}. Indeed,
the rate of change of the conserved quantity $u$ is obtained as
\begin{equation}
  \od{}{t} \int_\Omega u
  =
  \int_{\Omega} \partial_t u
  =
  - \frac{1}{3} \int_{\Omega} u \partial_x u
  - \frac{1}{3} \int_{\Omega} \partial_x u^2
  =
  - \frac{1}{2} u^2 \big|_{\partial\Omega}
  + 0,
\end{equation}
where integration by parts has been used for half of 
$\int_{\Omega} u \partial_x u$. Similarly, the rate of change of
the $L^2$ entropy $U(u) = \frac{1}{2} u^2$ is obtained by
\begin{equation}
  \od{}{t} \int_\Omega \frac{1}{2} u^2
  =
  \int_{\Omega} u \partial_t u
  =
  - \frac{1}{3} \int_{\Omega} u^2 \partial_x u
  - \frac{1}{3} \int_{\Omega} u \partial_x u^2
  =
  - \frac{1}{3} u^3 \big|_{\partial\Omega},
\end{equation}
using again integration by parts.

\subsection{Semidiscretisations Using General SBP Operators}

Based on the results of \cite{fisher2013discretely, gassner2013skew} for diagonal
norm SBP operators including the boundaries, a semidiscretisation for general SBP
operators has been proposed in \cite{ranocha2016summation, ranocha2017extended},
which can be written as
\begin{equation}
\label{eq:burgers-semidiscretisation}
  \partial_t \vec{u}
  + \frac{1}{3} \mat{u}[^*] \mat{D} \vec{u}
  + \frac{1}{3} \mat{D} \mat{u} \vec{u}
  + \mat{M}[^{-1}] \mat{R}[^T] \mat{B} \left(
    \vecfnum 
    - \frac{1}{3} \mat{R} \mat{u} \vec{u}
    - \frac{1}{6} \bigl( \mat{R} \vec{u} \bigr) \cdot \bigl( \mat{R} \vec{u} \bigr)
  \right)
  = 0.
\end{equation}
Here, $\mat{u} = \diag{\vec{u}}$ is a multiplication operator, $\mat{u}[^*]
= \mat{M}[^{-1}] \mat{u}[^T] \mat{M}$ is its $\mat{M}$-adjoint, and $\vecfnum = 
\bigl(\fnum_L, \fnum_R \bigr)^T$ is the vector of the numerical fluxes. Thus, 
this semidiscretisation is similar to the split form \eqref{eq:standard-general}.

Similarly to the continuous setting, the semidiscretisation
\eqref{eq:burgers-semidiscretisation} is conservative across elements, since
\begin{equation}
\begin{aligned}
  \vec{1}^T \mat{M} \partial_t \vec{u}
  &=
  - \frac{1}{3} \vec{1}^T \mat{M} \mat{u}[^*] \mat{D} \vec{u}
  - \frac{1}{3} \vec{1}^T \mat{M} \mat{D} \mat{u} \vec{u}
  - \vec{1}^T \mat{R}[^T] \mat{B} \left(
    \vecfnum 
    - \frac{1}{3} \mat{R} \mat{u} \vec{u}
    - \frac{1}{6} \bigl( \mat{R} \vec{u} \bigr) \cdot \bigl( \mat{R} \vec{u} \bigr)
  \right)
  \\
  &=
  - \frac{1}{3} \vec{u}^T \mat{M} \mat{D} \vec{u}
  - \vec{1}^T \mat{R}[^T] \mat{B} \left(
    \vecfnum
    - \frac{1}{6} \bigl( \mat{R} \vec{u} \bigr) \cdot \bigl( \mat{R} \vec{u} \bigr)
  \right)
  =
  \vec{0} - \vec{1}^T \mat{R}[^T] \mat{B} \vecfnum,
\end{aligned}
\end{equation}
due to the SBP property \eqref{eq:SBP}, cf. \cite{ranocha2017extended}.
Similarly, the semidiscretisation \eqref{eq:burgers-semidiscretisation} is stable
across elements if an entropy stable numerical flux $\fnum$ is used, since
\begin{equation}
\begin{aligned}
  \vec{u}^T \mat{M} \partial_t \vec{u}
  &=
  - \frac{1}{3} \vec{u}^T \mat{u}[^T] \mat{M} \mat{D} \vec{u}
  - \frac{1}{3} \vec{u}^T \mat{M} \mat{D} \mat{u} \vec{u}
  - \vec{u}^T \mat{R}[^T] \mat{B} \left(
    \vecfnum 
    - \frac{1}{3} \mat{R} \mat{u} \vec{u}
    - \frac{1}{6} \bigl( \mat{R} \vec{u} \bigr) \cdot \bigl( \mat{R} \vec{u} \bigr)
  \right)
  \\
  &=
  - \vec{u}^T \mat{R}[^T] \mat{B} \left(
    \vecfnum 
    - \frac{1}{6} \bigl( \mat{R} \vec{u} \bigr) \cdot \bigl( \mat{R} \vec{u} \bigr)
  \right),
\end{aligned}
\end{equation}
again due to the SBP property \eqref{eq:SBP}. Adding the contributions of two
neighbouring elements with indices $+$ and $-$ at the same boundary yields
\begin{equation}
  (u_+ - u_-) \fnum(u_-, u_+) - \frac{1}{6} \left( u_+^3 - u_-^3 \right) \leq 0,
\end{equation}
since the numerical flux $\fnum$ is entropy stable for the $L^2$ entropy $U(u)
= \frac{1}{2} u^2$ \cite{tadmor1987numerical, tadmor2003entropy}, cf.
\cite{ranocha2017extended}.

\begin{remark}
  In general, $L^2$ stability is not sufficient to obtain convergence for nonlinear
  conservation laws, since weak convergence and nonlinear operations are not
  compatible in general. Similarly, the semidiscretisations described above
  should be considered as entropy stable baseline schemes that have to be enhanced
  by additional techniques if discontinuities arise, such as artificial viscosity
  \cite{vonneumann1950method, tadmor1989convergence, grahs2002entropy, grahs2002image, 
  breuss2006numerical, bianchini2005vanishing, mattsson2004stable,
  nordstrom2006conservative, persson2006sub, barter2010shock, guermond2011suitable, 
  guermond2011entropy}, modal filtering \cite{grahs2002data, burgel2002continuous,
  hesthaven2008filtering, meister2012application, meister2012comparison,
  meister2013extended}, finite volume subcells \cite{huerta2012simple,
  dumbser2014posteriori, sonntag2014shock, meister2016positivity, sonntag2016efficient,
  ranocha2017shallow}, ENO type dissipation \cite{fjordholm2012arbitrarily,
  fjordholm2013high}, or comparison with WENO methods \cite{fisher2012high,
  fisher2013highJCP}.
\end{remark}

\begin{remark}
  Curvilinear coordinates have been investigated inter alia in \cite{kopriva2006metric,
  wintermeyer2017entropy, svard2004coordinate}. However, this setting is out of
  scope of this article, since the main target is the investigation of
  semidiscretisation using general SBP operators for linear equations with
  variable coefficients.
\end{remark}

\subsection{Numerical Results}

Here, the periodic initial value problem
\begin{equation}
\label{eq:burgers}
\begin{aligned}
  \partial_t u(t,x) + \partial_x \frac{u(t,x)^2}{2} &= 0,
  && t \in (0,T), x \in [0,2],
  \\
  u(0,x) &= u_0(x),
\end {aligned}
\end{equation}
with initial condition $u_0(x) = \sin(\pi x)$ will be solved up to $T = 0.3$.
The analytical solution can be computed by solving the implicit equation
$u(t,x) = u_0(x - t u(t,x))$ \cite[Example I.1.3]{lefloch2002hyperbolic}.

The domain $[0,2]$ is divided into $N$ elements of width $\frac{2}{N}$ and
polynomials of degree $p$ are used on each element, either in a nodal Lobatto
or a nodal Gauss basis. Godunov's flux \cite{osher1984riemann}
\begin{equation}
  \fnum(u_-,u_+)
  =
  \begin{cases}
    \displaystyle \min_{u \in [u_-,u_+]} \frac{u^2}{2}, &\text{if } u_- \leq u_+,
    \\
    \displaystyle \max_{u \in [u_+,u_-]} \frac{u^2}{2}, &\text{else},
  \end{cases}
\end{equation}
is used as numerical flux and the semidiscrete scheme is advanced in time by
the ten-stage, fourth order, strong stability preserving Runge-Kutta method
of \cite{ketcheson2008highly} with time step $\frac{2}{(2p+1)N}$.

The errors $\norm{u - u_\mathrm{ana}}_M$ of the numerical solutions $u_\mathrm{num}$
compared to the analytical solution $u_\mathrm{ana}$ are computed using the SBP
mass matrix $\mat{M}$. Together with the corresponding experimental order of
convergence (EOC), representative results are shown in \autoref{tab:convergence-burgers}.

\begin{table}[!ht]
\centering
  \caption{Convergence results for \eqref{eq:burgers} with initial condition 
           $u_0(x) = \sin(\pi x)$ and Godunov's flux.}
  \label{tab:convergence-burgers}
  \begin{tabular}{rr | cccc}
    \toprule
    && \multicolumn{2}{c}{Lobatto nodes}
      & \multicolumn{2}{c}{Gauss nodes}
    \\
    $p$  & $N$
    & $\norm{u - u_\mathrm{ana}}_M$ & EOC & $\norm{u - u_\mathrm{ana}}_M$ & EOC
    \\
    \midrule
    $ 2$ & $  100$ & 4.89e-03 &       & 3.16e-04 &       \\
         & $  200$ & 1.44e-03 &  1.77 & 8.45e-05 &  1.90 \\
         & $  400$ & 3.04e-04 &  2.24 & 2.20e-05 &  1.94 \\
         & $  800$ & 5.39e-05 &  2.49 & 2.78e-06 &  2.99 \\
         & $ 1600$ & 9.29e-06 &  2.54 & 5.30e-07 &  2.39 \\
         & $ 3200$ & 1.47e-06 &  2.66 & 8.77e-08 &  2.60 \\
    \midrule 
    $ 3$ & $  100$ & 1.08e-03 &       & 8.84e-05 &       \\
         & $  200$ & 1.48e-04 &  2.86 & 3.08e-05 &  1.52 \\
         & $  400$ & 2.08e-05 &  2.83 & 1.57e-06 &  4.29 \\
         & $  800$ & 4.84e-06 &  2.10 & 2.54e-07 &  2.63 \\
         & $ 1600$ & 5.75e-07 &  3.07 & 2.91e-08 &  3.13 \\
         & $ 3200$ & 5.90e-08 &  3.29 & 2.19e-09 &  3.73 \\
    \midrule 
    $ 4$ & $  100$ & 2.15e-04 &       & 7.55e-05 &       \\
         & $  200$ & 4.38e-05 &  2.29 & 7.47e-06 &  3.34 \\
         & $  400$ & 7.83e-06 &  2.48 & 1.81e-07 &  5.37 \\
         & $  800$ & 4.54e-07 &  4.11 & 2.52e-08 &  2.85 \\
         & $ 1600$ & 1.82e-08 &  4.64 & 8.43e-10 &  4.90 \\
         & $ 3200$ & 7.74e-10 &  4.55 & 3.64e-11 &  4.53 \\
    \midrule 
    $ 5$ & $  100$ & 8.42e-05 &       & 3.84e-05 &       \\
         & $  200$ & 1.93e-05 &  2.12 & 6.93e-07 &  5.79 \\
         & $  400$ & 8.98e-07 &  4.43 & 7.63e-08 &  3.18 \\
         & $  800$ & 1.88e-08 &  5.58 & 1.12e-09 &  6.09 \\
         & $ 1600$ & 1.14e-09 &  4.04 & 5.23e-11 &  4.42 \\
         & $ 3200$ & 3.17e-11 &  5.17 & 1.18e-12 &  5.47 \\
    \midrule 
  \end{tabular}
\end{table}

As can be seen there, all schemes converge numerically in the investigated range
of parameters. For polynomials of degree $p$, the experimental order of convergence
is mostly between $p-\frac{1}{2}$ and $p+1$ if the resolution is good enough,
both for Lobatto and Gauss nodes.

Since a periodic problem is solved, the total mass $\int u$ should remain constant.
In the numerical experiments using $64$ bit floating point numbers, the total mass
(computed via the SBP mass matrix) is conserved up to machine precision and
therefore not shown explicitly in a table.

\section{Summary and Conclusions}
\label{sec:summary}

Conservation laws with variable coefficients have been discussed. At the continuous
level, conservation and stability are important properties that should be mimicked
discretely. Using classical summation-by-parts operators with diagonal norms and
including the boundary nodes, these can be obtained in a straightforward way,
mimicking the manipulations for the continuous estimates exactly.

However, in the broad setting of generalised summation-by-parts operators, the
corresponding results are less clear. Translating the schemes of classical SBP
discretisations to the generalised ones, additional care has to be taken.
Otherwise, the resulting methods will not be conservative and stable, as discussed
recently \cite{nordstrom2017conservation, manzanero2017insights}.
Nevertheless, by constructing new correction terms for general SBP operators,
both conservative and stable discretisations can be created, as shown in 
sections~\ref{sec:standard-L2} and \ref{sec:weighted-L2} and the corresponding
numerical results in section~\ref{sec:numerical-results}.

Of course, there are still many open problems. Starting with Burgers' equation as
an example of a nonlinear conservation law, conservative and stable
semidiscretisations using generalised SBP operators have been constructed in
\cite{ranocha2016summation, ranocha2017extended}, cf. \autoref{sec:burgers}.
Moreover, generalised SBP operators have been applied to some nonlinear systems
of balance laws, resulting in conservative and kinetic energy preserving or 
entropy stable schemes \cite{ortleb2016kinetic, ranocha2017shallow}. There,
techniques and ideas similar to the ones presented here have been used.
However, the creation of adequate formulations can become complicated and it may
even be unknown whether suitable corrections exist \citep{ranocha2017comparison}.
Furthermore, studying the advantages and disadvantages of classical and generalised 
SBP operators is still worthwhile.

\appendix

\section{Translation Rules}
\label{sec:translation-rules}

Here, some rules to translate the notation used in this article into other
conventions are given.

\subsection{Finite Difference Notations}

In order to translate the notation of this article into the one used in the
finite difference framework of \cite{nordstrom2017conservation},
\autoref{tab:translation-FD} can be used. There, the indices $L,R$ have been
used for the left and right boundary instead of the indices $\alpha,\beta$ found
in \cite{nordstrom2017conservation}.

\begin{table}[!ht]
\caption{Translation rules for the finite difference setting.}
\label{tab:translation-FD}
\begin{center}
\begin{tabular}{lcc}
  \toprule
  & Notation of this article & Finite difference notation 
  \\ \midrule
  Numerical solution & $\vec{u}$ & $\mathbf{u}$
  \\
  Multiplication operator & $\mat{u}$, $\mat{a}$ & $U$, $A$
  \\
  Mass / norm matrix & $\mat{M}$ & $P$
  \\
  Derivative matrix & $\mat{D}$ & $D$
  \\
  Restriction matrix & $\mat{R}$ & $\vect{ \mathbf{t}_L^T \\ \mathbf{t}_R^T }$
  \\
  \multirow{2}{*}{SBP property}
    & \multirow{2}{*}{
    $\mat{M} \mat{D} + \mat{D}[^T] \mat{M} = \mat{R}[^T] \mat{B} \mat{R}$}
    & $D = \inv{P} Q$, $Q + Q^T = E$,
    \\&& $E = \mathbf{t}_R \mathbf{t}_R^T - \mathbf{t}_L \mathbf{t}_L^T$
  \\
  Correction matrix & $\mat{M}[^{-1}] \mat{R}[^T] \mat{B}$
                    & $P^{-1} \left( -\mathbf{t}_L, \mathbf{t}_R \right)$
  \\ \bottomrule
\end{tabular}
\end{center}
\end{table}

\subsection{Finite Element Notations}

In order to translate the notation of this article into the one used in the
discontinuous Galerkin spectral element (DGSEM) framework of
\cite{kopriva2014energy, manzanero2017insights}, \autoref{tab:translation-DGSEM}
can be used.

\begin{table}[!ht]
\caption{Translation rules for the DGSEM setting.}
\label{tab:translation-DGSEM}
\begin{center}
\begin{tabular}{lcc}
  \toprule
  & Notation of this article & DGSEM notation 
  \\ \midrule
  Polynomial degree  & $p$ & $N$
  \\
  Numerical functions 
    & $\vec{u}$, $\vec{a}$ 
    & $U = I^N(u)$, $\{\vec{\boldsymbol{A}}^{el}\}$
  \\
  Multiplication operator 
    & $\mat{a}$
    & $[\boldsymbol{A}^{el}] = \diag{ \{\vec{\boldsymbol{A}}^{el}\} }$
  \\
  Discrete product
    & $\mat{a} \vec{u}$
    & $I^N[A U]$, $[\boldsymbol{A}^{el}] \{\vec{\boldsymbol{U}}^{el}\}$
  \\
  Derivative matrix & $\mat{D}$ & $[\boldsymbol{D}]$
  \\
  Mass matrix & $\mat{M}$ & $[\boldsymbol{M}]$
  \\
  Discrete scalar product
    & $\vec{u}^T \mat{M} \vec{v}$
    & $\langle U^{el}, V^{el} \rangle_N$, 
      $\{\vec{\boldsymbol{U}}^{el}\}^T [\boldsymbol{M}] \{\vec{\boldsymbol{V}}^{el}\}$
  \\
  \multirow{2}{*}{SBP property}
  & $\vec{u}^T \mat{M} \mat{D} \vec{v} + \vec{u}^T \mat{D}[^T] \mat{M} \vec{v}$
  & \multirow{2}{*}{$\langle U, V_\xi \rangle_N + \langle U_\xi, V \rangle = U V \big|_{-1}^1$}
  \\
  & $= \vec{u}^T \mat{R}[^T] \mat{B} \mat{R} \vec{v}$
  \\
  Numerical flux & $\fnum$ & $F^*$
  \\ \bottomrule
\end{tabular}
\end{center}
\end{table}

\section{Comparison with Results of Nordström and Ruggiu}
\label{subsec:comparison-nordstrom-ruggiu}

The results of \autoref{sec:standard-L2} might seem to contradict the results of
\cite[Proposition 4.6 and 5.3]{nordstrom2017conservation} concerning generalised
SBP operators.
There, the authors investigated semidiscretisations of the conservative linear 
advection equation~\eqref{eq:var-lin-adv} using standard and generalised SBP 
operators in a finite difference framework. They proved that their semidiscretisations
are conservative and stable if classical SBP operators (having diagonal norms
and including the boundary nodes) are used. Contrary, using generalised SBP operators,
they proved that their semidiscretisations are in general not conservative and stable.

However, the results of \autoref{sec:standard-L2} concerning conservative and stable
semidiscretisations using generalised SBP operators do not contradict the ones
of \citep{nordstrom2017conservation}, since the semidiscretisation using
generalised SBP operators proposed in \cite[Equation (18)]{nordstrom2017conservation}
is different from equation \eqref{eq:standard-general}, the semidiscretisation
investigated in this article. Using a notation similar to \citet{nordstrom2017conservation},
their semidiscretisation is
\begin{equation}
\label{eq:nordstrom-ruggiu}
  \partial_t \mathbf{u}
  =
  - \frac{1}{2} \inv{P} \left( Q A + A Q \right) \mathbf{u}
  - \frac{1}{2} U \, D \mathbf{a}
  + \sigma_L^{I}  A \inv{P} \mathbf{t}_L \left(
      \mathbf{t}_L^T \mathbf{u} - g_L
    \right)
  + \sigma_L^{II} \inv{P} \mathbf{t}_L \left(
      \mathbf{t}_L^T A \mathbf{u} - a(\xmin) g_L
    \right),
\end{equation}
whereas \eqref{eq:standard-general} with the upwind numerical flux
$\fnum(u_-, u_+) = a_- \, u_-$ can be written using the translation rules of
\autoref{tab:translation-FD} as
\begin{equation}
\label{eq:standard-general-nordstrom}
\begin{aligned}
  \partial_t \mathbf{u}
  =&
  - \frac{1}{2} \inv{P} Q A \mathbf{u}
  - \frac{1}{2} \inv{P} A^T Q \mathbf{u}
  - \frac{1}{2} \inv{P} U^T Q \mathbf{a}
  \\&
  + \frac{1}{2} \inv{P} \mathbf{t}_L \left(
      2 \left(\mathbf{t}_L^T \mathbf{a}\right) g_L 
      - \mathbf{t}_L^T A \mathbf{u} 
      - \left(\mathbf{t}_L^T \mathbf{a}\right) 
        \left(\mathbf{t}_L^T \mathbf{u}\right) 
    \right)
  - \frac{1}{2} \inv{P} \mathbf{t}_R \left(
      \left(
      \mathbf{t}_R^T \mathbf{a}\right) \left(\mathbf{t}_R^T \mathbf{u}\right)
      - \mathbf{t}_R^T A \mathbf{u} 
    \right).
\end{aligned}
\end{equation}
Here, the upwind numerical fluxes $\fnum(u_-, u_+) = a_- \, u_-$ has been inserted via 
\begin{equation}
\label{eq:translate-upwind-numerical-flux-to-nordstrom}
  -\mat{M}[^{-1}] \mat{R}[^T] \mat{B} \vecfnum
  =
  -P^{-1} \left( -\mathbf{t}_L, \mathbf{t}_R \right) 
  \begin{pmatrix}
    (\mathbf{t}_L^T \mathbf{a}) g_L \\
    (\mathbf{t}_R^T \mathbf{a}) (\mathbf{t}_R^T \mathbf{u})
  \end{pmatrix}
  =
  P^{-1} \mathbf{t}_L (\mathbf{t}_L^T \mathbf{a}) g_L
  - P^{-1} \mathbf{t}_R (\mathbf{t}_R^T \mathbf{a}) (\mathbf{t}_R^T \mathbf{u}).
\end{equation}
In the case of diagonal norms considered in \cite{nordstrom2017conservation},
the multiplication operators are diagonal and self-adjoint, i.e.
$\inv{P} A^T = A \inv{P}$. Hence, \eqref{eq:standard-general-nordstrom} can be
rewritten as
\begin{equation}
\label{eq:standard-diagonal-nordstrom}
\begin{aligned}
  \partial_t \mathbf{u}
  =&
  - \frac{1}{2} \inv{P} \left( Q A + A Q \right) \mathbf{u}
  - \frac{1}{2} U \, D \mathbf{a}
  \\&
  + \frac{1}{2} \inv{P} \mathbf{t}_L \left(
      2 \left(\mathbf{t}_L^T \mathbf{a}\right) g_L 
      - \mathbf{t}_L^T A \mathbf{u} 
      - \left(\mathbf{t}_L^T \mathbf{a}\right) 
        \left(\mathbf{t}_L^T \mathbf{u}\right) 
    \right)
  - \frac{1}{2} \inv{P} \mathbf{t}_R \left(
      \left(
      \mathbf{t}_R^T \mathbf{a}\right) \left(\mathbf{t}_R^T \mathbf{u}\right)
      - \mathbf{t}_R^T A \mathbf{u} 
    \right).
\end{aligned}
\end{equation}
Thus, the volume terms are the same as in \eqref{eq:nordstrom-ruggiu} but the
boundary terms are different, resulting in different properties regarding
conservation and stability.

If diagonal norms are considered and the boundaries are included, the vectors
$\mathbf{t}_{L/R}$ used for the interpolations are $\mathbf{t}_{L} = (1, 0, \dots,
0)^T$ and $\mathbf{t}_{R} = (0, \dots, 0, 1)^T$. Thus, interpolation to the boundary
and multiplication commute, i.e. $\mathbf{t}_{L/R}^T A \mathbf{u} = 
\left(\mathbf{t}_{L/R}^T \mathbf{a} \right) \left(\mathbf{t}_{L/R}^T \mathbf{u} \right)$,
since $A = \diag{\mathbf{a}}$. Therefore, the semidiscretisation
\eqref{eq:standard-diagonal-nordstrom} can be simplified as
\begin{equation}
\label{eq:standard-diagonal-incl-boundaries-nordstrom}
  \partial_t \mathbf{u}
  =
  - \frac{1}{2} \inv{P} \left( Q A + A Q \right) \mathbf{u}
  - \frac{1}{2} U \, D \mathbf{a}
  + \inv{P} \mathbf{t}_L \left(
        \left(\mathbf{t}_L^T \mathbf{a}\right) g_L 
      - \mathbf{t}_L^T A \mathbf{u}
    \right)
\end{equation}
if diagonal norm SBP operators including the boundary nodes are used. This is
the translation of equation \eqref{eq:standard-simplified} to the notation of
\cite{nordstrom2017conservation} if the upwind numerical flux is used.

\subsection*{\texorpdfstring{$a$}{a}-Generalised SBP Operators}

\citet[Section 8]{nordstrom2017conservation} proposed $a$-generalised SBP operators
$P^{-1} \Theta$ approximating the anti-symmetric part of $\partial_x (a \, u)$
via $P^{-1} \Theta \mathbf{u} \approx \frac{1}{2} \left( \partial_x (a \, u)
+ a \, \partial_x u \right)$. 
Their new discretisation is \citep[Equation (37)]{nordstrom2017conservation}
\begin{equation}
\label{eq:a-ganeralised-op-nordstrom}
  \partial_t \mathbf{u}
  =
  - P^{-1} \Theta \mathbf{u}
  - U \, P^{-1} \Theta \mathbf{1}
  - a(x_L) P^{-1} \mathbf{t}_L \left( \mathbf{t}_L^T \mathbf{u} - g_L \right).
\end{equation}
However, there is an important difference between the $a$-generalised SBP
operators of \cite{nordstrom2017conservation} and the setting described here.
Indeed, the analytical value $a(\xmin)$ of the speed at the boundary is used in
\eqref{eq:a-ganeralised-op-nordstrom}, whereas interpolated speeds 
$\mathbf{t}_L^T \mathbf{a}$ are used in \eqref{eq:standard-diagonal-nordstrom}.
In view of the conditions of \autoref{thm:standard-general}, the values 
$a(x_{L/R})$ and $\mathbf{t}_{L/R}^T \mathbf{a}$ should be the same.
For nodal DG methods, this property can be achieved as described in
Remark~\ref{rem:use-Lobatto-nodes-and-interpolate-to-Gauss}.

Comparing equation \eqref{eq:standard-diagonal-nordstrom} with the ``possible 
remedy for generalised SBP operators'' \citep[Section 8]{nordstrom2017conservation}
results in the following translation rules concerning the $a$-generalised SBP
operator $\inv{P} \Theta$
\begin{equation}
\label{eq:translation-rules-a-generalised-SBP-op}
\begin{aligned}
  \Theta
  &\longleftrightarrow 
  \begin{cases}
    \frac{1}{2} \left( Q A + A Q \right)
    - \frac{1}{2} \left(\mathbf{t}_L^T \mathbf{a}\right) \mathbf{t}_L \mathbf{t}_L^T
    + \frac{1}{2} \mathbf{t}_L \mathbf{t}_L^T A
    \\
    + \frac{1}{2} \left(\mathbf{t}_R^T \mathbf{a}\right) \mathbf{t}_R \mathbf{t}_R^T
    - \frac{1}{2} \mathbf{t}_R \mathbf{t}_R^T A,
  \end{cases}
  \\
  - \inv{P} \Theta \mathbf{u}
  &\longleftrightarrow 
  \begin{cases}
  - \frac{1}{2} \inv{P} \left( Q A + A Q \right) \mathbf{u}
  \\
  + \frac{1}{2} \inv{P} \mathbf{t}_L \left(
      \left(\mathbf{t}_L^T \mathbf{a}\right) \left(\mathbf{t}_L^T \mathbf{u}\right)
      - \mathbf{t}_L^T A \mathbf{u} 
    \right)
  \\
  - \frac{1}{2} \inv{P} \mathbf{t}_R \left(
      \left(\mathbf{t}_R^T \mathbf{a}\right) \left(\mathbf{t}_R^T \mathbf{u}\right)
      - \mathbf{t}_R^T A \mathbf{u} 
    \right),
  \end{cases}
  \\
  -U \inv{P} \Theta \mathbf{1} 
  &\longleftrightarrow 
  -\frac{1}{2} U \, D \mathbf{a},
  \\
  - a(\xmin) \inv{P} \mathbf{t}_L \left( \mathbf{t}_L^T \mathbf{u} - g_L \right)
  &\longleftrightarrow 
  \inv{P} \mathbf{t}_L \left( 
    \left(\mathbf{t}_L^T \mathbf{a}\right) g_L 
    - \left(\mathbf{t}_L^T \mathbf{a}\right) \left(\mathbf{t}_L^T \mathbf{u}\right)
    \right).
\end{aligned}
\end{equation}
Using these translations, equation \eqref{eq:standard-diagonal-nordstrom} is the
same as equation \eqref{eq:a-ganeralised-op-nordstrom}, i.e. equation (37) of
\cite{nordstrom2017conservation}. Thus, the construction of $a$-generalised SBP
operators as done in \cite[Section 8.2]{nordstrom2017conservation} can be
simplified significantly in the framework presented here.

As described above, there is an important difference between the setting of 
$a$-generalised SBP operators in \cite{nordstrom2017conservation} and the one
described here. Instead of \citep[Property iii) of Definition 8.1]{nordstrom2017conservation}
\begin{equation}
  \Theta + \Theta^T 
  = 
  a(\xmax) \mathbf{t}_R \mathbf{t}_R^T - a(\xmin) \mathbf{t}_L \mathbf{t}_L^T,
\end{equation}
the operator $\Theta$ given in the translation rule 
\eqref{eq:translation-rules-a-generalised-SBP-op} satisfies
\begin{equation}
  \Theta + \Theta^T 
  = 
  (\mathbf{t}_R^T \mathbf{a}) \mathbf{t}_R \mathbf{t}_R^T 
  - (\mathbf{t}_L^T \mathbf{a}) \mathbf{t}_L \mathbf{t}_L^T.
\end{equation}

\begin{remark}
\label{rem:translation-rules-sign}
  There might seem to be an error in the translation rules, since one could expect
  that terms corresponding to the same boundary should have the same sign.
  However, the signs in \eqref{eq:translation-rules-a-generalised-SBP-op} are
  correct. Indeed, $P^{-1} \Theta \mathbf{u}$ is an approximation of
  $\frac{1}{2} \left( \partial_x (a \, u) + a \, \partial_x u \right)$, since
  the volume terms $\frac{1}{2} \inv{P} \left( Q A + A Q \right) \mathbf{u}$
  approximate $\frac{1}{2} \left( \partial_x (a \, u) + a \, \partial_x u \right)$
  and the boundary terms are consistent with zero  due to the different signs.
  Thus, the boundary terms can be seen as corrections to inexact multiplication,
  similarly to the splitting of the volume terms, see also \cite{ranocha2016summation,
  ranocha2017extended}. Indeed, considering the boundary term for the left-hand side,
  $\left(\mathbf{t}_L^T \mathbf{a}\right) \left(\mathbf{t}_L^T \mathbf{u}\right)
  - \mathbf{t}_L^T A \mathbf{u}$ is the difference of the product of the interpolations
  and the interpolation of the product.
\end{remark}

\section{Weighted \texorpdfstring{$L^2$}{L²} Estimates for a Nonconservative Equation}
\label{sec:weighted-nonconservative}

Similar to \autoref{sec:weighted-L2}, the weighted $L^2$ framework of 
\cite{manzanero2017insights} will be applied to semidiscretisations of a 
nonconservative equation using generalised SBP operators.

Consider the nonconservative linear advection equation
\begin{equation}
\label{eq:nonconservative}
\begin{aligned}
  \partial_t u(t,x) + a(x) \, \partial_x u(t,x) &= 0,
  \qquad&& t > 0, \; x \in (\xmin, \xmax),
  \\
  u(t,\xmin) &= g_L(t),
  \qquad&& t \geq 0,
  \\
  u(0,x) &= u_0(x),
  \qquad&& x \in [\xmin,\xmax],
\end{aligned}
\end{equation}
with variable speed $a(x) > 0$ and compatible initial and boundary conditions
$u_0$, $g_L$. Here, stability will be investigated in a weighted $L^2$ norm
given by the scalar product
\begin{equation}
\label{eq:L2-inva}
  \scp{u}{v}_{L^2_{\inv a}}
  =
  \int_\xmin^\xmax \frac{1}{a} \, u \, v.
\end{equation}

\subsection{Continuous Estimates}

Multiplying equation \eqref{eq:nonconservative} with $u \, \inv{a}$ and integrating
results due to integration-by-parts in
\begin{equation}
\label{eq:weighted-nonconservative-stability}
\begin{gathered}
  \od{}{t} \int_\xmin^\xmax \inv{a} \, u^2
  =
  2 \int_\xmin^\xmax \inv{a} \, u \, \partial_t u
  =
  - 2 \int_\xmin^\xmax u \, \partial_x u
  =
  - u^2 \big|_\xmin^\xmax,
  \\
  \implies
  \od{}{t} \norm{u}_{L^2_{\inv a}}^2 
  = g_L^2 - u(\xmax)^2 
  \leq g_L^2.
\end{gathered}
\end{equation}
Thus, the weighted $L^2$ norm fulfils
\begin{equation}
  \norm{u(t)}_{L^2_{\inv a}}^2
  \leq
  \norm{u_0}_{L^2_{\inv a}}^2 + \int_0^t g_L(\tau)^2 \dif \tau.
\end{equation}
As described in \cite[Section 4]{manzanero2017insights}, this can be translated
by the equivalence of norms
\begin{equation}
  \frac{\norm{u}_{L^2}^2}{\max_{x \in [\xmin,\xmax]} \set{a(x)}} 
  \leq
  \norm{u(t)}_{L^2_{\inv a}}^2
  \leq
  \frac{\norm{u}_{L^2}^2}{\min_{x \in [\xmin,\xmax]} \set{a(x)}} 
\end{equation}
to the following $L^2$ bound on the solution $u$
\begin{equation}
  \frac{\norm{u(t)}_{L^2}^2}{\max_{x \in [\xmin,\xmax]} \set{a(x)}} 
  \leq
  \frac{\norm{u_0}_{L^2}^2}{\min_{x \in [\xmin,\xmax]} \set{a(x)}} 
  + \int_0^t a(\xmin)^2 g_L(\tau)^2 \dif \tau.
\end{equation}

\subsection{Semidiscrete Estimates}
\label{subsec:weighted-nonconservative-semidiscrete}

In this section, general SBP semidiscretisations of 
are considered. Since no product rule has been used for the estimate
\eqref{eq:weighted-nonconservative-stability} in the continuous setting, the 
following (unsplit) form of the semidiscretisation will be considered, where
the $\mat{M}$-adjoint is still given as 
$\mat{a}[^*] = \mat{M}[^{-1}] \mat{a}[^T] \mat{M}$.
\begin{equation}
\label{eq:semidiscretisation-nonconservative-general}
  \partial_t \vec{u}
  =
  - \mat{a}[^*] \mat{D} \vec{u}
  - \mat{a}[^*] \mat{M}[^{-1}] \mat{R}[^T] \mat{B} \widetilde\vecfnum 
  + \mat{a}[^*] \mat{M}[^{-1}] \mat{R}[^T] \mat{B} \mat{R} \vec{u}.
\end{equation}
Since equations \eqref{eq:nonconservative} is not conservative, it is advantageous
to consider modified numerical fluxes not depending on the velocity and treat the
speed~$a$ separately. Thus, the modified fluxes considered in this section are
\begin{align}
  \label{eq:fnum-modified-central}
  \text{Modified central flux}\qquad
  \widetilde\fnum(u_-, u_+)
  &=
  \frac{u_- + u_+}{2},
  \\
  \label{eq:fnum-modified-upwind}
  \text{Modified upwind flux}\qquad
  \widetilde\fnum(u_-, u_+)
  &=
  u_-.
\end{align}
Assume that $\mat{a}[^{-T}] \mat{M}$ induces a norm, where 
$\mat{a}[^{-T}] := \left( \mat{a}[^T] \right)^{-1}$, cf. Remarks~\ref{rem:norms} 
and \ref{rem:weighted-conservative-compared-to-standard} for the analogous discussion
of $\mat{a} \mat{M}$. Then, the semidiscrete energy rate can be obtained via
multiplication of \eqref{eq:semidiscretisation-nonconservative-general} with 
$\vec{u}^T \mat{a}[^{-T}] \mat{M}$, i.e.
\begin{equation}
\label{eq:semidiscretisation-nonconservative-general-stability}
\begin{aligned}
  \od{}{t} \norm{\vec{u}}_{\mat{a}[^{-1}]\mat{M}}^2
  &=
  2 \vec{u}^T \mat{a}[^{-T}] \mat{M} \partial_t \vec{u}
  =
  - 2 \vec{u}^T \mat{M} \mat{D} \vec{u}
  - 2 \vec{u}^T \mat{R}[^T] \mat{B} \widetilde\vecfnum 
  + 2 \vec{u}^T \mat{R}[^T] \mat{B} \mat{R} \vec{u}
  \\
  &=
  \vec{u}^T \mat{R}[^T] \mat{B} \mat{R} \vec{u}
  - 2 \vec{u}^T \mat{R}[^T] \mat{B} \widetilde\vecfnum.
\end{aligned}
\end{equation}
Using the central flux~\eqref{eq:fnum-modified-central}, the contribution of one 
boundary to the rate of change of the energy becomes
\begin{equation}
  2 \left( u_+ - u_- \right) \widetilde\fnum
  - \left( u_+^2 - u_-^2 \right)
  =
  \left( u_+ - u_- \right) \left( u_+ + u_- \right)
  - \left( u_+^2 - u_-^2 \right)
  =
  0.
\end{equation}
Similarly, the modified upwind numerical flux \eqref{eq:fnum-modified-upwind}
can be used, resulting in the energy rate contribution
\begin{equation}
  2 \left( u_+ - u_- \right) \widetilde\fnum
  - \left( u_+^2 - u_-^2 \right)
  =
  2 \left( u_+ - u_- \right) u_-
  - \left( u_+^2 - u_-^2 \right)
  =
  - (u_- - u_+)^2
  \leq 0.
\end{equation}
Considering now the total rate of change of the energy in a bounded domain, using 
the modified upwind flux \eqref{eq:fnum-modified-upwind} at the exterior boundaries
results due to \eqref{eq:semidiscretisation-nonconservative-general-stability} in
\begin{equation}
\label{eq:semidiscretisation-nonconservative-general-stability-global}
\begin{aligned}
  \od{}{t} \norm{\vec{u}}_{\mat{a}[^{-1}]\mat{M}}^2
  &=
  \vec{u}^T \mat{R}[^T] \mat{B} \mat{R} \vec{u}
  - 2 \vec{u}^T \mat{R}[^T] \mat{B} \widetilde\vecfnum
  =
  u_R^2 - u_L^2 - 2 u_R^2 + 2 u_L g_L
  \\
  &=
  g_L^2 - u_R^2 - \left( u_L - g_L \right)^2
  \leq g_L^2 - u_R^2
  \leq g_L^2.
\end{aligned}
\end{equation}
This is again consistent with the global estimate
\eqref{eq:weighted-nonconservative-stability} with an additional stabilising term
$- (u_L - g_L)^2 \leq 0$.
This proves
\begin{theorem}
\label{thm:weighted-nonconservative}
  Using general SBP discretisations, the semidiscretisation
  \eqref{eq:semidiscretisation-nonconservative-general} of
  \eqref{eq:nonconservative} is both conservative and stable across elements if the
  upwind numerical flux \eqref{eq:fnum-modified-upwind} is used at the 
  exterior boundaries.
  
  If multiple elements are used, the numerical flux at inter-element boundaries
  can be chosen to be the upwind one (adding additional dissipation) or the central
  flux \eqref{eq:fnum-modified-central} (without additional dissipation).
\end{theorem}
In order to give reliable estimates, $\mat{a}[^{-1}] \mat{M}$ has to be symmetric 
and positive definite, see also Remarks~\ref{rem:norms} and 
\ref{rem:weighted-conservative-compared-to-standard}.

\begin{remark}
\label{rem:weighted-nonconservative-simplified}
  If boundary nodes are included, the semidiscretisation
  \eqref{eq:semidiscretisation-nonconservative-general}
  can be simplified as follows. The term 
  $- \mat{a}[^*] \mat{M}[^{-1}] \mat{R}[^T] \mat{B} \widetilde\vecfnum
  =- \mat{M}[^{-1}] \mat{a}[^T] \mat{R}[^T] \mat{B} \widetilde\vecfnum$
  can be paraphrased as ``If this term is integrated with a vector $\vec{v}$ via
  the mass matrix $\mat{M}$, the result is
  $\vec{v}^T \mat{a}[^T] \mat{R}[^T] \mat{B} \widetilde\vecfnum$, the difference
  of the product of $v$, $a$, and the modified numerical flux at the boundaries.''.
  Since boundary nodes are included, the computation of the product can be
  modified. Thus, the semidiscretisation can be written as
  \begin{equation}
  \label{eq:semidiscretisation-nonconservative-boundaries-included}
    \partial_t \vec{u}
    =
    - \mat{a}[^*] \mat{D} \vec{u}
    - \mat{M}[^{-1}] \mat{R}[^T] \mat{B} \vecfnum 
    + \mat{a}[^*] \mat{M}[^{-1}] \mat{R}[^T] \mat{B} \mat{R} \vec{u},
  \end{equation}
  where a standard numerical flux $\fnum$ (approximating $a u$) accounts for the 
  multiplication with $a$. Thus, the numerical fluxes of the previous sections
  can be used. Furthermore, if the norm matrix $\mat{M}$ is diagonal, the
  semidiscretisation \eqref{eq:semidiscretisation-nonconservative-general} becomes
  \begin{equation}
  \label{eq:semidiscretisation-nonconservative-simplified}
    \partial_t \vec{u}
    =
    - \mat{a} \mat{D} \vec{u}
    - \mat{M}[^{-1}] \mat{R}[^T] \mat{B} \left( \vecfnum - \mat{R} \mat{a} \vec{u} \right).
  \end{equation}
  Thus, all the stability results of \autoref{thm:weighted-nonconservative}
  can be transferred to \eqref{eq:semidiscretisation-nonconservative-simplified},
  if the modified numerical fluxes \eqref{eq:fnum-modified-central} and
  \eqref{eq:fnum-modified-upwind} are replaced by their counterparts
  \eqref{eq:fnum-edge-central}, \eqref{eq:fnum-split-central}, 
  \eqref{eq:fnum-unsplit-central} and
  \eqref{eq:fnum-edge-upwind}, \eqref{eq:fnum-split-upwind}, 
  \eqref{eq:fnum-unsplit-upwind}, respectively. The variants of the unmodified
  fluxes are again equivalent if boundary nodes are included.
\end{remark}

A comparison of these results with the ones of \cite{manzanero2017insights}
is given in \ref{subsec:comparison-manzanero-nonconservative}.

\section{Comparison with Results of Manzanero, Rubio, Ferrer, Valero, and Kopriva}
\label{sec:comparison-manzanero}

In this section, a comparison of the weighted estimates in \autoref{sec:weighted-L2} 
and \ref{sec:weighted-nonconservative} with some results of
\cite{manzanero2017insights} will be given.

\subsection{Conservative Equation}
\label{subsec:comparison-manzanero-conservative}

The results of \autoref{sec:weighted-L2} concerning general SBP 
discretisations without boundary nodes might seem to contradict the numerical 
results of \cite{manzanero2017insights}.
After investigating a nodal DG method using Lobatto nodes, i.e. a diagonal norm
SBP discretisation including the boundary nodes, they presented the eigenvalues
of a periodic problem using both Lobatto and Gauss nodes (not including the boundaries).
There, the discretisation using Lobatto nodes resulted in eigenvalues on the
imaginary axis \citep[Figure~1]{manzanero2017insights}, whereas the discretisation
using Gauss nodes yielded eigenvalues with positive real part
\citep[Figure~2]{manzanero2017insights}. However, they used the split central 
flux \eqref{eq:fnum-split-central} instead of the unsplit one 
\eqref{eq:fnum-unsplit-central}, resulting in a stable scheme.

Despite of the different numerical fluxes, the scheme of \cite{manzanero2017insights}
is equivalent to the semidiscretisation 
\eqref{eq:weighted-conservative-semidiscretisation}
investigated here. Indeed, their semidiscretisation of the conservative equation
\eqref{eq:conservative} in the standard element is given as
\citep[Equation (35) with $\Delta x = 2$, $\alpha = 1$, and $\theta = 0$,%
dropping the superscript $el$]{manzanero2017insights}
\begin{equation}
  \langle \Phi, U_t \rangle_N
  =
  - F^* \Phi \big|_{-1}^1
  + \langle \Phi_\xi, A U \rangle_N.
\end{equation}
Using the translation rules of \autoref{tab:translation-DGSEM}, this can be 
rewritten as
\begin{equation}
  \vec{\Phi}^T \mat{M} \partial_t \vec{u}
  =
  - \vec{\Phi}^T \mat{R}[^T] \mat{B} \vecfnum
  + \vec{\Phi}^T \mat{D}[^T] \mat{M} \mat{a} \vec{u}.
\end{equation}
Due to the SBP property \eqref{eq:SBP}, this can be reformulated as
\begin{equation}
  \vec{\Phi}^T \mat{M} \partial_t \vec{u}
  =
  - \vec{\Phi}^T \mat{R}[^T] \mat{B} \vecfnum
  + \vec{\Phi}^T \mat{R}[^T] \mat{B} \mat{R} \mat{a} \vec{u}
  - \vec{\Phi}^T \mat{M} \mat{D} \mat{a} \vec{u}.
\end{equation}
Since this equation should hold for all test functions $\Phi$, i.e. all vectors
$\vec{\Phi}$, it is equivalent to
\begin{equation}
  \partial_t \vec{u}
  =
  - \mat{M}[^{-1}] \mat{R}[^T] \mat{B} \vecfnum
  + \mat{M}[^{-1}] \mat{R}[^T] \mat{B} \mat{R} \mat{a} \vec{u}
  - \mat{D} \mat{a} \vec{u},
\end{equation}
which is the same as \eqref{eq:weighted-conservative-semidiscretisation}.

\subsection{Nonconservative Equation}
\label{subsec:comparison-manzanero-nonconservative}

As in \ref{subsec:comparison-manzanero-conservative}, the results of
\ref{sec:weighted-nonconservative} concerning general SBP discretisations 
without boundary nodes might seem to contradict the numerical results of 
\cite{manzanero2017insights}. Analogously to the conservative case, they investigated
a nodal DG method using Lobatto nodes for the nonconservative equation
\eqref{eq:nonconservative}.  Again, the discretisation using Lobatto nodes resulted 
in eigenvalues on the imaginary axis \citep[Figure~1]{manzanero2017insights},
whereas the discretisation using Gauss nodes yielded eigenvalues with positive 
real part \citep[Figure~2]{manzanero2017insights}. 
However, they used the split central flux \eqref{eq:fnum-split-central} instead 
of the unsplit one \eqref{eq:fnum-unsplit-central} and the semidiscretisation
\eqref{eq:semidiscretisation-nonconservative-boundaries-included} adapted to 
SBP discretisations including the boundary nodes instead of the general discretisation
\eqref{eq:semidiscretisation-nonconservative-general}. 
Indeed, their semidiscretisation of the nonconservative equation 
\eqref{eq:conservative} in the standard element is given as
\citep[Equation (35) with $\Delta x = 2$, $\alpha = 0$, and $\theta = 1$,%
dropping the superscript $el$]{manzanero2017insights}
\begin{equation}
  \langle \Phi, U_t \rangle_N
  =
  - F^* \Phi \big|_{-1}^1
  + \langle U, (I^N[A \Phi])_\xi \rangle_N.
\end{equation}
Using the translation rules of \autoref{tab:translation-DGSEM}, this can be 
rewritten as
\begin{equation}
  \vec{\Phi}^T \mat{M} \partial_t \vec{u}
  =
  - \vec{\Phi}^T \mat{R}[^T] \mat{B} \vecfnum
  + \vec{\Phi}^T \mat{a}[^T] \mat{D}[^T] \mat{M} \vec{u}.
\end{equation}
Due to the SBP property \eqref{eq:SBP}, this can be reformulated as
\begin{equation}
  \vec{\Phi}^T \mat{M} \partial_t \vec{u}
  =
  - \vec{\Phi}^T \mat{R}[^T] \mat{B} \vecfnum
  + \vec{\Phi}^T \mat{a}[^T] \mat{R}[^T] \mat{B} \mat{R} \vec{u}
  - \vec{\Phi}^T \mat{a}[^T] \mat{M} \mat{D} \vec{u}.
\end{equation}
Since this equation should hold for all test function $\Phi$, i.e. all vectors
$\vec{\Phi}$, it is equivalent to
\begin{equation}
  \partial_t \vec{u}
  =
  - \mat{M}[^{-1}] \mat{R}[^T] \mat{B} \vecfnum
  + \mat{M}[^{-1}] \mat{a}[^T] \mat{R}[^T] \mat{B} \mat{R} \vec{u}
  - \mat{M}[^{-1}] \mat{a}[^T] \mat{M} \mat{D} \vec{u},
\end{equation}
which is the same as \eqref{eq:semidiscretisation-nonconservative-boundaries-included}.

\subsection{Numerical Results}
\label{subsec:numerical-results-eigenvalues-manzanero-et-al}

Here, the setup for the numerical experiments of \cite{manzanero2017insights}
will be used. The advection speed is given by $a(x) = 1 + (1-x^2)^5$ on the domain
$[\xmin,\xmax] = [-1,1]$ equipped with periodic boundary conditions.

Applying the discontinuous Galerkin spectral element (DGSEM) semidiscretisation
\cite{kopriva2009implementing}
with either Gauss-Lobatto or Gauss nodes, i.e. diagonal norm SBP
discretisations including the boundary nodes (Lobatto) and not including boundary
nodes (Gauss), the resulting ordinary differential equations can be written
as $\od{}{t} u = L u$, where the linear operator $L$ is defined via the
semidiscretisation. Investigating its eigenvalues, stability results can be
deduced.

\begin{figure}[!ht]
  \begin{subfigure}{0.49\textwidth}
    \includegraphics[width=\textwidth]{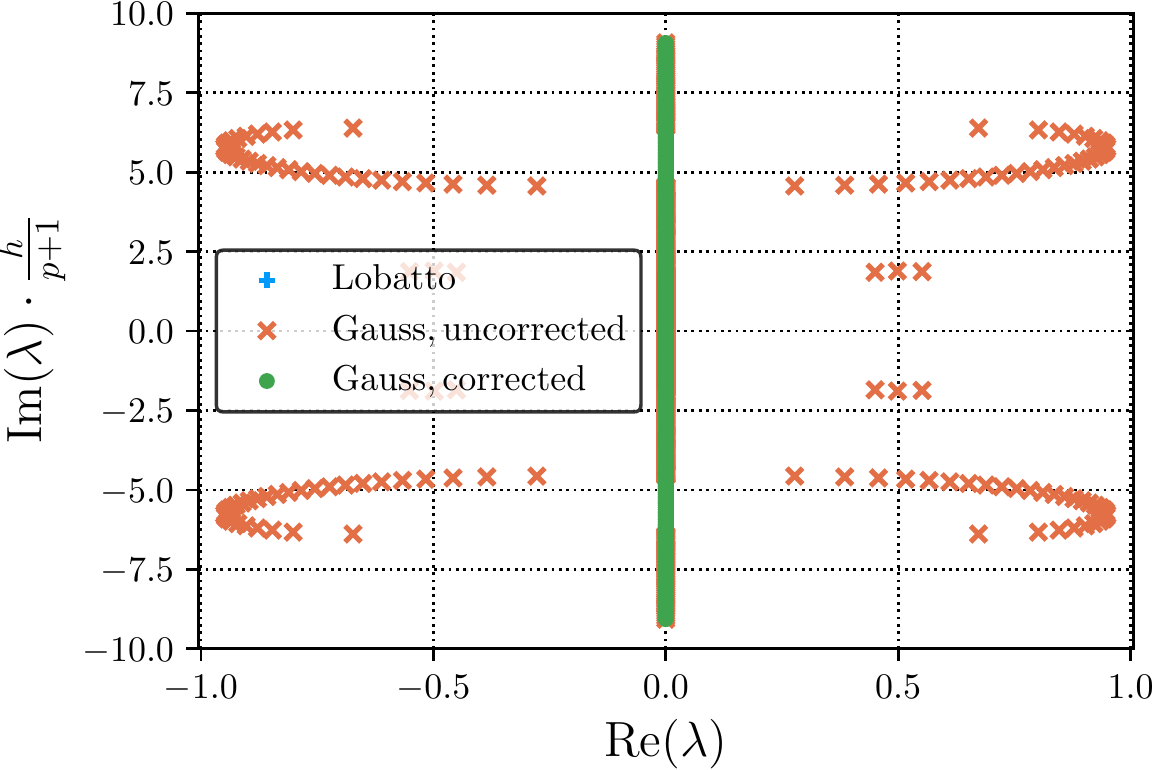}
    \caption{Conservative equation \eqref{eq:conservative}.}
    \label{fig:eigenvalues-conservative}
  \end{subfigure}%
  ~
  \begin{subfigure}{0.49\textwidth}
    \includegraphics[width=\textwidth]{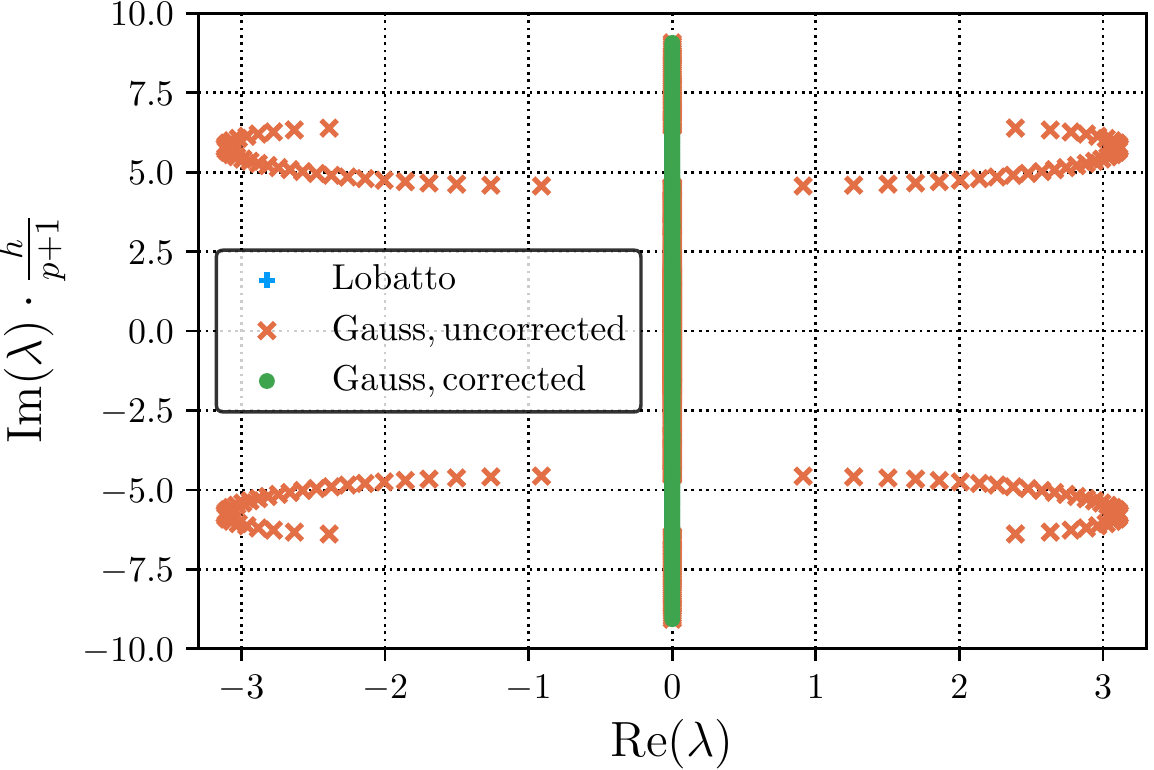}
    \caption{Nonconservative equation \eqref{eq:nonconservative}.}
    \label{fig:eigenvalues-nonconservative}
  \end{subfigure}
  \caption{Eigenvalues of the semidiscretisation using the central numerical fluxes,
           $200$ elements of width $h = \frac{2}{N}$ and polynomials of degree
           $p = 5$ for the linear advection equation with variable speed
           $a(x) = 1 + \left( 1 - x^2 \right)^5$ and periodic boundary conditions
           in $[-1,1]$.}
  \label{fig:eigenvalues}
\end{figure}

Using the DGSEM semidiscretisation \eqref{eq:weighted-conservative-semidiscretisation}
for the conservative equation \eqref{eq:conservative} with $200$ uniform elements
in $[-1,1]$ and polynomials of degree $\leq 5$, the eigenvalues of the operator
$L$ are given in Figure~\ref{fig:eigenvalues-conservative}. As can be seen there,
Lobatto nodes and the central numerical flux
\eqref{eq:fnum-split-central} result in a purely imaginary
spectrum. Using the same semidiscretisation and Gauss nodes results in some
eigenvalues in the half plane with positive real part. However, using the
corrected numerical flux \eqref{eq:fnum-unsplit-central}
results again in a purely imaginary spectrum.

Similarly, the spectrum of the DGSEM operators for the nonconservative equation 
\eqref{eq:nonconservative} with $200$ uniform elements in $[-1,1]$, polynomials 
of degree $\leq 5$, and different bases are shown in
\autoref{fig:eigenvalues-nonconservative}.
Again, the semidiscretisation \eqref{eq:semidiscretisation-nonconservative-simplified}
using Lobatto nodes and the central flux \eqref{eq:fnum-split-central} 
results in eigenvalues with vanishing real parts. Using Gauss nodes without 
modification results in eigenvalues with positive real parts. The corrected 
form \eqref{eq:semidiscretisation-nonconservative-general} with the modified
numerical flux \eqref{eq:fnum-modified-central} yields a purely imaginary spectrum.

\section*{Acknowledgements}
The author would like to thank the anonymous reviewers very much for their
helpful comments.

\printbibliography

\end{document}